\documentclass[graybox]{svmult}

\usepackage{mathptmx}       % selects Times Roman as basic font
\usepackage{mathptmx,amsmath,amsfonts,amssymb}
\usepackage{helvet}         % selects Helvetica as sans-serif font
\usepackage{courier}        % selects Courier as typewriter font
\usepackage{type1cm}        % activate if the above 3 fonts are
\usepackage{makeidx}         % allows index generation
\usepackage{graphicx}        % standard LaTeX graphics tool
\usepackage{multicol}        % used for the two-column index
\usepackage[bottom]{footmisc}% places footnotes at page bottom

 \newcommand{\ds}{\displaystyle }

\newcommand{\vc}[1]{\vec{#1}}

\newcommand{\sr}[1]{{\cal #1}}
 \newcommand{\dd}[1]{\mathbb{#1}}
 \newcommand{\edge}{{\rm r}}
 \newcommand{\cp}{{\Gamma}}
 \newcommand{\rmn}[1]{\if#11I\else {\if#12I\hspace{-0.12ex}I\hspace{-0.5ex}\else {\if #13I\hspace{-0.16ex}I\hspace{-0.16ex}I\hspace{-1.0ex}\else I\hspace{-1.2ex}V \fi} \fi} \fi}

\newcommand{\eqn}[1]{(\ref{eqn:#1})}
\newcommand{\lem}[1]{Lemma~\ref{lem:#1}}

\newcommand{\thr}[1]{Theorem~\ref{thr:#1}}

\newcommand{\rem}[1]{Remark~\ref{rem:#1}}
\newcommand{\fig}[1]{Figure~\ref{fig:#1}}
\newcommand{\tab}[1]{Table~\ref{tab:#1}}

\newcommand{\sectn}[1]{Section~\ref{sect:#1}}

\newcommand{\lemt}[1]{\ref{lem:#1}}

\newcommand{\thrt}[1]{\ref{thr:#1}}

\newcommand{\figt}[1]{\ref{fig:#1}}

\newcommand{\sect}[1]{\ref{sect:#1}}

\newcommand{\br}[1]{\langle #1 \rangle}

\newcommand{\ol}{\overline}
\newcommand{\ul}{\underline}
\newcommand{\pend}{\hfill \qed} %\thicklines \framebox(4.7,4.7)[l]{}}

\newenvironment{proof*}[1]{\noindent {\it  #1}. \rm}{\pend}

\newenvironment{mylist}[1]{\begin{list}{}
{\setlength{\itemindent}{#1mm}}
{\setlength{\itemsep}{0ex plus 0.2ex}}
{\setlength{\parsep}{0.5ex plus 0.2ex}}
{\setlength{\labelwidth}{10mm}}
}{\end{list}}

\makeindex             % used for the subject index
\begin{document}

\title*{Revisit to the tail asymptotics of the double QBD process: Refinement and complete solutions for the coordinate and diagonal directions}
\titlerunning{Revisit to tail asymptotics of double QBD}
\author{Masahiro Kobayashi and Masakiyo Miyazawa}
\authorrunning{M. Kobayashi and M. Miyazawa}
\institute{Masahiro Kobayashi \at Tokyo University of Science, Noda, Chiba 278-8510, Japan \email{masahilow@gmail.com}
\and Masakiyo Miyazawa \at Tokyo University of Science, Noda, Chiba 278-8510, Japan \email{miyazawa@is.noda.tus.ac.jp}}
\maketitle

\abstract*{We consider a two dimensional skip-free reflecting random walk on a nonnegative integer quadrant. We are interested in the tail asymptotics of its stationary distribution, provided its existence is assumed. We derive exact tail asymptotics for the stationary probabilities on the coordinate axis. This refines the asymptotic results in the literature, and completely solves the tail asymptotic problem on the stationary marginal distributions in the coordinate and diagonal directions. For this, we use the so-called analytic function method in such a way that either generating functions or moment generating functions are suitably chosen. The results are exemplified by a two node network with simultaneous arrivals.}

\abstract{We consider a two dimensional skip-free reflecting random walk on a nonnegative integer quadrant. We are interested in the tail asymptotics of its stationary distribution, provided its existence is assumed. We derive exact tail asymptotics for the stationary probabilities on the coordinate axis. This refines the asymptotic results in the literature, and completely solves the tail asymptotic problem on the stationary marginal distributions in the coordinate and diagonal directions. For this, we use the so-called analytic function method in such a way that either generating functions or moment generating functions are suitably chosen. The results are exemplified by a two node network with simultaneous arrivals.}

\section{Introduction} 
\label{sect:Introduction}

  We are concerned with a two dimensional reflecting random walk on a nonnegative integer quadrant, which is the set of two dimensional vectors $(i,j)$ such that $i,j$ are nonnegative integers. We assume that it is skip free in all directions, that is, its increments in each coordinate direction are at most one in absolute value. The boundary of the quadrant is partitioned into three faces, the origin and the two coordinate axes in the quadrant. We assume that the transition probabilities of this random walk is homogeneous on each boundary face, but they may change on different faces or the interior of the quadrant, that is, inside of the boundary.
  
  This reflecting random walk is referred to as a double quasi-birth-and-death process, a double QBD process for short, in \cite{Miyazawa 2009}. This process can be used to describe a two node queueing network under various setting such as server collaboration and simultaneous arrivals and departures, and its stationary distribution is important for the performance evaluation of such a network model. The existence of the stationary distribution, that is, stability, is nicely characterized, but the stationary distribution is hard to analytically get except for some special cases. Because of this as well as its own importance, research interest has been directed to its tail asymptotics.
  
  Until now, the tail asymptotics for the double QBD have been obtained in terms of its modelling primitives under the most general setting in Miyazawa \cite{Miyazawa 2009}, while less explicit results have been obtained for more general two dimensional reflecting random in Borovkov and Mogul'skii \cite{BM 2001}. Foley and McDonald \cite{FM 2005a,FM 2005b} studied the double QBD under some limitations. Recently, Kobayashi and Miyazawa \cite{KM 2011} modified the double QBD process in such a way that upward jumps may be unbounded, and studied its tail asymptotics. This process, called a double $M/G/1$ type, includes the double QBD process as a special case. For special cases such as tandem and priority queues, the tail asymptotics have been recently investigated in Guillemin and Leeuwaarden \cite{Guillemin-Leeuwaarden 2011} and Li and Zhao \cite{Li-Zhao 2009,Li-Zhao 2011a}. Recently, Li and Zhao \cite{Li-Zhao 2011b} challenged the general double QBD (see Additional note at the end of this section).
  
  The tail asymptotic problems have also been studied for a semi-martingale reflecting Brownian motion, SRBM for short, which is a continuous time and state counterpart of a reflecting random walk. For the two dimensional SRBM, the rate function for large deviations has been obtained under a certain extra assumption in Avram Dai and Hasenbein \cite{Avram-Dai-Hasenbein 2001}. Dai and Miyazawa \cite{Dai-Miyazawa 2011a} derived more complete answers but for the stationary marginal distributions. 
  
  Thus, we now have many studies on the tail asymptotics for two dimensional reflecting and related processes (see, e.g., \cite{Miyazawa 2011} for survey). Nevertheless, there still remain many problems unsolved even for the double QBD. The exact tail asymptotics of the stationary marginal distributions in the coordinate directions are one of such problems. Here, a sequence of nonnegative numbers $\{p(n); n=0,1,2\}$ is said to have exact tail asymptotic $\{h(n); n=0,1,\ldots\}$ if their ratio $p(n)/h(n)$ converges to a positive constant as $n$ goes to infinity. We also write this asymptotic as
\begin{eqnarray*}
  p(n) \sim h(n)
\end{eqnarray*}
  We will find $h(n) = n^{\kappa} a^{-n}$ or $n^{\kappa} (1+b(-1)^{n}) a^{-n}$ with constants $\kappa = -\frac 32, - \frac 12, 0, 1$, $a > 1$ and $|b| \le 1$ for the marginal distributions (also for the stationary probabilities on the boundaries). 

  We aim to completely solve the exact tail asymptotics of the stationary marginal distributions in the coordinate and diagonal directions, provided the stationary distribution exists. It is known that the tail asymptotics of the stationary probabilities on each coordinate axis are a key for them (e.g., see \cite{Miyazawa 2011}). These asymptotics have been studied in \cite{BM 2001,Miyazawa 2009}. They used Markov additive processes generated by removing one of the boundary face which is not the origin, and related their asymptotics. However, there are some limitations in that approach.
  
  In this paper, we revisit the double QBD process using a different approach, recently developed in \cite{Dai-Miyazawa 2011a,KM 2011,MR 2009}. This approach is purely analytic, and called an analytic function method. It is closely related to the kernel method used in \cite{Guillemin-Leeuwaarden 2011,Li-Zhao 2009,Li-Zhao 2011a}. Their details and related topics are reviewed in \cite{Miyazawa 2011}.
  
  The analytic function method in \cite{Dai-Miyazawa 2011a,KM 2011,MR 2009} only uses moment generating functions because they have nice analytic properties including convexity. However, a generating function is more convenient for a distribution on integers because they are polynomials. Thus, generating functions have been used in the kernel method.
  
  In this paper, we use both of generating functions and moment generating functions. We first consider the convergence domain of the moment generating function of the stationary distribution, which is two dimensional. This part mainly refers to recent results due to \cite{KM 2011}. Once the domain is obtained, we switch from moment generating function to generating function, and consider analytic behaviors around its dominant singular points. A key is the so called kernel function. We derive inequalities for it (see \lem{extension 1}), adapting the idea used in \cite{Dai-Miyazawa 2011a}. This is a crucial step in the present approach, which enables us to apply analytic extensions not using the Riemann surface which has been typically used in the kernel method. We then apply the inversion technique for generating functions, and derive the exact tail asymptotics of the stationary tail probabilities on the coordinate axes.
  
  The asymptotic results are exemplified by a two node queueing network with simultaneous arrivals. This model is an extension of a two parallel queues with simultaneous arrivals. For the latter, the tail asymptotics of its stationary distribution in the coordinate directions are obtained in \cite{Flatto-Hahn 1984,Flatto-McKean 1977}. We modify this model in such a way that a customer who has completed service may be routed to another queue with a given probability. Thus, our model is more like a Jackson network, but it does not have a product form stationary distribution because of simultaneous arrivals. We will discuss how we can see the tail asymptotics from the modeling primitives.
  
  This paper is made up by seven sections. In \sectn{Double QBD}, we introduce the double QBD process, and summarize existence results using moment generating functions. \sectn{Analytic function} considers the generating functions for the stationary probabilities on the coordinate axes. Analytic behaviors around their dominant singular points are studied. We then apply the inversion technique and derive exact asymptotics in Sections \sect{exact asymptotic v-a} and \sect{exact asymptotic not v-a}. The example for simultaneous arrivals is considered in \sectn{Application to}. We discuss some remaining problems in \sectn{Concluding remarks}.

\medskip
  
  \noindent ({\it Additional note}) After the first submission of this paper, we have known that Li and Zhao \cite{Li-Zhao 2011b} studied the same exact tail asymptotic problem, including the case that the tail asymptotics is periodic. This periodic case was lacked in our original submission, and was added in the present paper. Thus, we benefited by them. However, our approach is different from theirs although both uses analytic functions and its asymptotic inversions. Namely, the crucial step in \cite{Li-Zhao 2011b} is analytic extensions on a Riemann surface studied in \cite{FIM 1999}, while we use the convergence domain obtained in \cite{KM 2011} and the key lemma. Another difference is sorting tail asymptotic results. Their presentation is purely analytic while we use the geometrical classifications of \cite{KM 2011,Miyazawa 2009} (see also \cite{Miyazawa 2011}).

\section{Double QBD process and convergence domain} 
\label{sect:Double QBD}

The double QBD process was introduced and studied in \cite{Miyazawa 2009}. We here briefly introduce it, and present results on the tail asymptotics of its stationary distribution. We will use following set of numbers.
\begin{eqnarray*}
 && \dd{Z} = \mbox{ the set of all integers},\qquad \qquad \dd{Z}_{+} = \{ j \in \dd{Z}; j \ge 0 \},\\
 && \dd{U} = \{(i,j) \in \dd{Z}^{2}; i, j =0,1,-1\},\\
 && \dd{R} = \mbox{ the set of all real numbers}, \qquad
 \dd{R}_{+} = \{ x \in \dd{R}; x \ge 0 \},\\
 && \dd{C} = \mbox{ the set of all complex numbers}.
\end{eqnarray*}
  Let $S = \dd{Z}_{+}^{2}$, which is a state space for the double QBD process. Define the boundary faces of $S$ as
\begin{eqnarray*}
 \partial S_{0} = \{ (0,0) \}, \quad
 \partial S_{1} = \{(i,0) \in \dd{Z}^{2}_{+}; i \geq 1 \},\quad
 \partial S_{2} = \{(0,i) \in \dd{Z}^{2}_{+}; i \geq 1 \}.
\end{eqnarray*}
  Let $\partial S = \cup_{i=0}^{2} \partial S_{i}$ and $S_{+} = S \setminus \partial S$. We refer to $\partial S$ and $S_{+}$ as the boundary and interior of $S$, respectively.
  
  Let $\{\vc{Y}_{\ell}; \ell=0,1,\ldots\}$ be a skip free random walk on $\dd{Z}^{2}$. That is, its increments $\vc{X}^{(+)}_{\ell} \equiv \vc{Y}_{\ell} - \vc{Y}_{\ell-1}$ take values in $\dd{U}$, and are independent and identically distributed. By $\vc{X}^{(+)}$, we simply denote a random vector which has the same distribution as $\vc{X}^{(+)}_{\ell}$. Define a discrete time Markov chain $\{\vc{L}_{\ell}\}$ with state space $S$ by the transition probabilities:
\begin{eqnarray*}
P(\vc{L}_{\ell+1} = \vc{j} |\vc{L}_{\ell} = \vc{i}) = \left\{ 
\begin{array}{ll}
P(\vc{X}^{(+)} = \vc{j}-\vc{i}),  & \vc{j} \in S, \vc{i} \in S_{+},\\
P(\vc{X}^{(k)} = \vc{j}-\vc{i}), & \vc{j} \in S, \vc{i} \in \partial S_{k}, k=0,1,2,
\end{array}
\right.
\end{eqnarray*}
where $\vc{X}^{(k)}$ is a random vector taking values in $\{ (i_{1},i_{2}) \in \dd{U}; i_{3-k} \ge 0\}$ for $k=1,2$ and in $\{ (i_{1},i_{2}) \in \dd{U}; i_{1}, i_{2} \ge 0\}$ for $k=0$. Hence, we can write as
\begin{eqnarray}
\label{eqn:transition dynamics}
  \vc{L}_{\ell+1} = \vc{L}_{\ell} + \sum_{k=0,1,2,+} \vc{X}^{(k)}_{\ell} 1(\vc{L}_{\ell} \in S_{k}), \qquad \ell =0,1,2,\ldots.
\end{eqnarray}
  where $1(\cdot)$ is the indicator function of the statement ``$\cdot$'', and $\vc{X}^{(k)}_{\ell}$ has the same distribution as that of $\vc{X}^{(k)}$ for each $k = 0,1,2,+$, and is independent everything else.

  Thus, $\{\vc{L}_{\ell}\}$ is a skip free reflecting random walk on the nonnegative integer quadrant $S$, which is called a double QBD process because its QBD transition structure is unchanged when level and background state are exchanged.
  
  We denote the moment generating functions of $\vc{X}^{(k)}$ by $\gamma_{k}$, that is, for $\vc{\theta} \equiv (\theta_{1},\theta_{2}) \in \dd{R}^{2}$,
\begin{eqnarray*}
 \gamma_{k}(\vc{\theta}) = E(e^{\br{\vc{\theta}, \vc{X}^{(k)}}}), \qquad k=0,1,2, +,
\end{eqnarray*} 
where $\br{\vc{a}, \vc{b}} = a_{1} b_{1} + a_{2} b_{2}$ for $\vc{a} = (a_{1}, a_{2})$ and $\vc{b} = (b_{1}, b_{2})$. As usual, $\dd{R}^{2}$ is considered to be a metric space with Euclidean norm $\| \vc{a} \| \equiv \sqrt{\br{\vc{a}, \vc{a}}}$. In particular, a vector $\vc{c}$ is called a directional vector if $\|\vc{c} \|=1$. In this paper, we assume that
\begin{mylist}{0}
\item [(i)] The random walk $\{\vc{Y}_{\ell}\}$ is irreducible.
\item [(ii)] The reflecting process $\{\vc{L}_{\ell}\}$ is irreducible and aperiodic.
\item [(iii)] Either $E(X^{(+)}_{1}) \ne 0$ or $E(X^{(+)}_{2}) \ne 0$ for $\vc{X}^{(+)} = (X^{(+)}_{1}, X^{(+)}_{2})$.
\end{mylist}

\begin{remark} {\rm
\label{rem:condition (iii)}
  If $E(X^{(+)}_{1}) = E(X^{(+)}_{2}) = 0$, then it is known that the stationary distribution of $\{\vc{L}_{\ell}\}$ cannot have a light tail, that is, cannot geometrically (or exponentially) decay in all directions . See \cite{FMM 1995} and Remark 3.1 of \cite{KM 2011}. Thus, (iii) is not a restrictive assumption for considering the light tail. 
}\end{remark}

  Under these assumptions, tractable conditions are obtained for the existence of the stationary distribution in the book \cite{FMM 1995}. They are recently corrected in \cite{KM 2011}. We refer to this corrected version below.
\begin{lemma}[Lemma 2.1 of \cite{KM 2011}] {\rm
\label{lem:stability d=2}
  Assume condition (i)--(iii), and let
\begin{eqnarray*}
  && \vc{m} = (E(X^{(+)}_{1}), E(X^{(+)}_{2})),\\
  && \vc{m}^{(1)}_{\bot} = (E(X^{(1)}_{2}), - E(X^{(1)}_{1})),\\
  && \vc{m}^{(2)}_{\bot} = (-E(X^{(2)}_{2}), E(X^{(2)}_{1})).
\end{eqnarray*}
  Then, the reflecting random walk $\{\vc{L}_{\ell}\}$ has the stationary distribution if and only if either one of the following three conditions hold (see \cite{KM 2011}).
\begin{eqnarray}
\label{eqn:stability 1}
 && m_{1} < 0, m_{2} < 0, \br{\vc{m}, \vc{m}^{(1)}_{\bot}} < 0, \br{\vc{m}, \vc{m}^{(2)}_{\bot}} < 0,\\
\label{eqn:stability 2}
 && m_{1} \ge 0, m_{2} < 0, \br{\vc{m}, \vc{m}^{(1)}_{\bot}} < 0. \mbox{ In addition, $m^{(2)}_{2} < 0$ is needed if $m^{(2)}_{1} = 0$.} \quad \\
\label{eqn:stability 3}
 && m_{1} < 0, m_{2} \ge 0, \br{\vc{m}, \vc{m}^{(2)}_{\bot}} < 0, \mbox{ In addition, $m^{(1)}_{1} < 0$ is needed if $m^{(1)}_{2} = 0$.} \quad 
\end{eqnarray}
}\end{lemma}

  Throughout the paper, we also assume this stability condition. That is,
\begin{mylist}{0}
\item [(iv)] Either one of \eqn{stability 1}, \eqn{stability 2} or \eqn{stability 3} holds.
\end{mylist}

  In addition to the conditions (i)--(iv), we will use the following conditions to distinguish some periodical nature of the tail asymptotics.

\begin{mylist}{0}
\item [(v-a)] $P(\vc{X}^{(+)} \in \{(1,1),(-1,1),(0,0),(1,-1),(-1,-1)\}) < 1$.
\item [(v-b)] $P(\vc{X}^{(1)} \in \{(1,1), (0,0), (-1,1)\}) < 1$.
\item [(v-c)] $P(\vc{X}^{(2)} \in \{(1,1), (0,0), (1,-1)\}) < 1$ .
\end{mylist}
These conditions are said to be non-arithmetic in the interior and boundary faces $1, 2$, respectively, while the conditions that they do not hold are called arithmetic. The remark below explains why they are so called.

\begin{remark} {\rm
\label{rem:condition (v-b)}
  To see the meaning of these conditions, let us consider random walk $\{\vc{Y}_{\ell}\}$ on $\dd{Z}^{2}$. We can view this random walk as a Markov additive process in the $k$-th coordinate direction if we consider the $k$-th entry of $\vc{Y}_{\ell}$ as an additive component and the other entry as a background state ($k=1,2$). Then, the condition (v-a) is exactly the non-arithmetic condition of this Markov additive process in each coordinate direction (see \cite{MZ 2004} for the definition of the period of a Markov additive process). It is notable that, for the random walk $\{\vc{Y}_{\ell}\}$, if the Markov additive process in one direction is non-arithmetic, then the one in the other direction is also non-arithmetic.
  
  We can give similar interpretations for (v-b) and (v-c). Namely, for each $k=1,2$, consider the random walk with increments subject to the same distribution as $\vc{X}^{(k)}$. This random walk is also viewed as a Markov additive process with an additive component in the $k$-th coordinate direction. Then, (v-b) and (v-c) are the non-arithmetic condition of this Markov additive process for $k=1,2$, respectively.
}\end{remark}
  
\begin{remark} {\rm
\label{rem:X-shaped}
  These conditions were recently studied in \cite{Li-Zhao 2011b}. They called a probability distribution on $\dd{U} \equiv \{(i,j); i,j=-1,0,1\}$ to be $X$-shaped if its support is included in 
  $$\{(1,1),(-1,1),(0,0),(1,-1),(-1,-1)\}.$$ 
  Thus, the conditions (v-a), (v-b) and (v-c) are for $\vc{X}^{(+)}$, $\vc{X}^{(1)}$ and $\vc{X}^{(2)}$, respectively, not to be $X$-shaped.
}\end{remark}

  We denote the stationary distribution of $\{\vc{L}_{\ell}; \ell=0,1,\ldots\}$ by $\nu$, and let $\vc{L}$ be a random vector subject to $\nu$. Then, it follows from \eqn{transition dynamics} that
\begin{eqnarray}
\label{eqn:stationary dynamics}
  \vc{L} \simeq \vc{L} + \sum_{k=0,1,2,+} \vc{X}^{(k)} 1(\vc{L} \in S_{k}) ,
\end{eqnarray}
  where ``$\simeq$'' stands for the equality in distribution. We introduce four moment generating functions concerning $\nu$. For $\vc{\theta} \in \dd{R}^{2}$,
\begin{eqnarray*}
 && \varphi(\vc{\theta}) = E(e^{\br{\vc{\theta}, \vc{L}}}), \\
 && \varphi_{+}(\vc{\theta}) = E(e^{\br{\vc{\theta}, \vc{L}}} 1(\vc{L} \in S_{+})),\\
  && \varphi_{k}(\theta_{k}) = E(e^{\theta_{k} L_{k}} 1(\vc{L} \in S_{k})), \qquad k=1,2.
\end{eqnarray*}
  Then, from \eqn{stationary dynamics} and the fact that
\begin{eqnarray*}
  \varphi(\vc{\theta}) = \varphi_{+}(\vc{\theta}) + \sum_{k=1}^{2} \varphi_{k}(\theta_{k}) + \nu(\vc{0}),
\end{eqnarray*}
we can easily derive the following stationary equation.
\begin{eqnarray}
\label{eqn:stationary equation 1}
  \lefteqn{(1 - \gamma_{+}(\vc{\theta})) \varphi_{+}(\vc{\theta}) + (1 - \gamma_{1}(\vc{\theta})) \varphi_{1}(\theta_{1}) } \nonumber \hspace{10ex}\\
  && + (1 - \gamma_{2}(\vc{\theta})) \varphi_{2}(\theta_{2}) + (1 - \gamma_{0}(\vc{\theta})) \nu(\vc{0}) = 0, \qquad
\end{eqnarray}
  as long as $\varphi(\vc{\theta})$ is finite. Clearly, this finiteness holds for $\vc{\theta} \le \vc{0}$.
  
  To find the maximal region for \eqn{stationary equation 1} to be valid, we define the convergence domain of $\varphi$ as
\begin{eqnarray*}
  \sr{D} = \mbox{the interior of } \{ \vc{\theta} \in \dd{R}^{2}; \varphi(\vc{\theta}) < \infty \}.
\end{eqnarray*}
   This domain is obtained by Kobayashi and Miyazawa \cite{KM 2011}. To present this result, we introduce notations. 
  
  From \eqn{stationary equation 1}, we can see that the curves $1 - \gamma_{k}(\vc{\theta}) = 0$ for $k=+,1,2$ are keys for $\varphi(\vc{\theta})$ to be finite. Thus, we let
\begin{eqnarray*}
 && \Gamma_{k} = \{ \vc{\theta} \in \dd{R}^{2}; \gamma_{k}(\vc{\theta}) < 1\},\\
 && \partial \Gamma_{k} = \{ \vc{\theta} \in \dd{R}^{2}; \gamma_{k}(\vc{\theta}) = 1\}, \qquad k=1,2,+.
\end{eqnarray*}
  We denote the closure of $\Gamma_{k}$ by $\ol{\Gamma}_{k}$. Since $\gamma_{k}$ is a convex function, $\Gamma_{k}$ and $\ol{\Gamma}_{k}$ are convex sets. Furthermore, condition (i) implies that $\Gamma_{+}$ is bounded, that is, it is included in a ball in $\dd{R}^{2}$. Let
\begin{eqnarray*}
 && \vc{\theta}^{(k,\edge)} = \arg_{\vc{\theta}\in \dd{R}^{2}} \sup \{ \theta_{k}; \vc{\theta} \in \Gamma_{+} \cap \Gamma_{k} \}, \quad k=1,2, \\
 && \vc{\theta}^{(k,\min)} = \arg_{\vc{\theta}\in \dd{R}^{2}} \inf \{ \theta_{k}; \vc{\theta} \in \Gamma_{+}\},\\
 && \vc{\theta}^{(k,\max)} = \arg_{\vc{\theta}\in \dd{R}^{2}} \sup \{ \theta_{k}; \vc{\theta} \in \Gamma_{+}\}.
\end{eqnarray*}
  These extreme points play key roles in obtaining the convergence domain. It is notable that $\vc{\theta}^{(k,\edge)}$ is not the zero vector $\vc{0}$ because the stability condition (iv) implies that, for each $k=1,2$, $\Gamma_{+} \cap \Gamma_{k}$ contains $\vc{\theta} = (\theta_{1}, \theta_{2})$ such that $\theta_{k} > 0$ (see Lemma 2.2 of \cite{KM 2011}).
  
   We further need the following points.
\begin{eqnarray*}
  \vc{\theta}^{(k,\cp)} = \left\{\begin{array}{ll}
  \vc{\theta}^{(k,\edge)}, \quad & \gamma_{k}(\vc{\theta}^{(k,\max)}) > 1,\\
  \vc{\theta}^{(k,\max)}, \quad & \gamma_{k}(\vc{\theta}^{(k,\max)}) \le 1,
  \end{array} \right. \qquad k=1,2.
\end{eqnarray*}
  According to Miyazawa \cite{Miyazawa 2009} (see also \cite{Dai-Miyazawa 2011a}), we classify the model into the following three categories.
\begin{mylist}{5}
\item [] Category \rmn{1} \hspace{1ex} $\theta^{(2,\cp)}_{1} < \theta^{(1,\cp)}_{1}$ and $\theta^{(1,\cp)}_{2} < \theta^{(2,\cp)}_{2}$,
\item [] Category \rmn{2} \hspace{1ex} $\theta^{(2,\cp)}_{1} < \theta^{(1,\cp)}_{1}$ and $\theta^{(1,\cp)}_{2} \ge \theta^{(2,\cp)}_{2}$,
\item [] Category \rmn{3} \hspace{1ex} $\theta^{(2,\cp)}_{1} \ge \theta^{(1,\cp)}_{1}$ and $\theta^{(1,\cp)}_{2} < \theta^{(2,\cp)}_{2}$.
\end{mylist}
  Note that it is impossible to have $\theta^{(2,\cp)}_{1} \ge \theta^{(1,\cp)}_{1}$ and $\theta^{(1,\cp)}_{2} \ge \theta^{(2,\cp)}_{2}$ at once because $\theta^{(2,\cp)}_{1} \ge \theta^{(1,\cp)}_{1}$ and the convexity of $\Gamma_{+}$ imply that $\theta^{(1,\cp)}_{2} \le \theta^{(2,\cp)}_{2}$ (see Section 4 of \cite{Miyazawa 2009}). We further note that $\theta^{(1,\cp)}_{2} \ge \theta^{(2,\cp)}_{2}$ can be replaced by $\theta^{(1,\cp)}_{2} = \theta^{(2,\cp)}_{2}$ in Category \rmn{2}. Similarly, $\theta^{(2,\cp)}_{1} \ge \theta^{(1,\cp)}_{1}$ can be replaced by $\theta^{(2,\cp)}_{1} = \theta^{(1,\cp)}_{1}$ in Category \rmn{3}.

  Define the vector $\tau$ as
\begin{eqnarray*}
  \vc{\tau} = \left\{\begin{array}{ll}
  (\theta^{(1,\cp)}_{1}, \theta^{(2,\cp)}_{2}), \quad & \mbox{ for category \rmn{1}}, \\
  (\ol{\xi}_{1}(\theta^{(2,\edge)}_{2}), \theta^{(2,\edge)}_{2}), \quad & \mbox{ for category \rmn{2}}, \\
  (\theta^{(1,\edge)}_{1}, \ol{\xi}_{2}(\theta^{(1,\edge)}_{1})), \quad & \mbox{ for category \rmn{3}}, 
  \end{array} \right.
\end{eqnarray*}
 where $\ol{\xi}_{k}(\theta_{3-k}) = \sup \{\theta_{k}; (\theta_{1}, \theta_{2}) \in \Gamma_{+} \}$. This definition of $\tau$ shows that categories \rmn{1}, \rmn{2} and \rmn{3} are convenient.
 
  We are now ready to present results on the convergence domain $\sr{D}$ and the tail asymptotics obtained by Kobayashi and Miyazawa \cite{KM 2011}. As we mentioned in \sectn{Introduction}, they are obtained for the more general reflecting random walk. Thus, some of their conditions automatically hold for the double QBD process.

\begin{lemma}[Theorem 3.1 of \cite{KM 2011}] {\rm
\label{lem:domain}
\begin{eqnarray}
\label{eqn:domain}
  \sr{D} = \{ \vc{\theta} \in \dd{R}^{2}; \vc{\theta} < \vc{\tau} \mbox{ and } \exists \vc{\theta}' \in \Gamma_{+} \mbox{ such that } \vc{\theta} < \vc{\theta}'\}.
\end{eqnarray}
}\end{lemma}
\begin{figure}[h]
	\caption{The light-green areas are domains $\sr{D}$ for categories \rmn{1} and \rmn{2}}
 	\centering
	\includegraphics[height=2.0in]{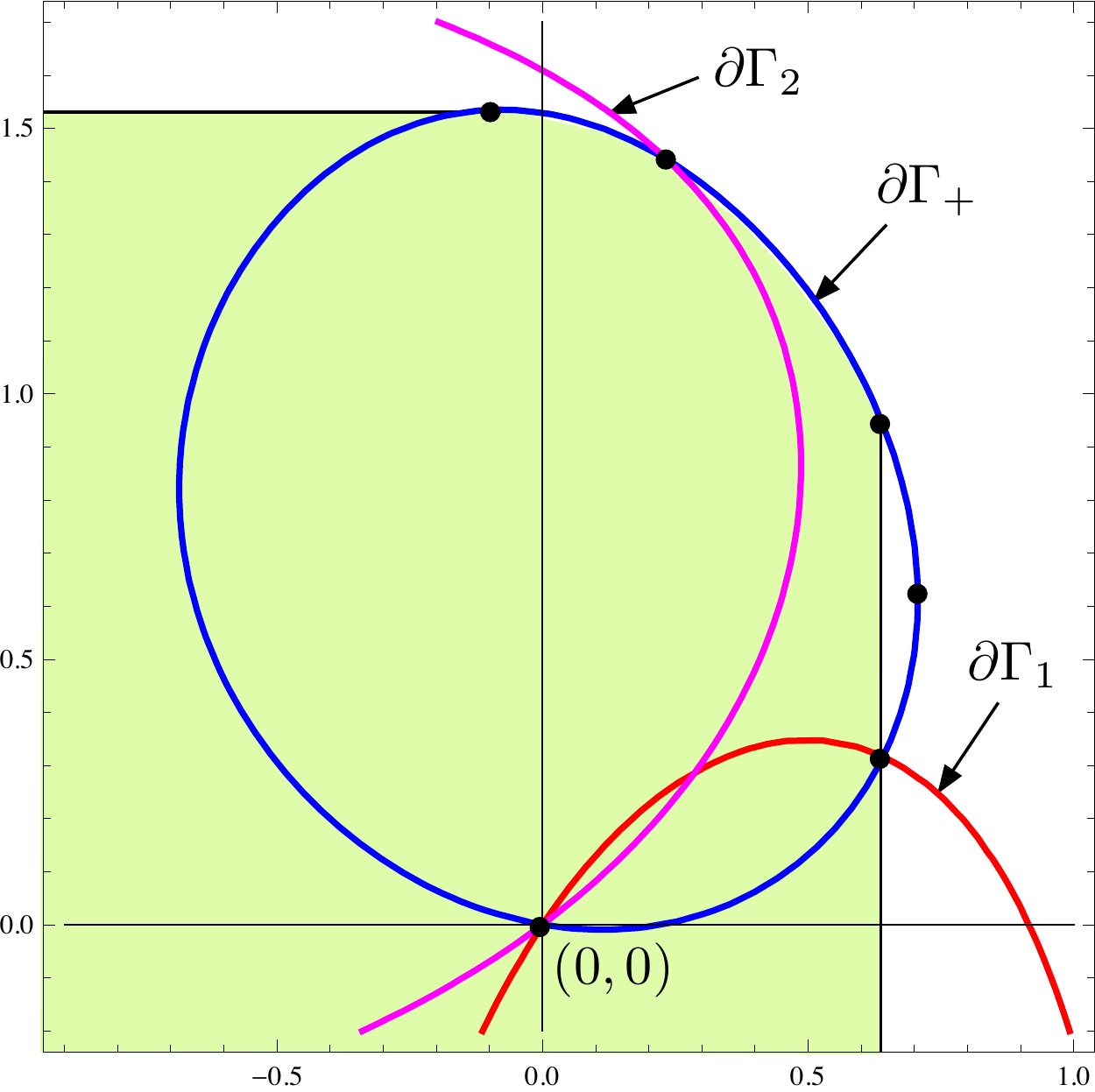} \hspace{5ex}
	\includegraphics[height=2.0in]{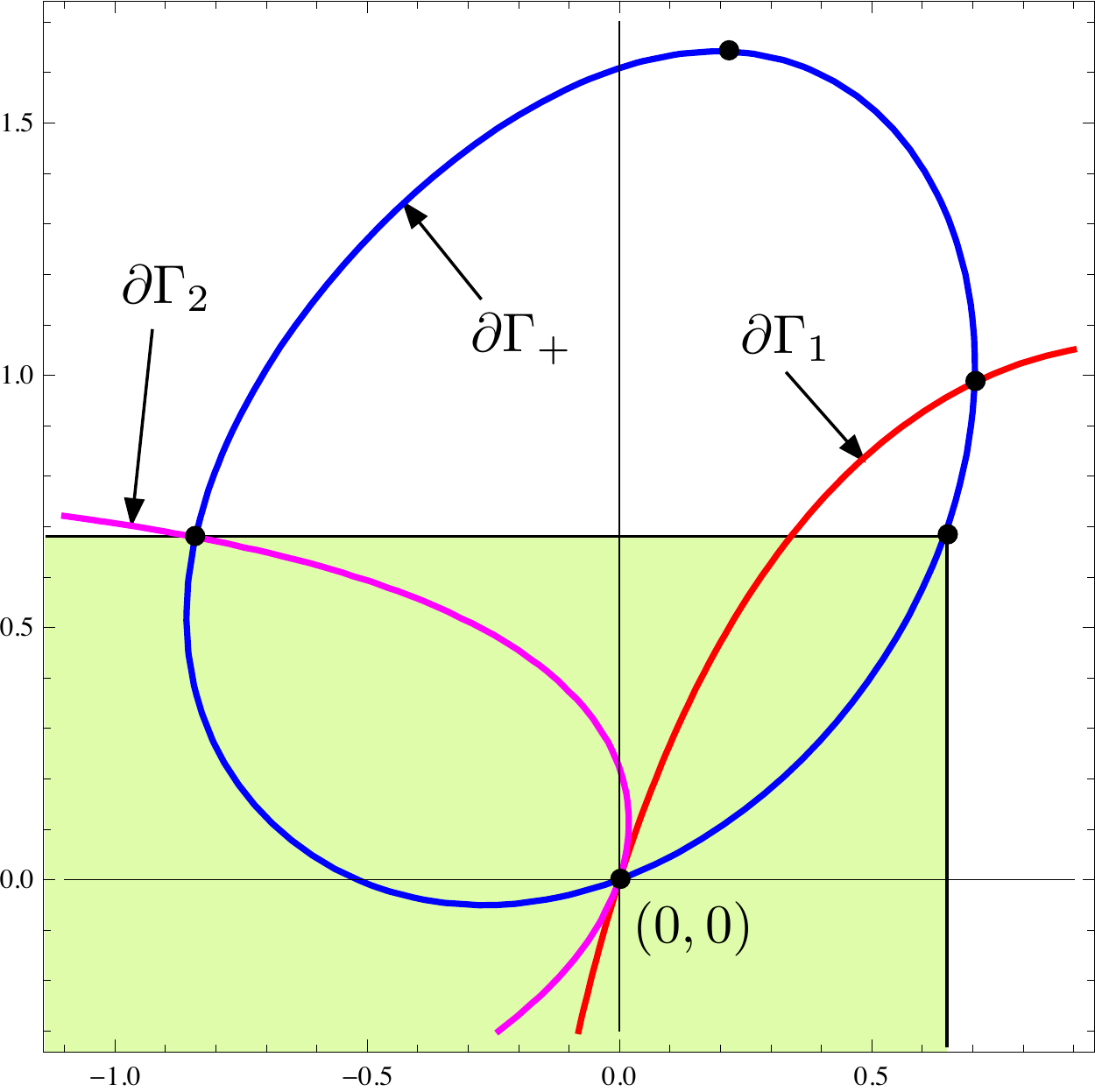} 
	\label{fig:domain 1}
\end{figure}

\begin{theorem}[Theorem 4.2 of \cite{KM 2011}] {\rm
\label{thr:decay rate 1}
  Under conditions (i)--(iv), we have, for $k=1,2$,
\begin{eqnarray}
\label{eqn:decay rate 0}
 \lim_{n \to \infty} \frac 1n \log P( L_{k} \ge n, L_{3-k} = 0 ) = - \tau_{k},
\end{eqnarray}
  and, for any directional vector $\vc{c} \ge \vc{0}$,
\begin{eqnarray}
\label{eqn:decay rate 1}
 \lim_{n \to \infty} \frac 1x \log P( \br{\vc{c}, \vc{L}} \ge x ) = - \alpha_{\vc{c}},
\end{eqnarray}
  where we recall that $\alpha_{\vc{c}} = \sup \{ x \ge 0; x \vc{c} \in \sr{D} \}$. Furthermore, if $\gamma(\alpha_{\vc{c}} \vc{c}) = 1$ and if $\gamma_{k}(\alpha_{\vc{c}} \vc{c}) \ne 1$ and $\alpha_{\vc{c}} c_{k} \ne \theta^{(\infty)}_{kk}$ for $k=1,2$, then we have the following exact asymptotics.
\begin{eqnarray}
\label{eqn:exact asymptotic 2}
  \lim_{x \to \infty} e^{\alpha_{\vc{c}}x} P( \br{\vc{c}, \vc{L}} \ge x ) = b_{\vc{c}}.
\end{eqnarray}
}\end{theorem}

  In this paper, we aim to refine these asymptotics to be exact when $\vc{c}$ is either  $(1,0)$, $(0,1)$ or $(1,1)$. Recall that a sequence of nonnegative number $\{p(n); n \in \dd{Z}_{+}\}$ is said to have the exact asymptotic $(1+b(-1)^{n}) n^{-\kappa} \alpha^{-n}$ for constants $\kappa$ and $\alpha > 1$ if there exist real number $b \in [-1,1]$ and a positive constant $c$ such that
\begin{eqnarray}
\label{eqn:exact f 1}
  \lim_{n \to \infty} (1+b(-1)^{n}) n^{\kappa} \alpha^{n} p(n) = c.
\end{eqnarray}
  We note that this asymptotic is equivalent to
\begin{eqnarray}
\label{eqn:exact f 2}
  \lim_{n \to \infty} (1+b(-1)^{n}) n^{\kappa} \alpha^{n} \sum_{\ell=n}^{\infty} p(\ell) = c'
\end{eqnarray}
  for some $c' > 0$. Thus, there is no difference on the exact asymptotic between $P( L_{k} \ge n )$ and $P( L_{k} = n )$. In what follows, we are mainly concerned with the latter type of exact asymptotics.

\section{Analytic function method} 
\label{sect:Analytic function}
  
  Our basic idea for deriving exact asymptotics is to adapt the method used in \cite{Dai-Miyazawa 2011a} which extends the moment generating functions to complex variable analytic functions, and gets the exact tail asymptotics from analytic behavior around their singular points. A similar method is called a kernel method in some literature \cite{Guillemin-Leeuwaarden 2011,Li-Zhao 2009,Li-Zhao 2011a,Li-Zhao 2011b}. We here call it an analytic function method because our approach heavily use the convergence domain $\sr{D}$, which is not the case for the kernel method. See \cite{Miyazawa 2011} for more details.
  
  There is one problem in adapting the method of \cite{Dai-Miyazawa 2011a} because the moment generating functions $\gamma_{k}(\vc{\theta})$ are not polynomials, while the corresponding functions of SRBM are polynomial. If they are not polynomials, the analytic function approach is hard to apply. This problem is resolved if we use generating functions instead of moment generating functions. We here thanks for the skip free assumption. 
  
\subsection{Convergence domain of a generating function}
\label{sect:Convergence domain}
  
  Let us convert results on moment generating functions to those on generating function, using a mapping from $\vc{z} \equiv (z_{1}, z_{2}) \in \dd{C}$ to $\vc{g}(\vc{z}) \equiv (e^{z_{1}}, e^{z_{2}}) \in \dd{C}$. In particular, for $\vc{\theta} \in \dd{R}^{2}$, $\vc{g}(\vc{\theta}) \in (\dd{R}_{+}^{\rm o})^{2}$, where $\dd{R}_{+}^{\rm o} = (0,\infty)$.  We use the following notations for $k=1,2$.
\begin{eqnarray*}
 && (u^{(k,\min)}_{1}, u^{(k,\min)}_{2})  = \vc{g}(\vc{\theta}^{(k,\min)}), \qquad (u^{(k,\max)}_{1}, u^{(k,\max)}_{2}) = \vc{g}(\vc{\theta}^{(k,\max)}),\\
 && (u^{(k,\edge)}_{1}, u^{(k,\edge)}_{2})  = \vc{g}(\vc{\theta}^{(k,\edge)}), \hspace{10.5ex} (u^{(k,\cp)}_{1}, u^{(k,\cp)}_{2})  = \vc{g}(\vc{\theta}^{(k,\cp)}).\\
 && (\tilde{\tau}_{1}, \tilde{\tau}_{2}) = \vc{g}(\vc{\tau}),
\end{eqnarray*}
    We now transfer the results on the moment generating functions in \sectn{Double QBD} to those on the generating functions. For this, we define
\begin{eqnarray*}
 && \tilde{\sr{D}} = \{ \vc{g}(\vc{\theta}) \in \dd{R}_{+}^{2}; \vc{\theta} \in \sr{D} \}, \\
 && \tilde{\Gamma}_{k} = \{ \vc{g}(\vc{\theta}) \in \dd{R}_{+}^{2}; \vc{\theta} \in \Gamma_{k} \} , \quad k=1,2,+.
\end{eqnarray*}

  Define the following generating functions. For $k = 0,1,2,+$,
\begin{eqnarray*}
  \tilde{\gamma}_{k}(\vc{z}) = E(z_{1}^{X^{(k)}_{1}} z_{2}^{X^{(k)}_{2}}), \quad \vc{z} \equiv (z_{1}, z_{2}) \in \dd{C}^{2},
\end{eqnarray*}
  which exists except for $z_{1} = 0$ or $z_{2} = 0$. Similarly,
\begin{eqnarray*}
 && \tilde{\varphi}(\vc{z}) = E(z_{1}^{L_{1}} z_{2}^{L_{2}}), \\
 && \tilde{\varphi}_{+}(\vc{z}) = E(z_{1}^{L_{1}} z_{2}^{L_{2}} 1(\vc{L} \in S_{+})),\\
  && \tilde{\varphi}_{k}(z_{k}) = E(z_{k}^{L_{k}} 1(\vc{L} \in S_{k})), \qquad k=1,2.
\end{eqnarray*}
as long as they exist.

  Obviously, these generating functions are obtained from the corresponding moment generating functions using the inverse mapping $\vc{g}^{-1}$.
\begin{eqnarray*}
  && \tilde{\gamma}_{k}(\vc{z}) = \gamma_{k}(\log z_{1}, \log z_{2}), \quad k=0,1,2,+, \qquad\\
 && \tilde{\varphi}(\vc{z}) = \varphi(\log z_{1}, \log z_{2}), \\
 && \tilde{\varphi}_{+}(\vc{z}) = \varphi_{+}(\log z_{1}, \log z_{2}),\\
 && \tilde{\varphi}_{k}(z) = \varphi_{k}(\log z), \qquad k=1,2,
\end{eqnarray*}
Then, the stationary equation \eqn{stationary equation 1} can be written as
\begin{eqnarray}
\label{eqn:stationary equation 2}
 \lefteqn{ (1 -  \tilde{\gamma}_{+}(\vc{z}))  \tilde{\varphi}_{+}(\vc{z}) + (1 -  \tilde{\gamma}_{1}(\vc{z}))  \tilde{\varphi}_{1}(z_{1})}\nonumber \\
 && \hspace{10ex}+ (1 - \tilde{\gamma}_{2}(\vc{z}))  \tilde{\varphi}_{2}(z_{2}) + (1 -  \tilde{\gamma}_{0}(\vc{z})) \nu(\vc{0}) = 0. \qquad
\end{eqnarray}

\begin{figure}[h]
	\caption{The examples of $\sr{D}$ and the corresponding $\tilde{\sr{D}}$, which may not be convex, where $(p_{21}, p_{01}, p_{11}) = (0.1, 0.1, 0.7),
(p_{20}, p_{00}, p_{10}) = (1.5, 0.5, 0.5),
(p_{22}, p_{02}, p_{12}) = (2, 3, 1)$ for $p_{ij} = P(\vc{X}^{(+)}=(i,j))$.
} 
 	\centering
	\includegraphics[height=2.0in]{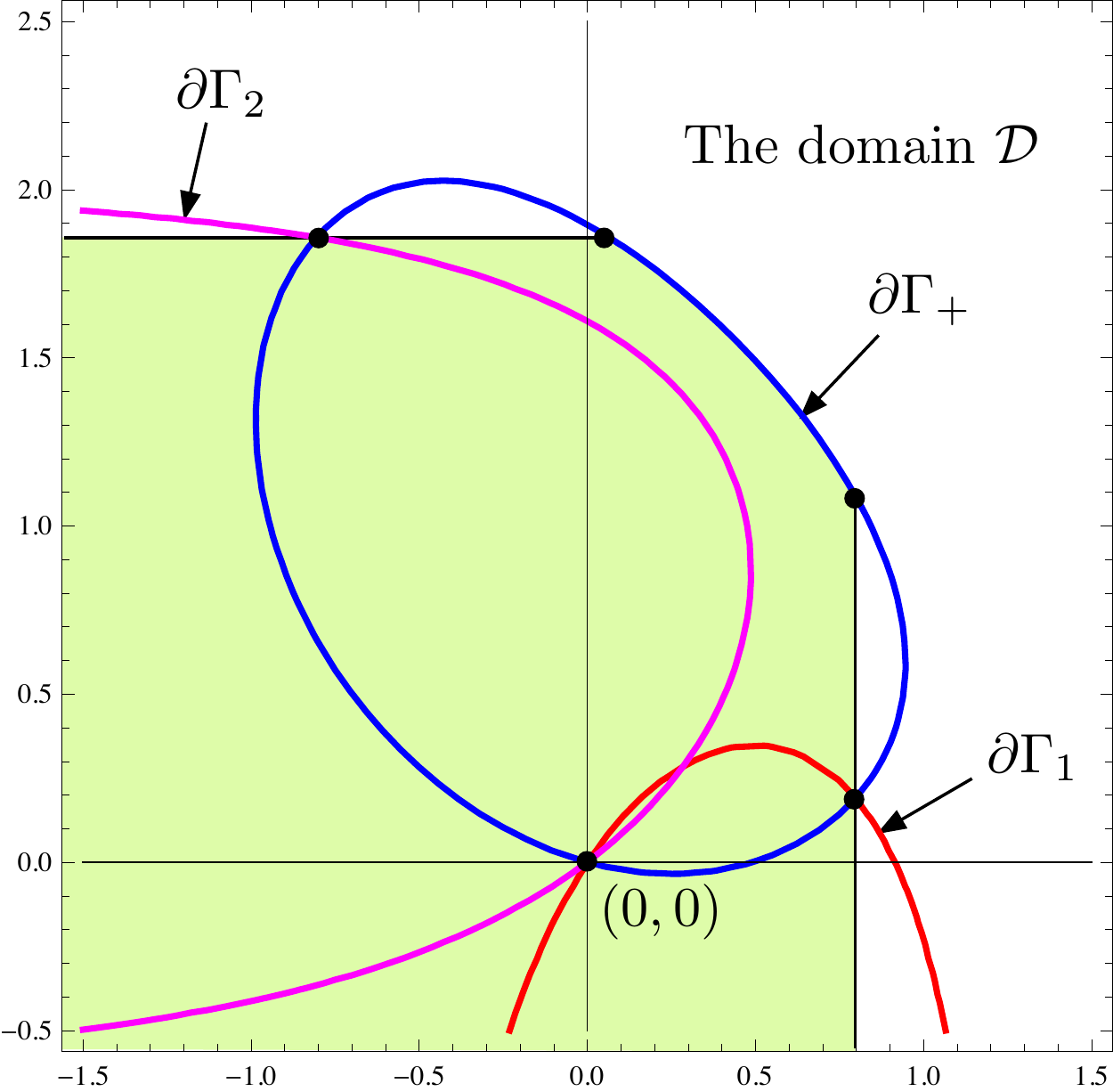} \hspace{5ex}
	\includegraphics[height=2.0in]{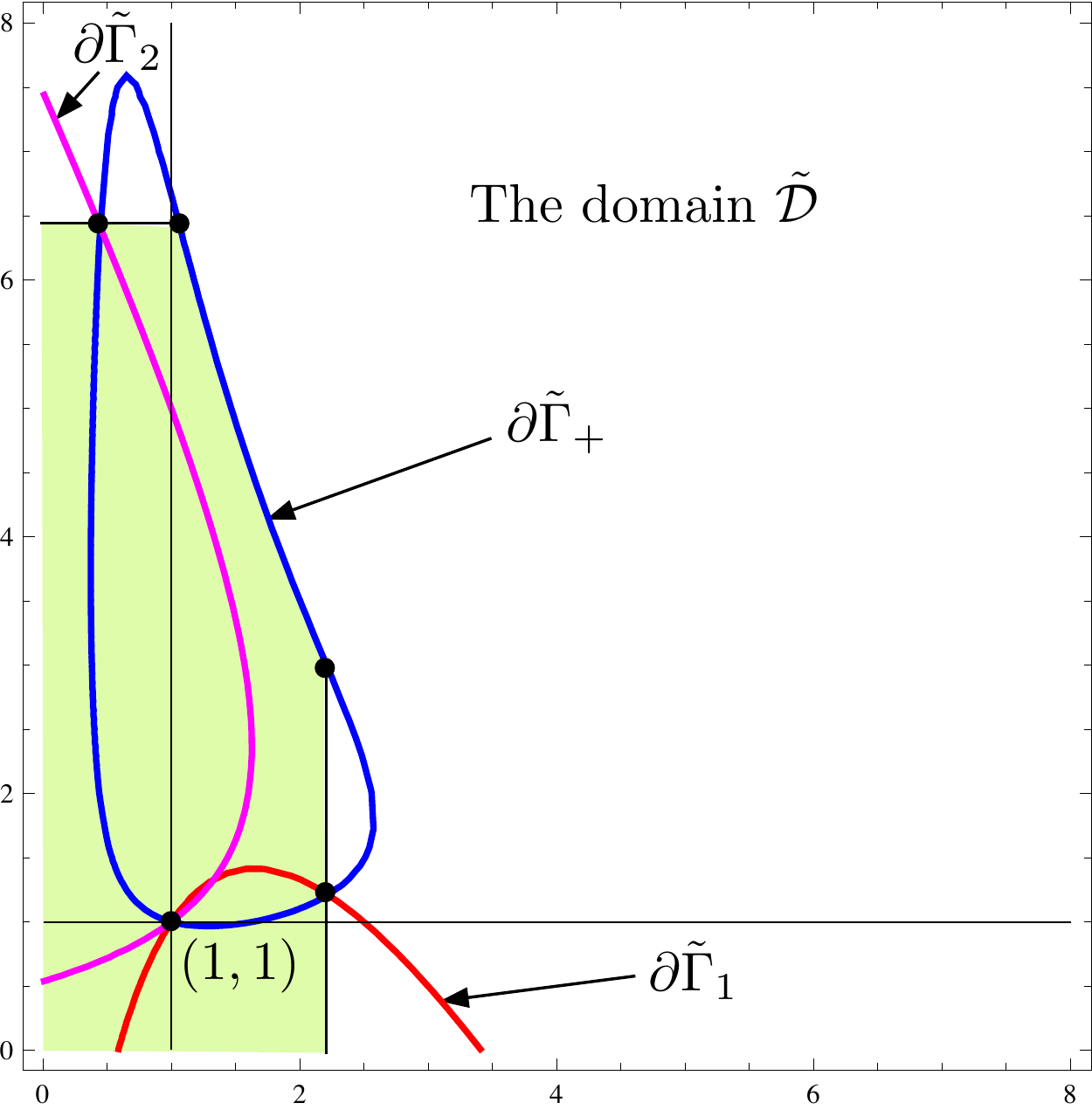} 
	\label{fig:non-convex}
\end{figure}

   It is easy to see that
\begin{eqnarray*}
 && \tilde{\Gamma}_{k} \equiv \{\vc{u} \in \dd{R}_{+}^{2}; \vc{u} > \vc{0}, \tilde{\gamma}_{k}(\vc{u}) < 1 \}, \quad k=1,2,+,\\
 && \tilde{\sr{D}} \equiv \{\vc{u} \in \dd{R}_{+}^{2}; \vc{u} > \vc{0}, \tilde{\varphi}(\vc{u}) < \infty \},
\end{eqnarray*}
  It is notable that these sets may not be convex because two dimensional generating functions may not be convex (see \fig{non-convex}). 
  Nevertheless, they still have nice properties because the generating functions are polynomials with nonnegative coefficients. To make this specific, we introduce the following terminology.
  
\begin{definition} {\rm
\label{dfn:directed-convex}
  A subset $A$ of $\dd{R}^{2}$ is said to be nonnegative-directed (or coordinate-directed) convex if $\lambda \vc{x} + (1-\lambda) \vc{y} \in A$ for any number $\lambda \in [0,1]$ and any $\vc{x}, \vc{y} \in A$ such that $\vc{y} - \vc{x} \ge \vc{0}$ (or $\vc{y} - \vc{x}$ in either one of the coordinate axes, respectively).
}\end{definition}

  We then immediately have the following facts.
\begin{lemma} {\rm
\label{lem:directed-convex}
  $\tilde{\sr{D}}$ is nonnegative-directed convex, and $\tilde{\Gamma}_{k}$ is coordinate-directed convex for for $k=+,0,1,2$.
}\end{lemma}

  Note that \eqn{stationary equation 2} is valid for $\vc{z} \in \dd{C}^{2}$ satisfying $|\vc{z}| \in \tilde{\sr{D}}$ because $|\tilde{\varphi}(\vc{z})| \le \tilde{\varphi}(|\vc{z}|)$. Furthermore, 
\begin{eqnarray*}
  \lefteqn{\{\vc{z} \in \dd{C}^{2}; |\vc{z}| \in \tilde{\sr{D}} \} = \{\vc{z} \in \dd{C}^{2}; E(e^{L_{1} \log|z_{1}| + L_{2}\log|z_{2}|}) < \infty \}}\\
  && = \{\vc{g}(\log|z_{1}| + i \arg z_{1}, \log|z_{2}| + i \arg z_{2}); \vc{z} \in \dd{C}^{2}, (\log|z_{1}|, \log|z_{2}|) \in \sr{D} \}\\
  && = \vc{g}(\{\vc{z} \in \dd{C}^{2}, (\Re{z_{1}}, \Re{z_{2}}) \in \sr{D} \}),
\end{eqnarray*}
  where $|\vc{z}| = (|z_{1}|, |z_{2}|)$. Hence, the domain $\sr{D}$ is well transferred to $\tilde{\sr{D}}$. We will work on $\tilde{\sr{D}}$ for finding the analytic behaviors of $\tilde{\varphi}_{1}(z)$ and $\tilde{\varphi}_{2}(z)$ around their dominant singular points. This is different from the kernel method, which directly works on the set of complex vectors $\vc{z}$ satisfying $\tilde{\gamma}_{+}(\vc{z}) = 1$, and applies deeper complex analysis such as analytic extension on a Riemann surface (e.g., see \cite{FIM 1999}). We avoid it using the domain $\tilde{\sr{D}}$. 

\subsection{A key function for analytic extension}
\label{sect:A key function}
  
  Once the domain $\tilde{\sr{D}}$ is obtained, the next step is to study analytic behaviors of the generating function $\tilde{\varphi}_{k}$ for $k=1,2$. For this, we use a relation between them by letting $\tilde{\gamma}_{+}(\vc{z}) - 1 = 0$ in the stationary equation \eqn{stationary equation 2}, which removes $\tilde{\varphi}_{+}(\vc{z})$. For this, let us consider the solution $u_{2} > 0$ of $\tilde{\gamma}_{+}(u_{1}, u_{2}) = 1$ for each fixed $u_{1} > 0$. Since this equation is quadratic concerning $u_{2}$ and $\tilde{\sr{D}} \subset (\dd{R}^{\rm o})_{+}^{2}$, it has two positive solutions for each $u_{1}$ satisfying
\begin{eqnarray*}
  u^{(1,\min)}_{1} \le u_{1} \le u^{(1,\max)}_{1}.
\end{eqnarray*}
  Denote these solutions by $\ul{\zeta}_{2}(u_{1})$ and $\ol{\zeta}_{2}(u_{1})$ such that $\ul{\zeta}_{2}(u_{1}) \le \ol{\zeta}_{2}(u_{1})$. Similarly, $\ul{\zeta}_{1}(u_{2})$ and $\ol{\zeta}_{1}(u_{2})$ are defined for $u_{2}$ satisfying
\begin{eqnarray*}
  u^{(2,\min)}_{2} \le u_{2} \le u^{(2,\max)}_{2}.
\end{eqnarray*}
  One can see these facts also applying the mapping $\vc{g}$ to the convex bounded set $\sr{D}$ (see \lem{directed-convex}).
  
We now adapt the arguments in \cite{Dai-Miyazawa 2011a}. For this, we first examine the function $\ul{\zeta}_{2}$. Let
\begin{eqnarray*}
  && p_{*k}(u) = E(u^{X^{(+)}_{1}} 1(X^{(+)}_{2} = k)), \\
 && p_{k*}(u) = E(u^{X^{(+)}_{2}} 1(X^{(+)}_{1} = k)), \quad k=0,1,-1. \qquad
\end{eqnarray*}
  Then, $\tilde{\gamma}_{+}(u_{1}, u_{2}) = 1$ can be written as
\begin{eqnarray}
\label{eqn:gamma + u}
  u_{2}^{2} p_{*1}(u_{1}) -  u_{2} (1 - p_{*0}(u_{1})) +  p_{*-1}(u_{1}) = 0.
\end{eqnarray}
  Hence, we have, for $u \in [u^{(1,\min)}_{1}, u^{(1,\max)}_{1}]$,
\begin{eqnarray}
\label{eqn:ul zeta 2}
  \ul{\zeta}_{2}(u) = \frac {1 - p_{*0}(u) - \sqrt{D_{2}(u)}} {2 p_{*1}(u)},
\end{eqnarray}
  where
\begin{eqnarray*}
  D_{2}(u) = (1 - p_{*0}(u))^{2} - 4 p_{*1}(u) p_{*-1}(u) \ge 0.
\end{eqnarray*}
  Since $D_{2}(u^{(1,\min)}_{1}) = D_{2}(u^{(1,\max)}_{1}) = 0$ and $u^{2} D_{2}(u)$ is a polynomial with order 4 at most and order 2 at least by condition (i), $u^{2} D_{2}(u)$ can be factorized as
\begin{eqnarray*}
  u^{2} D_{2}(u) = (u - u^{(1,\min)}_{1}) (u^{(1,\max)}_{1} - u) h_{2}(u),
\end{eqnarray*}
  where $h_{2}(u) \ne 0$ for $u \in (u^{(1,\min)}_{1}, u^{(1,\max)}_{1})$. This fact can be verified by the mapping $\vc{g}$ from $\Gamma_{+}$ to $\tilde{\Gamma}_{+}$.
  
  To get tail asymptotics, we will use analytic functions. So far, we like to analytically extend the function $\ul{\zeta}_{2}$ from the real interval to a sufficiently large region in the complex plane $\dd{C}$. For this, we prepare a series of lemmas. We first refer to the following fact.

\begin{lemma}[Lemma 2.3.8 of \cite{FIM 1999}] {\rm
\label{lem:D null points}
  All the solutions of $z^{2} D_{2}(z) = 0$ for $z \in \dd{C}$ are real numbers.
}\end{lemma}

In the light of the above arguments, this lemma immediately leads to the following fact.

\begin{lemma} {\rm
\label{lem:real root}
  $z^{2} D_{2}(z) = 0$ for $z \in \dd{C}$ has no solution in the region such that $|z| \in (u^{(1,\min)}_{1}, u^{(1,\max)}_{1})$.
}\end{lemma}

We will also use the following two lemmas, which show how the periodic nature of the random walk $\{\vc{Y}_{\ell}\}$ is related to the branch points (see \rem{condition (v-b)} for the periodic nature). They are proved in Appendices A and B, respectively.

\begin{lemma} {\rm
\label{lem:D root}
  The equation:
\begin{eqnarray}
\label{eqn:D root}
    D_{2}(z) = 0, \qquad |z| = u^{(1,\max)}_{1}, \; z \in \dd{C},
\end{eqnarray}
  has only one solution $z = u^{(1,\max)}_{1}$ if and only if (v-a) holds. Otherwise, it has two solutions $z = \pm u^{(1,\max)}_{1}$, and $u^{2} D_{2}(u)$ is an even function.
}\end{lemma}

\begin{lemma} {\rm
\label{lem:arithmetic condition 1}
  For each fixed $x,y > 0$, we have
  \begin{eqnarray}
\label{eqn:arithmetic condition 1}
  \tilde{\gamma}_{+}(x, y) = 1, \qquad \tilde{\gamma}_{+}(-x, -y) = 1,
\end{eqnarray}
 if and only if (v-a) does not hold.
}\end{lemma}

\begin{remark} {\rm
\label{rem:D root 1}
  \lem{D root} is essentially the same as Remark 3.1 of \cite{Li-Zhao 2011b}, which is obtained as a corollary of their Lemma 3.1, which is immediate from Lemmas 2.3.8 of  \cite{FIM 1999}.
}\end{remark}

  By Lemmas \lemt{D null points} and \lemt{real root}, $\ul{\zeta}_{2}(u)$ on $(u^{(1,\min)}_{1}, u^{(1,\max)}_{1})$ is extendable as an analytic function of a complex variable to the region $\tilde{\sr{G}}_{0}(u^{(1,\min)}_{1}, u^{(1,\max)}_{1})$, where
\begin{eqnarray*}
   \tilde{\sr{G}}_{0}(a,b) = \{ z \in \dd{C}; z \not \in (-\infty, a] \cup [b, \infty) \}, \qquad a, b \in \dd{R},
\end{eqnarray*}
  and has a single branch point $u^{(1,\max)}_{1}$ on $|z| = u^{(1,\max)}_{1}$ if (v-a) hold, and two branch points $\pm u^{(1,\max)}_{1}$ there otherwise by Lemmas \lemt{D root} and \lemt{real root}. Both branch points have order two. We denote this extended analytic function by $\ul{\zeta}_{2}(z)$. That is, we use the same notation for an analytically extended function. We identify it by its argument. The following lemma is a key for our arguments. The idea of this lemma is similar to Lemma 6.3 of \cite{Dai-Miyazawa 2011a}, but its proof is entirely different from that lemma.

\begin{lemma} {\rm
\label{lem:extension 1}
  (a) $\ul{\zeta}_{2}$ of \eqn{ul zeta 2} is analytically extended on $\tilde{\sr{G}}_{0}(u^{(1,\min)}_{1}, u^{(1,\max)}_{1})$.\\ 
  (b) For $z \in \dd{C} \mbox{ satisfying } |z| \in (u^{(1,\min)}_{1}, u^{(1,\max)}_{1}]$,
\begin{eqnarray}
\label{eqn:bound zeta}
  |\ul{\zeta}_{2}(z)| \le \ul{\zeta}_{2}(|z|) \le u^{(1,\max)}_{2},
\end{eqnarray}
where the second inequality is strict if $|z| < u^{(1,\max)}_{1}$.\\
  (c) If either $m^{(1)}_{2} = 0$ or (v-b) holds, then
\begin{eqnarray}
\label{eqn:no solution 1}
  \tilde{\gamma}_{1}(z, \ul{\zeta}_{2}(z)) = 1, \qquad |z| = u^{(1,\edge)}_{1},
\end{eqnarray}
has no solution other than $z = u^{(1,\edge)}_{1}$.\\
(d) The equation \eqn{no solution 1} has two solutions $z = \pm u^{(1,\edge)}_{1}$ if and only if neither $m^{(1)}_{2} = 0$, (v-a) nor (v-b) holds.

\begin{proof}
  We have already proved (a). Thus, we only need to prove (b), (c) and (d). We first prove (b). For this, it is sufficient to prove \eqn{bound zeta} for $|z| < u^{(1,\max)}_{1}$ by the continuity of $\ul{\zeta}_{2}(z)$ for $|z| \le u^{(1,\max)}_{1}$ at $z = u^{(1,\max)}_{1}$. Substituting complex numbers $z_{1}$ and $z_{2}$ into $u_{1}$ and $u_{2}$ of \eqn{gamma + u}, we have
\begin{eqnarray}
\label{eqn:gamma z 1}
   z_{2}^{2} p_{*1}(z_{1}) + z_{2} p_{*0}(z_{1}) +  p_{*-1}(z_{1}) = z_{2}.
\end{eqnarray}
 Obviously, this equation has the following solutions for each fixed $z_{1}$ such that $|z_{1}| \in (u^{(1,\min)}_{1}, u^{(1,\max)}_{1})$.
\begin{eqnarray}
\label{eqn:gamma z 2}
  z_{2} = \ul{\zeta}_{2}(z_{1}), \quad \ol{\zeta}_{2}(z_{1}).
\end{eqnarray}

  We next take the absolute values of both sides of \eqn{gamma z 1}, then
\begin{eqnarray*}
  |z_{2}|^{2} p_{*1}(|z_{1}|) + |z_{2}| p_{*0}(|z_{1}|) +  p_{*-1}(|z_{1}|) \ge |z_{2}|.
\end{eqnarray*}
  Thus, we get
\begin{eqnarray*}
  |z_{2}| (\tilde{\gamma}_{+}(|z_{1}|, |z_{2}|) - 1) \ge 0.
\end{eqnarray*}
  By the definitions of $\ul{\zeta}_{2}(|z_{1}|)$ and $\ol{\zeta}_{2}(|z_{1}|)$, this inequality can be written as
\begin{eqnarray*}
 (|z_{2}|- \ul{\zeta}_{2}(|z_{1}|)) (|z_{2}|- \ol{\zeta}_{2}(|z_{1}|)) = |z_{2}|(\tilde{\gamma}_{+}(|z_{1}|, |z_{2}|) - 1) \ge 0.
\end{eqnarray*}
  Hence, $\ul{\zeta}_{2}(|z_{1}|) \le \ol{\zeta}_{2}(|z_{1}|)$ implies
\begin{eqnarray}
\label{eqn:z 2 zeta}
  |z_{2}| \le \ul{\zeta}_{2}(|z_{1}|) \quad \mbox{ or } \quad \ol{\zeta}_{2}(|z_{1}|) \le |z_{2}|.
\end{eqnarray}
  By \eqn{gamma z 2}, we can substitute $z_{2} = \ul{\zeta}_{2}(z_{1})$ into \eqn{z 2 zeta}, and get
\begin{eqnarray}
\label{eqn:z 1 zeta}
  |\ul{\zeta}_{2}(z_{1})| \le \ul{\zeta}_{2}(|z_{1}|) \quad \mbox{or} \quad \ol{\zeta}_{2}(|z_{1}|) \le |\ul{\zeta}_{2}(z_{1})|, \qquad |z_{1}| \in (u^{(1,\min)}_{1}, u^{(1,\max)}_{1}). \quad
\end{eqnarray}

  Thus, (b) is obtained if we show that $\ol{\zeta}_{2}(|z_{1}|) \le |\ul{\zeta}_{2}(z_{1})|$ is impossible. Suppose the contrary of this, that is, there is a $z_{1}^{(0)}$ such that 
\begin{eqnarray}
\label{eqn:z (0)}
  \ol{\zeta}_{2}(|z_{1}^{(0)}|) \le |\ul{\zeta}_{2}(z_{1}^{(0)})|, \qquad |z_{1}^{(0)}| \in (u^{(1,\min)}_{1}, u^{(1,\max)}_{1}).
\end{eqnarray}
  Since $|\ul{\zeta}_{2}(z)|$ is continuous and converges to $\ul{\zeta}_{2}(|z^{(0)}_{1}|)$ as $z$ goes to $|z_{1}^{(0)}|$ on the path that $|z| = |z_{1}^{(0)}|$, there must be a $z_{1}^{(1)}$ such that $|z_{1}^{(1)}| = |z_{1}^{(0)}|$ and
\begin{eqnarray*}
   \ul{\zeta}_{2}(|z_{1}^{(1)}|) < |\ul{\zeta}_{2}(z_{1}^{(1)})| < \ol{\zeta}_{2}(|z_{1}^{(1)}|)
\end{eqnarray*}
  Since $|z_{1}^{(1)}| = |z_{1}^{(0)}| \in (u^{(1,\min)}_{1}, u^{(1,\max)}_{1})$, this contradicts \eqn{z 1 zeta}, which proves (b).
  
  We next prove (c). Let
\begin{eqnarray*}
  p^{(1)}_{*k}(z) = E(z^{X^{(1)}_{1}} 1(X^{(1)}_{2} = k)), \qquad k = 0, 1.
\end{eqnarray*}
  First, assume that $m^{(1)}_{2} = 0$. This implies $p^{(1)}_{*1}(z) = 0$, and therefore \eqn{no solution 1} is reduced to $p^{(1)}_{*0}(z) = 1$. Hence, its solution is $z=1$ or $z=p^{(1)}_{-10}/ p^{(1)}_{10} \ge 0$ if $p^{(1)}_{10} \ne 0$ (otherwise, $z=1$ is only the solution). Both are nonnegative numbers, and therefore \eqn{no solution 1} has no solution $z$ such that
\begin{eqnarray}
\label{eqn:z u 1}
  |z| = u^{(1,\edge)}_{1}, \qquad z \ne u^{(1,\edge)}_{1}.
\end{eqnarray}
  
  We next assume that $m^{(1)}_{2} \ne 0$, which implies $p^{(1)}_{*1}(z) \ne 0$. 
  Since \eqn{no solution 1} can be written as
\begin{eqnarray}
\label{eqn:no solution 2}
   \ul{\zeta}_{2}(z) p^{(1)}_{*1}(z) + p^{(1)}_{*0}(z) = 1
\end{eqnarray}
  and $1 \le |w| + |1-w|$ for any $w \in \dd{C}$, we have
\begin{eqnarray}
\label{eqn:zeta inequality 1}
  |\ul{\zeta}_{2}(z)| = \left| \frac {1 - p^{(1)}_{*0}(z)} {p^{(1)}_{*1}(z)} \right| \ge \frac {1 - |p^{(1)}_{*0}(z)|} {|p^{(1)}_{*1}(z)|} \ge \frac {1 - p^{(1)}_{*0}(|z|)} {p^{(1)}_{*1}(|z|)} = \ul{\zeta}_{2}(|z|).
\end{eqnarray}
  If \eqn{z u 1} holds, then both sides of this inequality are identical if and only if (v-b) does not hold. Hence,  if (v-b) holds, then $|\ul{\zeta}_{2}(z)| > \ul{\zeta}_{2}(|z|)$, and therefore \eqn{no solution 1} has no solution satisfying \eqn{z u 1} because of \eqn{bound zeta}.
  
  We finally prove (d). For this, we assume that both of $m^{(2)}_{1} = 0$ and (v-b) do not hold. In this case, $p^{(1)}_{01} = p^{(1)}_{(-1)0} = p^{(1)}_{10} = 0$, so it follows from \eqn{no solution 2} that
\begin{eqnarray*}
  \ul{\zeta}_{2}(z) = \frac {(1-p^{(1)}_{00})z} {p^{(1)}_{-11} + p^{(1)}_{11} z^{2}}.
\end{eqnarray*}
  Hence, if \eqn{z u 1} holds, then we must have $z = - u^{(1,\edge)}_{1}$ because of \eqn{bound zeta} and \eqn{zeta inequality 1}. By the above equation, we also have $ \ul{\zeta}_{2}(-u^{(1,\edge)}_{1}) = - \ul{\zeta}_{2}(u^{(1,\edge)}_{1})$. Hence, we need to check whether $(- u^{(1,\edge)}_{1}, - \ul{\zeta}_{2}(u^{(1,\edge)}_{1}))$ be the solution of $\gamma_{+}(x,y)=1$. By \lem{arithmetic condition 1}, $z = - u^{(1,\edge)}_{1}$ is the solution of \eqn{no solution 1} if and only if (v-a) does not hold. Combining this with (b) and (c) completes the proof of (d). \pend
\end{proof}
}\end{lemma}

\begin{table}[t] \normalsize
\caption{The solutions of \eqn{D root} and \eqn{no solution 1}, where $\bigcirc$, {\large $\times$} and {\large $-$} indicate ``yes'', ``no'' and `"irrelevant''.}
\begin{center}
\begin{tabular}{|c|c|c|c|c|c|c|}
\hline

Non-arithmetic: (v-a) & $\bigcirc$ & {\large $\hspace{1.5ex} \times$ \hspace{1ex} } & {\large $\hspace{1.5ex} \times \hspace{1.5ex}$} & {\large $\hspace{1.5ex} \times \hspace{1.5ex} $} & {\large $\times$} \\ \hline
Non-arithmetic: (v-b) &  $-$ & $\bigcirc$ & $\bigcirc$ & {\large $\times$} & {\large $\times$} \\ \hline
$m^{(1)}_{2} = 0$ & $-$ & $\bigcirc$ & {\large $\times$} & $\bigcirc$ & {\large $\times$} \\ \hline \hline
The solutions of \eqn{D root} & $\hspace{3ex} u^{(1,\max)}_{1} \hspace{3ex}$ & \multicolumn{4}{|c|}{$\pm u^{(1,\max)}_{1}$} \\ \hline
$\;$ The solutions of \eqn{no solution 1} $\;$ & $u^{(1,\edge)}_{1}$ & \multicolumn{3}{|c|}{\hspace{3ex}  $u^{(1,\edge)}_{1}$ \hspace{3ex} } & $\pm u^{(1,\edge)}_{1}$ \\ \hline
\end{tabular}
\end{center}
\label{tab:period}
\end{table}%

  For convenience of later reference, we summarize the results in (c) and (d) of \lem{extension 1} in \tab{period}. Similar results can be obtained in the direction of the 2nd axes using (v-b) and $m^{(2)}_{1} = 0$ instead of (v-a) and $m^{(1)}_{2} = 0$. Since the results are symmetric, we omit them. We remark that Li and Zhao \cite{Li-Zhao 2011b} have not considered the cases $m^{(2)}_{1} = 0$ and $m^{(1)}_{2} = 0$, which seems to be overlooked.

\subsection{Nature of the dominant singularity}
\label{sect:Nature of}

  We consider complex variable functions $\tilde{\varphi}_{1}(z_{1})$ and $\tilde{\varphi}_{2}(z_{2})$. Recall that
\begin{eqnarray}
\label{eqn:decomposition 1}
 \tilde{\varphi}(\vc{z}) = \tilde{\varphi}_{+}(\vc{z}) + \tilde{\varphi}_{1}(z_{1}) + \tilde{\varphi}_{2}(z_{2}) + \nu(\vc{0}).
\end{eqnarray}
  Obviously, $\tilde{\varphi}(\vc{z})$ is analytic for $\vc{z} \in \dd{C}^{2}$ such that $(|z_{1}|,|z_{2}|) \in \tilde{\sr{D}}$, and singular on the boundary of $\tilde{\sr{D}}$. This implies that $\tilde{\varphi}_{i}(z_{i})$ is analytic for $|z_{i}| < \tilde{\tau}_{i}$ and has a point on the circle $|z| = \tilde{\tau}_{i}$. This is easily seen from \eqn{decomposition 1} with $z_{j} = 0$ for $j = 3-i$. Furthermore, $z_{i} = \tilde{\tau}_{i}$ must be a singular point for $i=1,2$ by Pringsheim's theorem (see, e.g., Theorem 17.13 in Volume 1 of Markushevich \cite{Markushevich 1977}). In addition to this point, we need to find all singular points on $|z| = \tilde{\tau}_{i}$ to get the tail asymptotics as we will see. As expected from \lem{D root}, $z = - \tilde{\tau}_{i}$ may be another singular point, which occurs only when (v-a) does not hold.
  
  We focus on these singular points instead of searching singular points on $|z| = \tilde{\tau}_{i}$, and show that there is no other singular point on the circle through analytic behavior of $\tilde{\varphi}_{i}(z)$. Since results are symmetric for $\tilde{\varphi}_{1}(z)$ and $\tilde{\varphi}_{2}(z)$, we only consider $\tilde{\varphi}_{1}(z)$ in this section.
  
\begin{figure}[t]
	\caption{The shaded area is $\tilde{\sr{G}}^{-}_{\delta}(-u^{(1,\max)}_{1}) \cap \tilde{\sr{G}}^{+}_{\delta}(u^{(1,\max)}_{1})$} 
 	\centering
	\includegraphics[height=5.5cm]{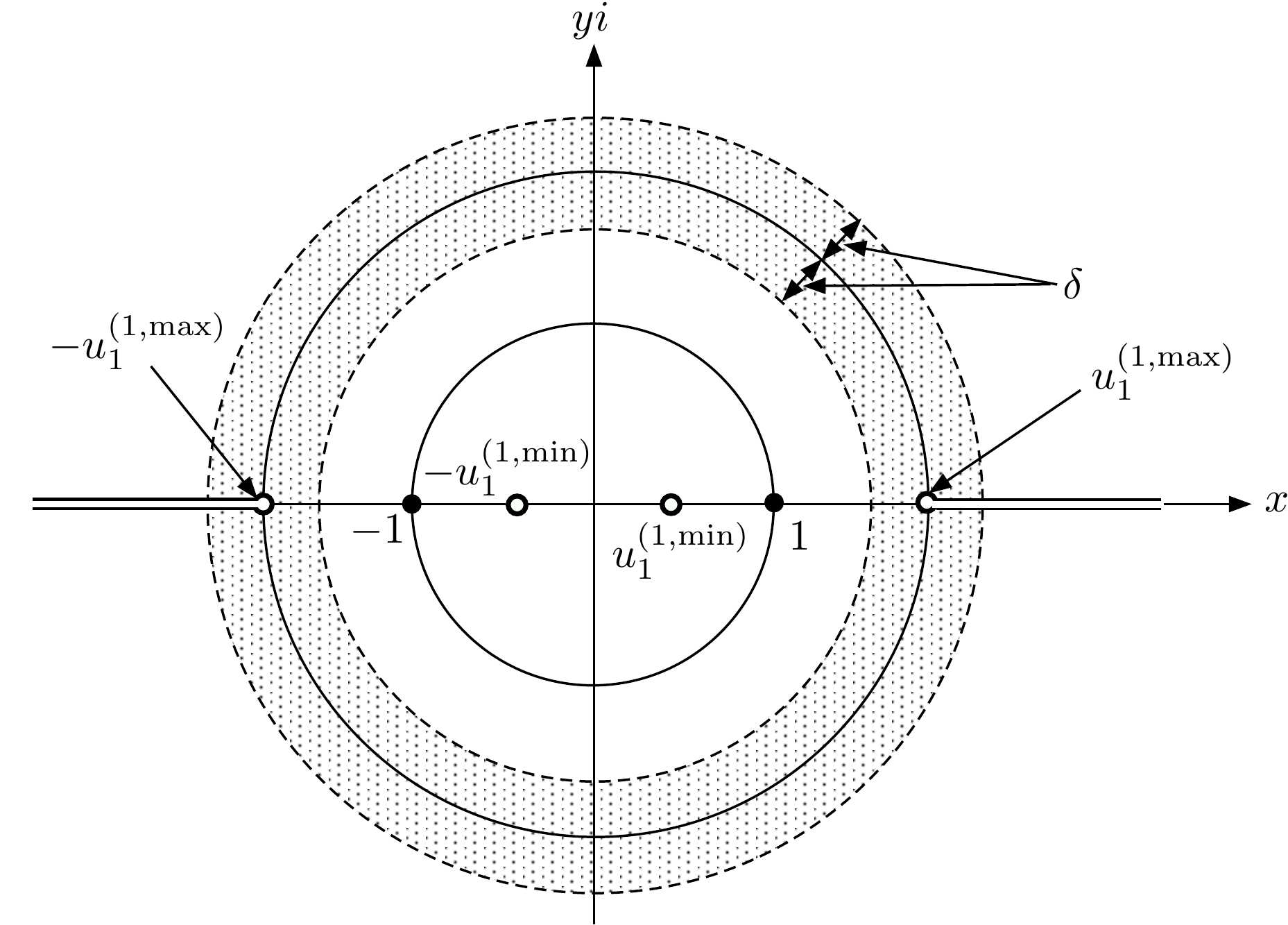} 
	\label{fig:region 1}
\end{figure}

  For this, we use the stationary equation \eqn{stationary equation 2}, which is valid on $\tilde{\sr{D}}$. Plugging $(z_{1}, z_{2}) = (z, \ul{\zeta}_{2}(z))$ into \eqn{stationary equation 2} yields, for $|z| \in (u^{(1,\min)}_{1}, \tilde{\tau}_{1})$,
\begin{eqnarray}
\label{eqn:stationary equation 3}
 \tilde{\varphi}_{1}(z) = \frac{(\tilde{\gamma}_{2}(z, \ul{\zeta}_{2}(z)) - 1)  \tilde{\varphi}_{2}(\ul{\zeta}_{2}(z))}{1 -  \tilde{\gamma}_{1}(z, \ul{\zeta}_{2}(z))} + \frac {(\tilde{\gamma}_{0}(z, \ul{\zeta}_{2}(z)) - 1) \nu(\vc{0})} {1 -  \tilde{\gamma}_{1}(z, \ul{\zeta}_{2}(z))}.
\end{eqnarray}
  In the light of this equation, the dominant singularity of $\tilde{\varphi}_{1}(z)$ is caused by $\ul{\zeta}_{2}(z)$, $\tilde{\varphi}_{2}(\ul{\zeta}_{2}(z))$ or
\begin{eqnarray}
\label{eqn:edge 1}
  \tilde{\gamma}_{1}(z, \ul{\zeta}_{2}(z)) = 1.
\end{eqnarray}
  In addition to $\tilde{\sr{G}}_{0}(a,b)$, we will use the following sets for considering analytic regions (see \fig{region 1}).
\begin{eqnarray*}
 && \tilde{C}_{\delta}(u) = \{z \in \dd{C}; u - \delta < |z| < u+\delta, z \ne u \}, \qquad u, \delta > 0,\\
 && \tilde{\sr{G}}^{+}_{\delta}(u) = \tilde{\sr{G}}_{0}(u^{(1,\min)}_{1}, u) \cap \tilde{C}_{\delta}(u), \qquad u^{(1,\min)}_{1} < u,\\
 && \tilde{\sr{G}}^{-}_{\delta}(u) = \tilde{\sr{G}}_{0}(u,-u^{(1,\min)}_{1}) \cap \tilde{C}_{\delta}(u), \qquad u < -u^{(1,\min)}_{1}.
\end{eqnarray*}

\begin{remark} {\rm
\label{rem:extension 1}
  One may wonder whether \eqn{bound zeta} in \lem{extension 1} is sufficient for verifying analyticity of $\tilde{\varphi}_{1}(z)$ in $\tilde{\sr{G}}^{+}_{\delta}(u^{(1,\max)}_{1})$ when $\tilde{\tau}_{1} = u^{(1,\max)}_{1}$. This will turn out to have no problem because of \eqn{stationary equation 3}.
}\end{remark}

  In what follows, we first consider the case when (v-a) holds, then consider the other case.

\subsubsection{Singularity for the non-arithmetic case}
\label{sect:non-arithmetic}
  Assume the non-arithmetic condition (v-a). We consider the analytic behavior of $\tilde{\varphi}_{1}(z)$ around the singular point $z = \tilde{\tau}_{1}$. This behavior will show that there is no other singular point on $|z| = \tilde{\tau}_{1}$. We separately consider the three causes which are discussed above. 
 
\medskip
  
\noindent (\sect{Analytic function}a) The solution of \eqn{edge 1}: This equation has six solutions at most because it can be written as a polynomial equation with order six. $z=1, u^{(1,\edge)}_{1}$ are clearly the solutions. Because $ \tilde{\varphi}_{1}(z)$ of \eqn{stationary equation 3} is analytic for $|z| < \tilde{\tau}_{1}$, \eqn{edge 1} can not have solution such that $|z| < \tilde{\tau}_{1}$ except for the points where the numerator of the right hand side of \eqn{stationary equation 3} vanishes. This must be finitely many because the numerator vanishes otherwise by the uniqueness of analytic extension. On the other hand, \eqn{edge 1} has no solution on the circle $|z| = u^{(1,\edge)}_{1}$ except for $z = u^{(1,\edge)}_{1}$ by \lem{extension 1}. 

Thus, the compactness of the circle implies that, if $\tilde{\tau}_{1} = u^{(1,\edge)}_{1} < u^{(1,\max)}_{1}$, then \eqn{edge 1} has no solution on $\tilde{C}_{\delta}(u^{(1,\edge)}_{1})$ for some $\delta > 0$. Hence, we have the following fact from \eqn{stationary equation 3}.
\begin{lemma} {\rm
\label{lem:pole 1}
  Assume that $\tilde{\tau}_{1} = u^{(1,\edge)}_{1} < u^{(1,\max)}_{1}$ and $\tilde{\varphi}_{2}(\ul{\zeta}_{2}(z))$ is analytic at $|z| = u^{(1,\edge)}_{1}$. Then, $\tilde{\varphi}_{1}(z)$ has a simple pole at $z = u^{(1,\edge)}_{1}$, and analytic on $\tilde{C}_{\delta}(u^{(1,\edge)}_{1})$.
}\end{lemma}
\begin{remark} {\rm
\label{rem:}
  For categories \rmn{1} and \rmn{3}, the analytic condition on $\tilde{\varphi}_{2}(\ul{\zeta}_{2}(z))$ in this lemma is always satisfied because \lem{extension 1} and the category condition, $\ul{\zeta}_{2}(\tilde{\tau}_{1}) < \tilde{\tau}_{2}$, imply, for $|z| = u^{(1,\edge)}_{1}$,
\begin{eqnarray*}
  |\tilde{\varphi}_{2}(\ul{\zeta}_{2}(z))| \le \tilde{\varphi}_{2}(|\ul{\zeta}_{2}(z)|) \le \tilde{\varphi}_{2}(\ul{\zeta}_{2}(|z|)) = \tilde{\varphi}_{2}(u^{(1,\edge)}_{2}) < \infty.
\end{eqnarray*}
}\end{remark}

  If $\tilde{\tau}_{1} = u^{(1,\edge)}_{1} = u^{(1,\max)}_{1}$, then the analytic behavior of $\tilde{\varphi}_{1}(z)$ around $z=u^{(1,\edge)}_{1}$ is a bit complicated because $\ul{\zeta}_{2}(z)$ is also singular there. We will consider this case in \sectn{exact asymptotic v-a}.

\medskip

\noindent (\sect{Analytic function}b) The singularity of $\ul{\zeta}_{2}(z)$: By \lem{extension 1}, this function is analytic on \linebreak $\tilde{\sr{G}}_{0}(u^{(1,\min)}_{1}, u^{(1,\max)}_{1})$ and singular at $z = u^{(1,\max)}_{1}$, which is a branch point.

\medskip
  
\noindent (\sect{Analytic function}c) The singularity of $\tilde{\varphi}_{2}(\ul{\zeta}_{2}(z))$: This function is singular at $z = \tilde{\tau}_{1}$ if $\ul{\zeta}_{2}(\tilde{\tau}_{1}) = \tilde{\tau}_{2}$. Otherwise, it is singular at $z = u^{(1,\max)}_{1}$ because $\ul{\zeta}_{2}(z)$ is singular there. Furthermore, we may simultaneously have $\ul{\zeta}_{2}(\tilde{\tau}_{1}) = \tilde{\tau}_{2}$ and $\tilde{\tau}_{1} = u^{(1,\max)}_{1}$. Thus, we need to consider these three cases: $\tilde{\tau}_{1} = u^{(1,\max)}_{1}$ for categories \rmn{1} and \rmn{3}, and $\tilde{\tau}_{1} < u^{(1,\max)}_{1}$ or $\tilde{\tau}_{1} = u^{(1,\max)}_{1}$ for category \rmn{2}. For this, we will use the following fact, which is essentially the same as Lemma 4.2 of \cite{MR 2009}.
\begin{lemma} {\rm
\label{lem:zeta 2 max}
  $\ol{\zeta}_{1}(e^{\theta})$ is a concave function of $\theta \in [\theta^{(2,\min)}_{2}, \theta^{(2,\max)}_{2}]$, $\ol{\zeta}_{1}'(u^{(1,\max)}_{2}) = 0$, $\ol{\zeta}_{1}''(u^{(1,\max)}_{2}) < 0$, and
\begin{eqnarray}
\label{eqn:zeta 2 max}
  \hspace{-10ex}\lim_{z \to u^{(1,\max)}_{1} \atop z \in \tilde{\sr{G}}_{0}(u^{(1,\min)}_{1}, u^{(1,\max)}_{1})} \frac {u^{(1,\max)}_{2} - \ul{\zeta}_{2}(z)} {(u^{(1,\max)}_{1} - z)^{\frac 12}} = \frac{\sqrt{2}} {\sqrt{-\ol{\zeta}_{1}''(u^{(1,\max)}_{2})}}.
\end{eqnarray}
\begin{proof}
  The first part is immediate from the facts that $\Gamma_{+}$ is a convex set and $u^{(1,\max)}_{1} = e^{\theta^{(1,\max)}_{1}}$. By Taylor expansion of $\ol{\zeta}_{1}(z_{2})$ at $z_{2} = u^{(1,\max)}_{2} < u^{(2,\max)}_{2}$,
\begin{eqnarray*}
  \ol{\zeta}_{1}(z_{2}) = u^{(1,\max)}_{1} + \frac 12 \ol{\zeta}_{1}''(u^{(1,\max)}_{2}) (z_{2} - u^{(1,\max)}_{2})^{2} + o(|z_{2} - u^{(1,\max)}_{2}|^{2}).
\end{eqnarray*}
  Letting $z_{2} = \ul{\zeta}_{2}(z)$ in this equation yields \eqn{zeta 2 max} since $\ol{\zeta}_{1}(\ul{\zeta}_{2}(z)) = z$ for $z$ to be sufficiently close to $u^{(1,\max)}_{1}$. \pend
\end{proof}
}\end{lemma}

Another useful asymptotic is:
\begin{lemma} {\rm
\label{lem:gamma 1 max}
  If $u^{(1,\max)}_{1} = u^{(1,\edge)}_{1}$, then for any $\delta > 0$,
\begin{eqnarray}
\label{eqn:gamma 1 max}
  \lim_{z \to u^{(1,\max)}_{1} \atop z \in \tilde{\sr{G}}^{+}_{\delta}(u^{(1,\max)}_{1})} \frac {(u^{(1,\max)}_{1} - z)^{\frac 12}} {1 - \tilde{\gamma}_{1}(z, \ul{\zeta}_{2}(z))} = \frac {\sqrt{-\ol{\zeta}_{1}''(u^{(1,\max)}_{2})}} {\sqrt{2} p^{(1)}_{*1}(u^{(1,\edge)}_{1})}.
\end{eqnarray}
\begin{proof}
By the condition $u^{(1,\max)}_{1} = u^{(1,\edge)}_{1}$, we have
\begin{eqnarray*}
  \lefteqn{1 - \tilde{\gamma}_{1}(z, \ul{\zeta}_{2}(z)) = \tilde{\gamma}_{1}(\vc{u}^{(1,\max)}) - \tilde{\gamma}_{1}(z, \ul{\zeta}_{2}(z))} \hspace{10ex}\\
  && = (u^{(1,\max)}_{2} - \ul{\zeta}_{2}(z)) p^{(1)}_{*1}(u^{(1,\edge)}_{1})\\
  && \qquad  + \ul{\zeta}_{2}(z) (p^{(1)}_{*1}(u^{(1,\edge)}_{1}) - p^{(1)}_{*1}(z)) + p^{(1)}_{*0}(u^{(1,\edge)}_{1}) - p^{(1)}_{*0}(z).
\end{eqnarray*}
  Hence, by dividing both sides by $(u^{(1,\max)}_{1} - z)^{\frac 12}$, \lem{zeta 2 max} yields \eqn{zeta 2 max} because $p^{(1)}_{*1}(z)$ and $p^{(1)}_{*0}(z)$ are analytic except for $z = 0$. 
\pend \end{proof}
 }\end{lemma}
 
 We now consider the three cases separately.

\medskip
  
\noindent (\sect{Analytic function}c-1) $\ul{\zeta}_{2}(\tilde{\tau}_{1}) < \tilde{\tau}_{2}$, equivalently, categories \rmn{1} or \rmn{3}, and $\tilde{\tau}_{1} = u^{(1,\max)}_{1}$: In this case, $\tilde{\varphi}_{2}(z)$ is analytic for $z \in \tilde{C}_{\delta}(u^{(1,\max)}_{2})$ for some $\delta > 0$ because $u^{(1,\max)}_{2} = \ul{\zeta}_{2}(u^{(1,\max)}_{1}) = \ul{\zeta}_{2}(\tilde{\tau}_{1}) < \tilde{\tau}_{2}$. Hence, by Taylor expansion, we have, for $|z| < \tilde{\tau}_{2}$,
\begin{eqnarray}
\label{eqn:Taylor 1}
  \tilde{\varphi}_{2}(z) = \tilde{\varphi}_{2}(u^{(1,\max)}_{2}) + \tilde{\varphi}_{2}'(u^{(1,\max)}_{2}) (z - u^{(1,\max)}_{2}) + o(|z - u^{(1,\max)}_{2}|).
\end{eqnarray}
  Thus the analytic behavior of $\tilde{\varphi}_{2}(\ul{\zeta}_{2}(z))$ around $z = u^{(1,\max)}_{1}$ is determined by that of $\ul{\zeta}_{2}(z) - u^{(1,\max)}_{2}$. Since $u^{(1,\max)}_{2} = \ul{\zeta}_{2}(u^{(1,\max)}_{1}) < \tilde{\tau}_{2}$ by the conditions of (\sect{Analytic function}c-1), \lem{zeta 2 max} yields
\begin{eqnarray}
\label{eqn:Branch 1}
  \lefteqn{\hspace{-5ex} \tilde{\varphi}_{2}(\ul{\zeta}_{2}(z)) = \tilde{\varphi}_{2}(u^{(1,\max)}_{2}) - \frac{\sqrt{2} \tilde{\varphi}_{2}'(u^{(1,\max)}_{2}) } {\sqrt{-\ol{\zeta}_{1}''(u^{(1,\max)}_{2})}} (u^{(1,\max)}_{1} - z)^{\frac 12}} \hspace{30ex} \nonumber\\
  && + o(|z - u^{(1,\max)}_{2}|^{\frac 12}).
\end{eqnarray}
Thus, $\tilde{\varphi}_{2}(\ul{\zeta}_{2}(z))$ has a branch point of order 2 at $z = \tilde{\tau}_{1} = u^{(1,\max)}_{1}$, and is analytic on $\tilde{\sr{G}}^{+}_{\delta}(u^{(1,\max)}_{1})$ for some $\delta > 0$.

\medskip

\noindent (\sect{Analytic function}c-2) $\ul{\zeta}_{2}(\tilde{\tau}_{1}) = \tilde{\tau}_{2}$ and $\tilde{\tau}_{1} < u^{(1,\max)}_{1}$: This is only for category \rmn{2}. Hence, $\tilde{\tau}_{2} = u^{(2,\edge)}_{2} < u^{(2,\max)}_{2}$, and therefore $\tilde{\varphi}_{2}$-version of \lem{pole 1} is available. Thus, $\tilde{\varphi}_{2}(z)$ has a simple pole at $z = u^{(2,\edge)}_{2}$. Here, $u^{(2,\edge)}_{2}$ is the solution of the equation:
\begin{eqnarray}
\label{eqn:edge 2}
  \tilde{\gamma}_{2}(\ul{\zeta}_{1}(z), z) = 1
\end{eqnarray}
is crucial. Furthermore, $\ul{\zeta}_{2}(z)$ is analytic at $z = \tilde{\tau}_{1}$. Hence, $\tilde{\varphi}_{2}(\ul{\zeta}_{2}(z))$ has a simple pole at $z = \tilde{\tau}_{1}$, and is analytic on $\tilde{C}_{\delta}(u^{(1,\max)}_{1})$ for some $\delta > 0$.

\medskip

\noindent (\sect{Analytic function}c-3) $\ul{\zeta}_{2}(\tilde{\tau}_{1}) = \tilde{\tau}_{2}$, and $\tilde{\tau}_{1} = u^{(1,\max)}_{1}$: This is also only for category \rmn{2}. This case is similar to (\sect{Analytic function}c-2) except that $\ul{\zeta}_{2}(z)$ has a branch point at $z = \tilde{\tau}_{1} = u^{(1,\max)}_{1}$. Since $\tilde{\varphi}_{2}(z)$ has a simple pole at $z = \tilde{\tau}_{2}$, we have, by \lem{zeta 2 max},
\begin{eqnarray*}
 \tilde{\varphi}_{2}(\ul{\zeta}_{2}(z)) \sim (u^{(1,\max)}_{1} - z)^{-\frac 12},
\end{eqnarray*}
  and $\tilde{\varphi}_{2}(\ul{\zeta}_{2}(z))$ is analytic on $\tilde{\sr{G}}^{+}_{\delta}(u^{(1,\max)}_{1})$ for some $\delta > 0$.
  
\subsubsection{Singularity for the arithmetic case}
\label{sect:arithmetic case}

  We next consider the case that (v-a) does not holds. That is, the Markov additive process for the interior is arithmetic. In this case, the singularity of $\tilde{\varphi}_{1}(z)$ at $z = \tilde{\tau}_{1}$ occurs similarly to in \sectn{non-arithmetic}. In addition to this singular point, we may have another singular point $- \tilde{\tau}_{1}$ as can be seen in \tab{period}. For this, we separately consider two sub-cases:
\begin{mylist}{3}
\item [(B1)] either (v-b) or $m^{(1)}_{2} = 0$ holds. \hspace{3ex} (B2) $\,$ neither (v-b) nor $m^{(1)}_{2} = 0$ holds.
\end{mylist}
 In some cases, we need further classification: 
\begin{mylist}{3}
\item [(C1)] either (v-c) or $m^{(2)}_{1} = 0$ holds. \hspace{3ex} (C2) $\,$ neither (v-b) nor $m^{(2)}_{1} = 0$ holds.
\end{mylist}

  Consider (B1). From \tab{period}, the solutions of \eqn{D root} are $z=\pm u^{(1,\max)}_{1}$, and the solution of \eqn{no solution 1} is $z= u^{(1,\edge)}_{1}$. There is no other solution. We consider cases similar to (\sect{Analytic function}a), (\sect{Analytic function}b), (\sect{Analytic function}c-2), (\sect{Analytic function}c-1) and (\sect{Analytic function}c-3) of \sectn{non-arithmetic}.

\begin{mylist}{3}
\item [(\sect{Analytic function}a')] The solution of \eqn{edge 1}: This case is exactly the same as in \sectn{non-arithmetic} because $z= - u^{(1,\edge)}_{1}$ is not the solution of \eqn{no solution 1}. Hence, \lem{pole 1} also holds true.

\medskip

\item [(\sect{Analytic function}b')] The singularity of $\ul{\zeta}_{2}(z)$ at $|z| = u^{(1,\max)}_{1}$:  It is singular at $z=\pm u^{(1,\max)}_{1}$.

\medskip

\item [(\sect{Analytic function}c')] The singularity of $\tilde{\varphi}_{2}(\ul{\zeta}_{2}(z))$ at $|z| = \tilde{\tau}_{1}$: For $z = \tilde{\tau}_{1}$, the story is the same as in \sectn{non-arithmetic}. Hence, we only consider the case that $z = -\tilde{\tau}_{1}$. From \eqn{ul zeta 2} and the condition that (v-a) does not hold, we have
\begin{eqnarray}
\label{eqn:ul zeta 3}
  \ul{\zeta}_{2}(-\tilde{\tau}_{1}) = - \frac {1 - p_{00}} {2 (p_{-11} + p_{11} \tilde{\tau}_{1}^{2})} \tilde{\tau}_{1} = - \ul{\zeta}_{2}(\tilde{\tau}_{1}),
\end{eqnarray}
Hence, $|\ul{\zeta}_{2}(-\tilde{\tau}_{1})| = \ul{\zeta}_{2}(\tilde{\tau}_{1}) > 0$, and
\begin{eqnarray*}
  |\ul{\zeta}_{1}(\ul{\zeta}_{2}(-\tilde{\tau}_{1}))| = |\ul{\zeta}_{1}(-\ul{\zeta}_{2}(\tilde{\tau}_{1}))| = \ul{\zeta}_{1}(\ul{\zeta}_{2}(\tilde{\tau}_{1})).
\end{eqnarray*}
  Since $\ul{\zeta}_{1}(\ul{\zeta}_{2}(\tilde{\tau}_{1}))  < \tilde{\tau}_{1}$, $\tilde{\varphi}_{1}(\ul{\zeta}_{1}(\ul{\zeta}_{2}(z)))$ is analytic around $z = -\tilde{\tau}_{1}$. Furthermore, \lem{zeta 2 max} and \eqn{Branch 1} are still valid if we replace $u^{(1,\max)}_{i}$ by $- u^{(1,\max)}_{i}$ for $i=1,2$. However, this $z= - \tilde{\tau}_{1}$ can not be the solution of \eqn{edge 1} because of (B1). Thus, we have to partially change the arguments in \sectn{non-arithmetic}. 
  
\begin{mylist}{0}
\item [(\sect{Analytic function}c'-1)] $\ul{\zeta}_{2}(\tilde{\tau}_{1}) < \tilde{\tau}_{2}$ and $\tilde{\tau}_{1} = u^{(1,\max)}_{1}$: This is only for categories \rmn{1} and \rmn{3}, and $\tilde{\varphi}_{2}(\ul{\zeta}_{2}(z))$ has a branch point of order 2 at $z = - u^{(1,\max)}_{1}$, and is analytic on $\tilde{\sr{G}}^{-}_{\delta}(-u^{(1,\max)}_{1}) \cap \tilde{\sr{G}}^{+}_{\delta}(u^{(1,\max)}_{1})$ for some $\delta > 0$ because it has also a branch point at $z = u^{(1,\max)}_{1}$.

\medskip

\item [(\sect{Analytic function}c'-2)] $\ul{\zeta}_{2}(\tilde{\tau}_{1}) = \tilde{\tau}_{2}$ and $\tilde{\tau}_{1} < u^{(1,\max)}_{1}$: This is only for category \rmn{2}. Since $\ul{\zeta}_{2}(z)$ is analytic at $z = \tilde{\tau}_{1}$, $\tilde{\varphi}_{2}(\ul{\zeta}_{2}(z))$ is analytic at $z = - \tilde{\tau}_{1}$ if (C1) holds. Otherwise, if (C2) holds, it has a simple pole at $z = - \tilde{\tau}_{1}$ because $\ul{\zeta}_{2}(-\tilde{\tau}_{1}) = - \ul{\zeta}_{2}(\tilde{\tau}_{1})$ is the solution of \eqn{edge 2}.

\medskip

\item [(\sect{Analytic function}c'-3)] $\ul{\zeta}_{2}(\tilde{\tau}_{1}) = \tilde{\tau}_{2}$ and $\tilde{\tau}_{1} = u^{(1,\max)}_{1}$: This is only for category \rmn{2}, and the situation is similar to (\sect{Analytic function}c'-2) except that the singularity is caused by $\ul{\zeta}_{2}(z)$ at $z = - \tilde{\tau}_{1}$. To verify this fact, we rework on $\tilde{\varphi}_{2}(\ul{\zeta}_{2}(z))$. Similarly to \eqn{stationary equation 3}, we have, for $|z| \in (u^{2,\min)}_{2}, \tilde{\tau}_{2})$,
\begin{eqnarray*}
 \tilde{\varphi}_{2}(z) = \frac{(\tilde{\gamma}_{1}(\ul{\zeta}_{1}(z),z) - 1)  \tilde{\varphi}_{1}(\ul{\zeta}_{1}(z))}{1 -  \tilde{\gamma}_{2}(\ul{\zeta}_{1}(z),z)} + \frac{ (\tilde{\gamma}_{0}(\ul{\zeta}_{2}(z), z) - 1) \nu(\vc{0})} {1 -  \tilde{\gamma}_{2}(\ul{\zeta}_{1}(z),z)}.
\end{eqnarray*}
  Substituting $\ul{\zeta}_{2}(z)$ into $z$ of this equation, we have
\begin{eqnarray}
\label{eqn:stationary equation 4}
 \lefteqn{\hspace{-8ex}\tilde{\varphi}_{2}(\ul{\zeta}_{2}(z)) = \frac{(\tilde{\gamma}_{1}(\ul{\zeta}_{1}(\ul{\zeta}_{2}(z)), \ul{\zeta}_{2}(z)) - 1)  \tilde{\varphi}_{1}(\ul{\zeta}_{1}(\ul{\zeta}_{2}(z)))} {1 -  \tilde{\gamma}_{2}(\ul{\zeta}_{1}(\ul{\zeta}_{2}(z)), \ul{\zeta}_{2}(z))}} \nonumber \\
 && \hspace{10ex}+ \frac{(\tilde{\gamma}_{0}(\ul{\zeta}_{2}(\ul{\zeta}_{2}(z)), \ul{\zeta}_{2}(z)) - 1) \nu(\vc{0})} {1 -  \tilde{\gamma}_{2}(\ul{\zeta}_{1}(\ul{\zeta}_{2}(z)), \ul{\zeta}_{2}(z))}.
\end{eqnarray}
  By the assumptions of \noindent (\sect{Analytic function}c-3), if (C2) holds, then $\tilde{\varphi}_{2}(z)$ has a simple pole at $z = - \tilde{\tau}_{2}$, and therefore $\tilde{\varphi}_{2}(\ul{\zeta}_{2}(z)) \sim (-u^{(1,\max)}_{1} - z)^{-\frac 12}$ around $z = - u^{(1,\max)}_{1}$ by \lem{zeta 2 max}. Otherwise, if (C1) holds, we need to consider $\tilde{\varphi}_{1}(\ul{\zeta}_{1}(\ul{\zeta}_{2}(z)))$ in \eqn{stationary equation 4} due to the singularity of $\ul{\zeta}_{2}(z)$ at $z = - \tilde{\tau}_{1} = - u^{(1,\max)}_{1}$, where $\tilde{\varphi}_{1}(\ul{\zeta}_{1}(z))$ is analytic at $z = - u^{(1,\max)}_{2} = - \ul{\zeta}_{2}(u^{(1,\max)}_{1})$ because $$|\ul{\zeta}_{1}(\ul{\zeta}_{2}(-\tilde{\tau}_{1}))| = |\ul{\zeta}_{1}(\ul{\zeta}_{2}(\tilde{\tau}_{1}))| < \tilde{\tau}_{1}.$$
  Hence, $\tilde{\varphi}_{1}(\ul{\zeta}_{1}(\ul{\zeta}_{2}(z))) - \tilde{\varphi}_{1}(-\ul{\zeta}_{1}(u^{(1,\max)}_{2})) \sim (-u^{(1,\max)}_{1} - z)^{\frac 12}$. On the other hand, $\tilde{\gamma}_{1}(\ul{\zeta}_{1}(\ul{\zeta}_{2}(z)), \ul{\zeta}_{2}(z)) - 1 \sim (-u^{(1,\max)}_{1} - z)^{\frac 12}$ because (v-a) does not hold. Combining these asymptotics in \eqn{stationary equation 4}, we have $\tilde{\varphi}_{2}(\ul{\zeta}_{2}(z)) - \tilde{\varphi}_{2}(-\tilde{\tau}_{2}) \sim (-u^{(1,\max)}_{1} - z)^{\frac 12}$ around $z = - u^{(1,\max)}_{1}$ by \lem{zeta 2 max}.
\end{mylist}
\end{mylist}

We next consider (B2). From \tab{period}, the solutions of \eqn{D root} are $z=\pm u^{(1,\max)}_{1}$, and the solutions of \eqn{no solution 1} are $z=\pm u^{(1,\edge)}_{1}$. In this case, the arguments for $z = - \tilde{\tau}_{1}$ is completely parallel to those for $z = \tilde{\tau}_{1}$ except for the cases (\sect{Analytic function}c'-2) and (\sect{Analytic function}c'-3). The latter two cases are also parallel if (C2) holds. Otherwise, $\tilde{\varphi}_{2}(z)$ is analytic at $z = - \tilde{\tau}_{2}$.

\subsection{Asymptotic inversion formula}
\label{sect:Asymptotic inversion}

  From these singularities, we derive exact tail asymptotics of the stationary distribution. For this, we use Tauberian type theorem for generating functions.

\begin{lemma}[Theorem VI.5 of \cite{Flajolet-Sedqewick 2009}] {\rm
\label{lem:G-asymptotic}
  Let $f$ be a generating function of a sequence of real numbers $\{p(n); n=0,1,\ldots\}$. If $f(z)$ is singular at finitely many points $a_{1}, a_{2}, \ldots, a_{m}$ on the circle $|z| = \rho$ for some $\rho > 0$ and positive integer $m$, and analytic on the set
\begin{eqnarray*}
  \Delta_{i} \equiv \{ z \in \dd{C}; |z| < r_{i}, z \ne a_{i}, |\arg(z - a_{i})| > \omega_{i} \}, \quad i=1,2,\ldots, m,
\end{eqnarray*}
   for some $\omega_{i}$ and $r_{i}$ such that $\rho < r_{i}$ and $0 \le \omega_{i} < \frac \pi 2$ and if
\begin{eqnarray}
\label{eqn:mgf a}
  \lim_{\Delta_{i} \ni z \to a_{i}} (a_{i} - z)^{\kappa_{i}} f(z) = b_{i}, \qquad i=1,2,\ldots, m,
\end{eqnarray}
  for $\kappa_{i} \not\in \{0, -1, -2, \ldots\}$ and some constant $b_{i} \in \dd{R}$, then 
\begin{eqnarray}
\label{eqn:fn a}
  \lim_{n \to \infty} \left( \sum_{i=1}^{m} \frac {n^{\kappa_{i}-1}}{\Gamma(\kappa_{i})} a_{i}^{-n}\right)^{-1} p(n) = b,
\end{eqnarray}
  for some real number $b$, where $\Gamma(z)$ is the gamma function for complex number $z$ (see Sec 52 of Volume II of \cite{Markushevich 1977}).
}\end{lemma}

  Recall that the asymptotic notation ``$\sim$'' introduced in \sectn{Introduction}. With this notation, \eqn{fn a} can be written as
\begin{eqnarray*}
  p(n) \sim \sum_{i=1}^{m} \frac {n^{\kappa_{i}-1}}{\Gamma(\kappa_{i})} a_{i}^{-n},
\end{eqnarray*}
  where $\Gamma(\frac 12) = \sqrt{\pi}$ and $\Gamma(- \frac 12) = - 2 \sqrt{\pi}$.

  We will apply \lem{G-asymptotic} in the following cases: For $m=1$, $a_{1} = u^{(1,\edge)}_{1}$ and $\kappa_{1} = 1,2$, $a_{1} = u^{(1,\max)}_{1}$ and $\kappa_{1} = \pm \frac 12$. For $m=2$, $a_{1} = \pm u^{(1,\edge)}_{1}$ and $\kappa_{1} = 1,2$, $a_{1} = \pm u^{(1,\max)}_{1}$ and $\kappa_{1} = - \frac 12$.

\section{Exact tail asymptotics for the non-arithmetic case} 
\label{sect:exact asymptotic v-a}

  Throughout this section, we assume the non-arithmetic condition (v-a). We first derive exact asymptotics for the stationary probabilities $\nu(n,0)$ and $\nu(0,n)$ on the boundary faces. Because of symmetry, we are only concerned with $\nu(n,0)$.  

\subsection{The boundary probabilities for the non-arithmetic case}
\label{sect:boundary probabilities v-a}

  We separately consider the two cases that $u^{(1,\cp)}_{2} < u^{(2,\cp)}_{2}$ and $u^{(1,\cp)}_{2} \ge u^{(2,\cp)}_{2}$, which correspond with categories \rmn{1} (or \rmn{3}) and \rmn{2}, respectively. In this subsection, we prove the following two theorems.
  
\begin{theorem} {\rm
\label{thr:exact nu 10n-1}
  Under the conditions (i)--(iv) and (v-a), for categories \rmn{1} and \rmn{3}, $\tilde{\tau}_{1} = u^{(1,\cp)}_{1}$, and $P(L_{1} = n, L_{2} = 0)$ has the following exact asymptotic $h_{1}(n)$.
\begin{eqnarray}
\label{eqn:exact nu 10n-1}
  h_{1}(n) = \left\{\begin{array}{ll}
  \tilde{\tau}_{1}^{-n}, \quad & u^{(1,\cp)}_{1} \not= u^{(1,\max)}_{1},\\
  n^{- \frac 12} \tilde{\tau}_{1}^{-n}, \quad & u^{(1,\cp)}_{1} = u^{(1,\max)}_{1} = u^{(1,\edge)}_{1},\\
  n^{- \frac 32} \tilde{\tau}_{1}^{-n}, \quad & u^{(1,\cp)}_{1} = u^{(1,\max)}_{1} \ne u^{(1,\edge)}_{1}.
  \end{array} \right.
\end{eqnarray}
  By symmetry, the corresponding results are also obtained for $P(L_{1} = 0, L_{2} = n)$ for categories \rmn{1} and \rmn{2}.
}\end{theorem}
  
\begin{theorem} {\rm
\label{thr:exact nu 10n-2}
  Under the conditions (i)--(iv) and (v-a), for category \rmn{2}, $\tilde{\tau}_{2} = u^{(2,\edge)}_{2}$, and $P(L_{1} = n, L_{2} = 0)$ has the following exact asymptotic $h_{1}(n)$.
\begin{eqnarray}
\label{eqn:exact nu 10n-2}
  h_{1}(n) = \left\{\begin{array}{ll}
  \tilde{\tau}_{1}^{-n}, \quad & \tilde{\tau}_{1} < u^{(1,\cp)}_{1}, \mbox{ or }\\
  & \tilde{\tau}_{1} = u^{(1,\cp)}_{1} = u^{(1,\max)}_{1} = u^{(1,\edge)}_{1}, \\
  n \tilde{\tau}_{1}^{-n}, \quad & \tilde{\tau}_{1} = u^{(1,\cp)}_{1} \not= u^{(1,\max)}_{1},\\
  n^{- \frac 12} \tilde{\tau}_{1}^{-n}, \; & \tilde{\tau}_{1} = u^{(1,\cp)}_{1} = u^{(1,\max)}_{1} \ne u^{(1,\edge)}_{1}.  \end{array} \right.
\end{eqnarray}
  By symmetry, the corresponding results are also obtained for $P(L_{1} = 0, L_{2} = n)$ for categories \rmn{3}.
}\end{theorem}

\begin{remark} {\rm
\label{rem:exact nu 10n-2 1}
  Theorems \thrt{exact nu 10n-1} and \thrt{exact nu 10n-2} are exactly corresponds with Theorem 6.1 of \cite{Dai-Miyazawa 2011b} (see also Theorems 2.1 and 2.3 of \cite{Dai-Miyazawa 2011a}). This is not surprising because of the similarity of the stationary equations although moment generating functions are used in \cite{Dai-Miyazawa 2011a,Dai-Miyazawa 2011b}.
}\end{remark}

\begin{remark} {\rm
\label{rem:exact nu 10n-2 2}
  These theorems fill missing cases for the exact asymptotics of Theorem 4.2 of \cite{Miyazawa 2009}. Furthermore, they correct two errors there. Both of them are for category \rmn{2}. The exact asymptotic is geometric for $\tilde{\tau}_{1} = u^{(1,\cp)}_{1} = u^{(1,\max)}_{1} = u^{(1,\edge)}_{1}$, and not geometric for $\tilde{\tau}_{1} = u^{(1,\cp)}_{1} \not= u^{(1,\max)}_{1}$ (see \thr{exact nu 10n-2}). However, in Theorem 4.2 of \cite{Miyazawa 2009}, they are not geometric (see (43d3) there) and geometric (see (4c) there), respectively. Thus, these should be corrected.
}\end{remark}

\begin{proof*}{Proof of \thr{exact nu 10n-1}}
  We assume that  either in category \rmn{1} or \rmn{3} occurs. This is equivalent to $u^{(1,\cp)}_{2} < u^{(2,\cp)}_{2}$, and $\tilde{\tau}_{1} = u^{(1,\cp)}_{1}$. Furthermore, we always have $\ul{\zeta}_{2}(u^{(1,\cp)}_{1}) = u^{(1,\cp)}_{2} < \tilde{\tau}_{2}$, and therefore $\tilde{\varphi}_{2}(\ul{\zeta}_{2}(z))$ is analytic at $z = u^{(1,\cp)}_{1}$. We consider three cases separately.
  
\noindent (\sect{exact asymptotic v-a}a)  $u^{(1,\cp)}_{1} < u^{(1,\max)}_{1}$: This case implies that $u^{(1,\edge)}_{1} < u^{(1,\max)}_{1}$ and $\tilde{\gamma}_{1}(\vc{u}^{(1,\max)}) > 1$, and therefore $\vc{u}^{(1,\edge)} = \vc{u}^{(1,\cp)}$. Hence, by \lem{pole 1}, $\tilde{\varphi}_{1}$ of \eqn{stationary equation 3} satisfies the conditions of \lem{G-asymptotic} under the setting \eqn{mgf a} with $a_{1} = u^{(1,\edge)}_{1}$, $\kappa = 1$. Thus, letting
\begin{eqnarray*}
  b = \frac{(\tilde{\gamma}_{2}(\vc{u}^{(1,\edge)}) - 1)  \tilde{\varphi}_{2}(u^{(1,\edge)}_{2}) + (\tilde{\gamma}_{0}(\vc{u}^{(1,\edge)}) - 1) \nu(\vc{0})} {\frac {d}{du} \tilde{\gamma}_{1}(u, \ul{\zeta}_{2}(u))|_{u = u^{(1,\edge)}_{1}}},
\end{eqnarray*}
  which must be positive by \eqn{fn a} and the fact that $\tilde{\varphi}_{1}(z)$ is singular at $z = u^{(1,\edge)}_{1}$, we have
\begin{eqnarray*}
  \lim_{n \to \infty} \tilde{\tau}_{1}^{n} P(L_{1} = n, L_{2} = 0) = b.
\end{eqnarray*}

\noindent (\sect{exact asymptotic v-a}b) $u^{(1,\cp)}_{1} = u^{(1,\max)}_{1}$, $u^{(1,\edge)}_{1} = u^{(1,\max)}_{1}$: In this case, category \rmn{3} is impossible, and $\tilde{\gamma}_{1}(\vc{u}^{(1,\max)}) = 1$. On the other hand, $\tilde{\varphi}_{2}(z)$ is analytic at $z = \ul{\zeta}_{2}(u^{(1,\max)}_{1}) < \tilde{\tau}_{2}$ because of Category \rmn{1}. Hence, we can use the Taylor expansion \eqn{Taylor 1}, and therefore \eqn{stationary equation 3}, \eqn{Branch 1} and \lem{gamma 1 max} yield, for some $\delta > 0$,
\begin{eqnarray}
\label{eqn:null asymptotic 1}
  \lim_{\tilde{\sr{G}}^{+}_{\delta}(u^{(1,\max)}_{1}) \ni z \to u^{(1,\max)}_{1}} (u^{(1,\max)}_{1} - z)^{\frac 12} \tilde{\varphi}_{1}(z) = b,
\end{eqnarray}
where  
\begin{eqnarray*}
  b = \Big( (\tilde{\gamma}_{2}(\vc{u}^{(1,\max)}) - 1)  \tilde{\varphi}_{2}(u^{(1,\max)}_{2}) + (\tilde{\gamma}_{0}(\vc{u}^{(1,\max)}) - 1) \nu(\vc{0}) \Big) \frac {\sqrt{-\ol{\zeta}_{1}''(u^{(1,\max)}_{2})}} {\sqrt{2} p^{(1)}_{*1}(u^{(1,\edge)}_{1})}.
\end{eqnarray*}
  Hence, $\tilde{\varphi}_{1}$ satisfies the conditions of \lem{G-asymptotic} under the setting \eqn{mgf a} with $a_{1} = u^{(1,\max)}_{1}$ and $\kappa_{1} = \frac 12$, and therefore we have
\begin{eqnarray*}
  \lim_{n \to \infty} n^{\frac 12} \tilde{\tau}_{1}^{n} P(L_{1} = n, L_{2} = 0) = \frac b{\sqrt{\pi}},
\end{eqnarray*}
  where the positivity of $b$ is checked similarly to case (\sect{exact asymptotic v-a}a) (see also case (\sect{exact asymptotic v-a}c) below).
  
\medskip

\noindent (\sect{exact asymptotic v-a}c) $u^{(1,\cp)}_{1} = u^{(1,\max)}_{1}$, $u^{(1,\edge)}_{1} \not= u^{(1,\max)}_{1}$: In this case, category \rmn{3} is also impossible, and $\tilde{\gamma}_{1}(\vc{u}^{(1,\max)}) \not= 1$. Thus, we consider the setting \eqn{mgf a} with $\kappa_{1} = - \frac 12$. From \eqn{stationary equation 3}, we have
\begin{eqnarray}
\label{eqn:difference 1}
  \lefteqn{\tilde{\varphi}_{1}(z) - \tilde{\varphi}_{1}(u^{(1,\max)}_{1})}\nonumber \\
 &&  = \frac{(\tilde{\gamma}_{2}(z, \ul{\zeta}_{2}(z)) - 1)  \tilde{\varphi}_{2}(\ul{\zeta}_{2}(z)) + (\tilde{\gamma}_{0}(z, \ul{\zeta}_{2}(z)) - 1) \nu(\vc{0})} {1 -  \tilde{\gamma}_{1}(z, \ul{\zeta}_{2}(z))} - \tilde{\varphi}_{1}(u^{(1,\max)}_{1}) \nonumber \\
  && = \frac{(\tilde{\gamma}_{2}(z, \ul{\zeta}_{2}(z)) - 1) (\tilde{\varphi}_{2}(\ul{\zeta}_{2}(z)) - \tilde{\varphi}_{2}(u^{(1,\max)}_{2}))} {1 -  \tilde{\gamma}_{1}(z, \ul{\zeta}_{2}(z))} \nonumber \\
 && \quad + \frac{(\tilde{\gamma}_{2}(z, \ul{\zeta}_{2}(z)) - \tilde{\gamma}_{2}(\vc{u}^{(1,\max)})) \tilde{\varphi}_{2}(u^{(1,\max)}_{2})}{1 -  \tilde{\gamma}_{1}(z, \ul{\zeta}_{2}(z))} + \frac{(\tilde{\gamma}_{0}(z, \ul{\zeta}_{2}(z)) - \tilde{\gamma}_{0}(\vc{u}^{(1,\max)})) \nu(\vc{0})} {1 -  \tilde{\gamma}_{1}(z, \ul{\zeta}_{2}(z))} \nonumber \\
  && \quad + \Bigg[ \frac {(\tilde{\gamma}_{2}(\vc{u}^{(1,\max)}) - 1) \tilde{\varphi}_{2}(u^{(1,\max)}_{2})}{(1 -  \tilde{\gamma}_{1}(z, \ul{\zeta}_{2}(z))) (1 -  \tilde{\gamma}_{1}(\vc{u}^{(1,\max)}))} + \frac{(\tilde{\gamma}_{0}(\vc{u}^{(1,\max)}) - 1) \nu(\vc{0})} {(1 -  \tilde{\gamma}_{1}(z, \ul{\zeta}_{2}(z))) (1 -  \tilde{\gamma}_{1}(\vc{u}^{(1,\max)}))} \Bigg] \nonumber \\
  && \hspace{5ex}\times \left( \tilde{\gamma}_{1}(z, \ul{\zeta}_{2}(z)) - \tilde{\gamma}_{1}(\vc{u}^{(1,\max)}) \right).
\end{eqnarray}
  We recall \eqn{Branch 1} that
\begin{eqnarray*}
  \tilde{\varphi}_{2}(\ul{\zeta}_{2}(z)) - \tilde{\varphi}_{2}(u^{(1,\max)}_{2}) = - (u^{(1,\max)}_{1} - z)^{\frac 12} \frac{\sqrt{2} \tilde{\varphi}_{2}'(u^{(1,\max)}_{2})} {\sqrt{-\ol{\zeta}_{1}''(u^{(1,\max)}_{2})}} + o(|u^{(1,\max)}_{1} - z|^{\frac 12}).
\end{eqnarray*}
  From \eqn{zeta 2 max}, we have
\begin{eqnarray*}
  \lefteqn{\tilde{\gamma}_{0}(z, \ul{\zeta}_{2}(z)) - \tilde{\gamma}_{0}(\vc{u}^{(1,\max)}) = (\ul{\zeta}_{2}(z) - \ul{\zeta}_{2}(u^{(1,\max)}_{1})) p^{(0)}_{*1}(z)} \hspace{10ex}\\
  && \quad + \ul{\zeta}_{2}(u^{(1,\max)}_{1}) (p^{(0)}_{*1}(z) - p^{(0)}_{*1}(u^{(1,\max)}_{1})) + p^{(0)}_{*0}(z) - p^{(0)}_{*0}(u^{(1,\max)}_{1})\\
  && = - \frac{\sqrt{2} p^{(0)}_{*1}(u^{(1,\max)}_{1})} {\sqrt{-\ol{\zeta}_{1}''(u^{(1,\max)}_{2})}} (u^{(1,\max)}_{1} - z)^{\frac 12} + o(|u^{(1,\max)}_{1} - z|^{\frac 12}).
\end{eqnarray*}
  Similarly,
\begin{eqnarray*}
 && \tilde{\gamma}_{1}(z, \ul{\zeta}_{2}(z)) - \tilde{\gamma}_{1}(\vc{u}^{(1,\max)})\\
 && \hspace{5ex} =  - \frac{\sqrt{2} p^{(1)}_{*1}(u^{(1,\max)}_{1})} {\sqrt{-\ol{\zeta}_{1}''(u^{(1,\max)}_{2})}} (u^{(1,\max)}_{1} - z)^{\frac 12} + o(|u^{(1,\max)}_{1} - z|^{\frac 12}),\\
 && \tilde{\gamma}_{2}(z, \ul{\zeta}_{2}(z)) - \tilde{\gamma}_{2}(\vc{u}^{(1,\max)})\\
 && \hspace{3ex}=  - \frac{\sqrt{2} \left(p^{(2)}_{*1}(u^{(1,\max)}_{1}) - \ds \frac {p^{(2)}_{*-1}(u^{(1,\max)}_{1})} {\big(u^{(1,\max)}_{2} \big)^{2}} \right)} {\sqrt{-\ol{\zeta}_{1}''(u^{(1,\max)}_{2})}} (u^{(1,\max)}_{1} - z)^{\frac 12} + o(|u^{(1,\max)}_{1} - z|^{\frac 12}).
\end{eqnarray*}

  With the following notation,
\begin{eqnarray*}
  && c_{1} = \frac{\sqrt{2}} {\big(1 -  \tilde{\gamma}_{1}(\vc{u}^{(1,\max)})\big) \sqrt{-\ol{\zeta}_{1}''(u^{(1,\max)}_{2})}}, \\
  && d_{k} = \left. \frac {\partial}{\partial v} \tilde{\gamma}_{k}(u^{(1,\max)}_{1}, v) \right|_{v=\ul{\zeta}_{2}(u^{(1,\max)}_{1})},
\end{eqnarray*}
  \eqn{difference 1} yields, as $z \to u^{(1,\max)}_{1}$ satisfying that $z \in \tilde{\sr{G}}^{+}_{\delta}(u^{(1,\max)}_{1})$ for some $\delta > 0$,
\begin{eqnarray}
\label{eqn:difference 2}
  \lefteqn{\hspace{-5ex}\tilde{\varphi}_{1}(z) - \tilde{\varphi}_{1}(u^{(1,\max)}_{1}) = - c_{1} (u^{(1,\max)}_{1} - z)^{\frac 12} \Big( \big( \tilde{\gamma}_{2}(\vc{u}^{(1,\max)}) - 1 \big) \tilde{\varphi}_{2}'(u^{(1,\max)}_{2})}\nonumber \\
  && \quad + d_{2} \tilde{\varphi}_{2}(u^{(1,\max)}_{2}) + d_{0} \nu(\vc{0}) + d_{1} \tilde{\varphi}_{1}(u^{(1,\max)}_{1}) \Big) + o(|u^{(1,\max)}_{1} - z|^{\frac 12}).
\end{eqnarray}
  Let
\begin{eqnarray*}
  b = - \Big( \big( \tilde{\gamma}_{2}(\vc{u}^{(1,\max)}) - 1 \big) \tilde{\varphi}_{2}'(u^{(1,\max)}_{2}) + d_{2} \tilde{\varphi}_{2}(u^{(1,\max)}_{2}) + d_{0} \nu(\vc{0}) + d_{1} \tilde{\varphi}_{1}(u^{(1,\max)}_{1}) \Big).
\end{eqnarray*}
  Then, taking $u_{1}$ which is sufficiently close $u^{(1,\max)}_{1}$ from below in \eqn{difference 2}, we can see that this $b$ must be negative because $\tilde{\varphi}_{1}(u_{1})$ is strictly increasing in $u_{1} \in [0,u^{(1,\max)}_{1})$. Thus, \eqn{mgf a} holds for the setting of \eqn{mgf a} with $\kappa_{1} = - \frac 12$, and therefore \eqn{fn a} leads to
\begin{eqnarray*}
    \lim_{n \to \infty} n^{\frac 32} \tilde{\tau}_{1}^{n} P(L_{1} = n, L_{2} = 0) = - \frac {b}{2\sqrt{\pi}} > 0.
\end{eqnarray*}
Thus, we have obtained all the cases of \eqn{exact nu 10n-1}, and the proof is completed.
\end{proof*}

\medskip

\begin{proof*}{Proof of \thr{exact nu 10n-2}}
  Assume category \rmn{2}. In this case, $\tilde{\tau}_{2} = \ul{\zeta}_{2}(\tilde{\tau}_{1})$, and $\tilde{\varphi}_{2}(z)$ has a simple pole at $z = \tilde{\tau}_{2}$ because of category \rmn{2} (see \noindent (\sect{Analytic function}c-2)). We need to consider the following cases.

\noindent (\sect{exact asymptotic v-a}a'): $\tilde{\tau}_{1} < u^{(1,\cp)}_{1}$: In this case, $\tilde{\varphi}_{2}(\ul{\zeta}_{2}(z))$ has a simple pole at $z = \tilde{\tau}_{1}$. Since $\tilde{\varphi}_{1}(z)$ has no other singularity on $|z| = \tilde{\tau}_{1}$. it has a simple pole at $z = \tilde{\tau}_{1}$.

\noindent (\sect{exact asymptotic v-a}b'): $\tilde{\tau}_{1} = u^{(1,\cp)}_{1}$: This case is further partitioned into the following subcases:
\begin{mylist}{3}
\item [(\sect{exact asymptotic v-a}b'-1)] $u^{(1,\cp)}_{1} \not= u^{(1,\max)}_{1}$: In this case, $\tilde{\tau}_{1} = u^{(1,\edge)}_{1} < u^{(1,\max)}_{1}$, and therefore it is easy to see from \eqn{stationary equation 3} that $\tilde{\varphi}_{1}(z)$ has a double pole at $z = \tilde{\tau}_{1}$. Hence, we can apply the setting \eqn{mgf a} with $a_{1} = u^{(1,\edge)}_{1}$ and $\kappa = 2$.
\item [(\sect{exact asymptotic v-a}b'-2)] $u^{(1,\cp)}_{1} = u^{(1,\max)}_{1} \ne u^{(1,\edge)}_{1}$: \eqn{edge 1} does not hold, and therefore \eqn{zeta 2 max} and the fact that $\tilde{\varphi}_{2}(z)$ has a simple pole at $z = \tilde{\tau}_{2}$ yield the same asymptotic as \eqn{null asymptotic 1} but with a different $b$. Hence, we apply \eqn{mgf a} with $\kappa_{1} = \frac 12$.
\item [(\sect{exact asymptotic v-a}b'-3)] $u^{(1,\cp)}_{1} = u^{(1,\max)}_{1} = u^{(1,\edge)}_{1}$: In this case, we note the following facts .
\begin{mylist}{0}
\item [(\sect{exact asymptotic v-a}b'-3-1)] $\tilde{\varphi}_{2}(z)$ has a simple pole at $z = \tilde{\tau}_{2}$, and therefore \lem{zeta 2 max} yield
\begin{eqnarray*}
  \tilde{\varphi}_{2}(\ul{\zeta}_{2}(z)) \sim (u^{(1,\max)}_{1} - z)^{- \frac 12}.
\end{eqnarray*}
\item [(\sect{exact asymptotic v-a}b'-3-2)] By \lem{gamma 1 max}, $1 - \tilde{\gamma}_{1}(z, \ul{\zeta}_{2}(z)) \sim (u^{(1,\max)}_{1} - z)^{\frac 12}$.
\end{mylist}
  Hence, \eqn{stationary equation 3} yields $\tilde{\varphi}_{1}(z) \sim (u^{(1,\max)}_{1} - z)^{-1}$, and therefore we apply \eqn{mgf a} with $a_{1} = u^{(1,\edge)}_{1}$ and $\kappa = 1$.
\end{mylist}
  
  Thus, similar to \thr{exact nu 10n-1}, we can obtain \eqn{exact nu 10n-2}, which completes the proof.
\end{proof*}
  
\subsection{The marginal distributions for the non-arithmetic case} 
\label{sect:marginal distributions v-a}

We consider the asymptotics of $P(\br{\vc{c}, \vc{L}} \ge x)$ as $x \to \infty$ for $\vc{c} = (1,0), (0,1), \linebreak(1,1)$. For them, we use the generating functions  $\tilde{\varphi}_{+}(z,1)$,  $\tilde{\varphi}_{+}(1,z)$ and $ \tilde{\varphi}_{+}(z, z)$. For simplicity, we denote them by $\psi_{10}(z)$, $\psi_{01}(z)$, $\psi_{11}(z)$, respectively. We note that generating functions are not useful for the other direction $\vc{c}$ because we can not appropriately invert them. For general $\vc{c} > 0$, we should use moment generating functions instead of generating functions. However, in this case, we need finer analytic properties to apply asymptotic inversion (e.g., see Appendix C of \cite{Dai-Miyazawa 2011a}). Thus, we leave it for future study.

From \eqn{stationary equation 2} and \eqn{decomposition 1}, we have, for $\vc{z} \in \dd{C}^{2}$ satisfying $(|z_{1}|,|z_{2}|) \in \tilde{\sr{D}}$,
\begin{eqnarray}
\label{eqn:stationary equation 6}
  \tilde{\varphi}(\vc{z}) = \left(1+\frac{\tilde{\gamma}_{1}(\vc{z}) - 1} {1 - \tilde{\gamma}_{+}(\vc{z}) }\right) \tilde{\varphi}_{1}(z_{1}) + \frac{\tilde{\gamma}_{2}(\vc{z}) - \tilde{\gamma}_{+}(\vc{z})} {1 - \tilde{\gamma}_{+}(\vc{z})} \tilde{\varphi}_{2}(z_{2}) + \frac{\tilde{\gamma}_{0}(\vc{z}) - \tilde{\gamma}_{+}(\vc{z})} {1 - \tilde{\gamma}_{+}(\vc{z}) } \nu(\vc{0}). \qquad
\end{eqnarray}
  Hence,  the asymptotics of $P(\br{\vc{c}, \vc{L}} \ge x)$ can be obtained for $\vc{c} = (1,0), (0,1), (1,1)$ by the analytic behavior of $\psi_{10}(z)$, $\psi_{01}(z)$, $\psi_{11}(z)$, respectively, around the singular points on the circles with radiuses $\rho_{\vc{c}}$, where
\begin{eqnarray*}
 && \rho_{(1,0)} = \sup \{ u \ge 0; (u,1) \in \tilde{\sr{D}} \}, \quad \rho_{(0,1)} = \sup \{ u \ge 0; (1,u) \in \tilde{\sr{D}} \},\\
 && \rho_{(1,1)} = \sup \{ u \ge 0; (u,u) \in \tilde{\sr{D}} \}.
\end{eqnarray*}
  
  Since $\psi_{10}(z)$ and $\psi_{01}(z)$ are symmetric, we only consider $\psi_{10}(z)$ and $\psi_{11}(z)$. From \eqn{stationary equation 6}, we have
\begin{eqnarray}
\label{eqn:psi 10}
 \lefteqn{\psi_{10}(z) = \left(1+\frac{\tilde{\gamma}_{1}(z,1) - 1} {1 - \tilde{\gamma}_{+}(z,1) }\right) \tilde{\varphi}_{1}(z) + \frac{\tilde{\gamma}_{2}(z,1) - \tilde{\gamma}_{+}(z,1)} {1 - \tilde{\gamma}_{+}(z,1)} \tilde{\varphi}_{2}(1)} \nonumber \hspace{1ex}\\
 && \hspace{30ex}+ \frac{\tilde{\gamma}_{0}(z,1) - \tilde{\gamma}_{+}(z,1)} {1 - \tilde{\gamma}_{+}(z,1) } \nu(\vc{0}),\\
\label{eqn:psi 11}
 \lefteqn{\psi_{11}(z)  = \left(1+\frac{\tilde{\gamma}_{1}(z,z) - 1} {1 - \tilde{\gamma}_{+}(z,z) }\right) \tilde{\varphi}_{1}(z) + \frac{\tilde{\gamma}_{2}(z,z) - \tilde{\gamma}_{+}(z,z)} {1 - \tilde{\gamma}_{+}(z,z)} \tilde{\varphi}_{2}(z)} \nonumber \hspace{1ex}\\
 && \hspace{30ex} + \frac{\tilde{\gamma}_{0}(z,z) - \tilde{\gamma}_{+}(z,z)} {1 - \tilde{\gamma}_{+}(z,z) } \nu(\vc{0}).
\end{eqnarray}

\begin{figure}[h]
	\caption{Left: $\ol{\zeta}_{2}(u^{(1,\Gamma)}_{1}) < 1$, Right: $\ol{\zeta}_{2}(u^{(1,\Gamma)}_{1}) > 1$ and $\ul{\zeta}_{2}(u^{(1,\Gamma)}_{1}) \ne 1$} 
 	\centering
	\includegraphics[height=4cm]{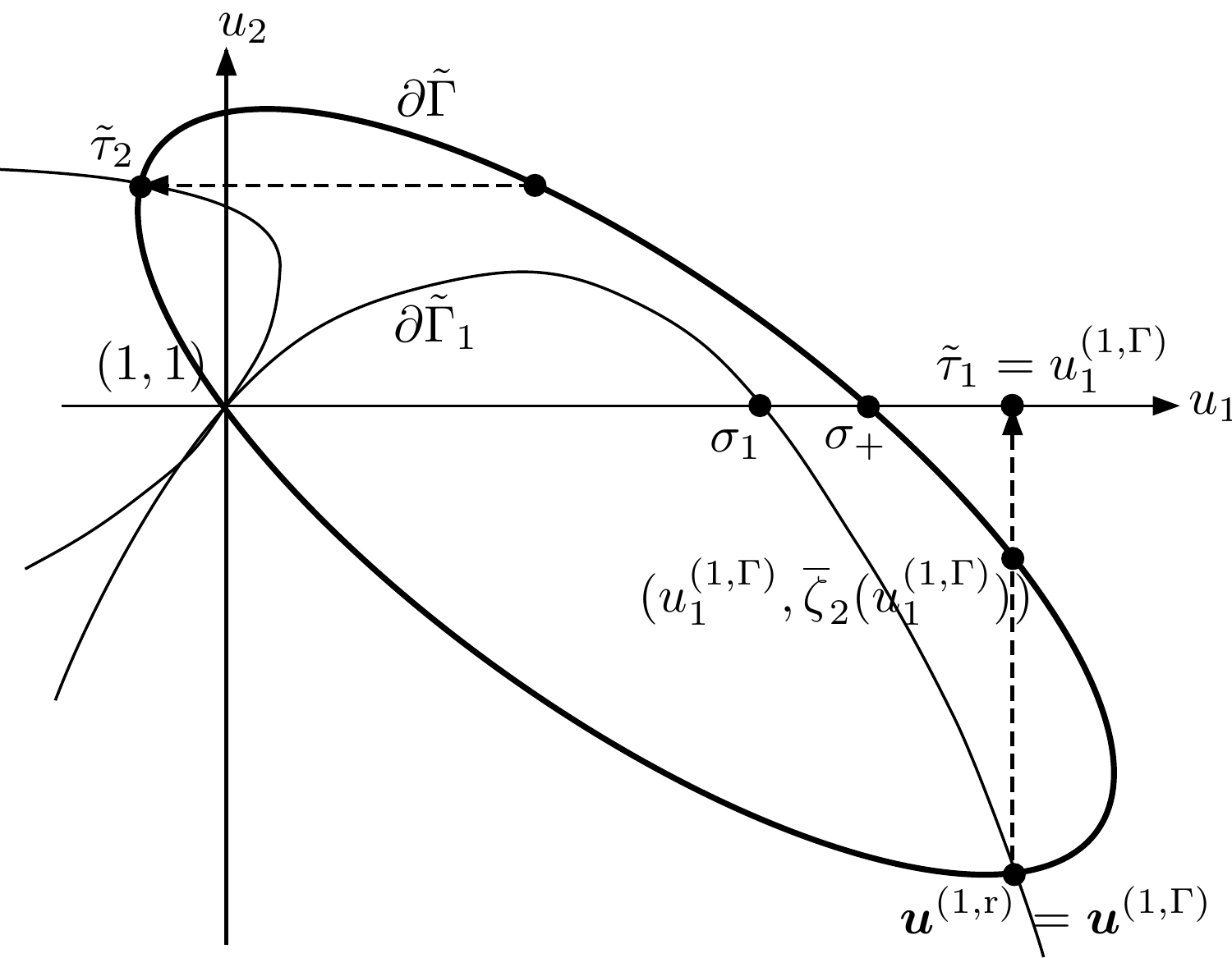} \hspace{2ex}
	\includegraphics[height=4cm]{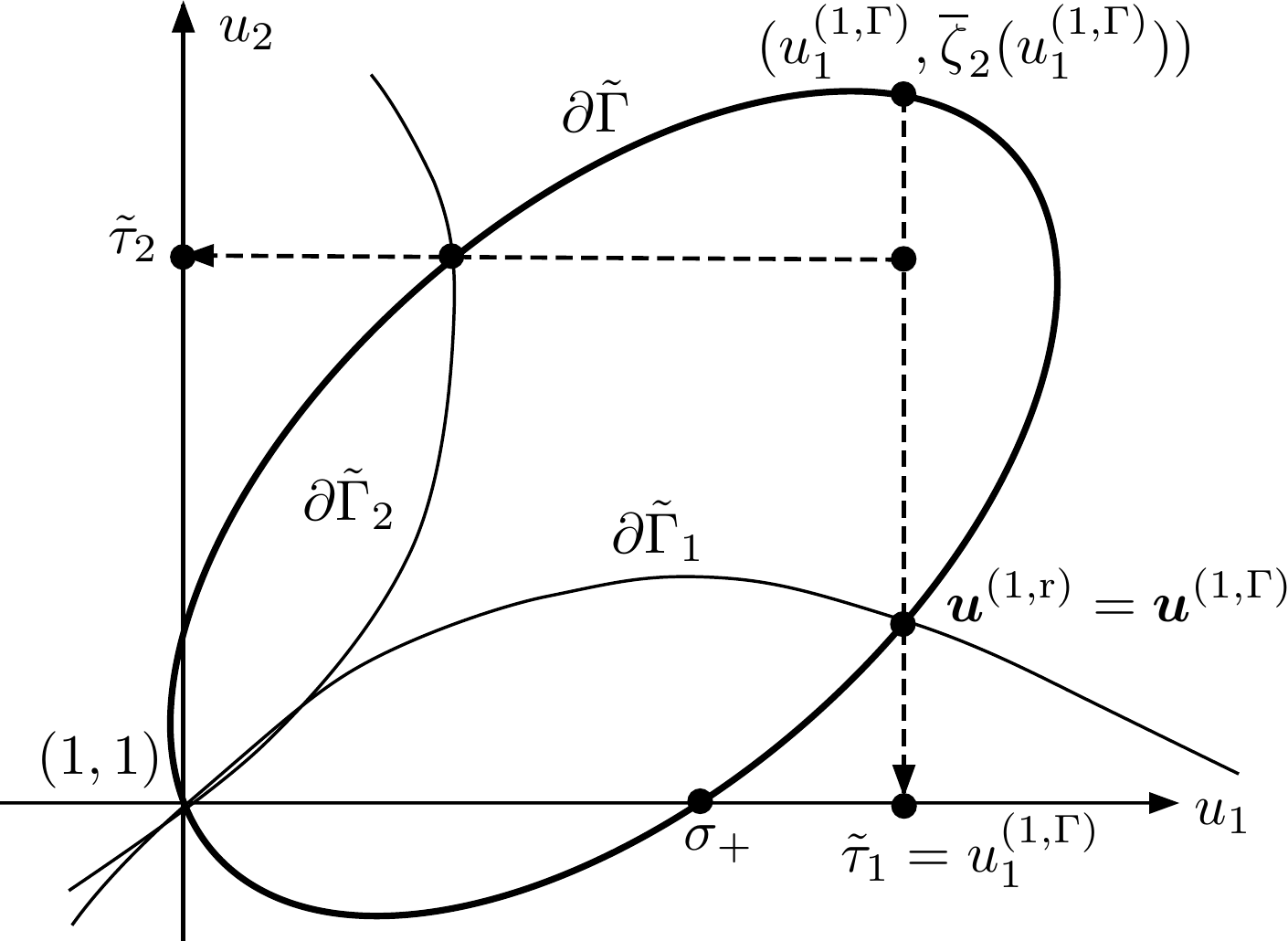} \\
	\label{fig:marginal 1}
	\caption{Left: $\ol{\zeta}_{2}(u^{(1,\Gamma)}_{1}) > 1$ and $\ul{\zeta}_{2}(u^{(1,\Gamma)}_{1}) = 1$, Right: $\ol{\zeta}_{2}(u^{(1,\Gamma)}_{1}) = \ul{\zeta}_{2}(u^{(1,\Gamma)}_{1}) = 1$} 
	\label{fig:marginal 2}
 	\centering
	\includegraphics[height=4cm]{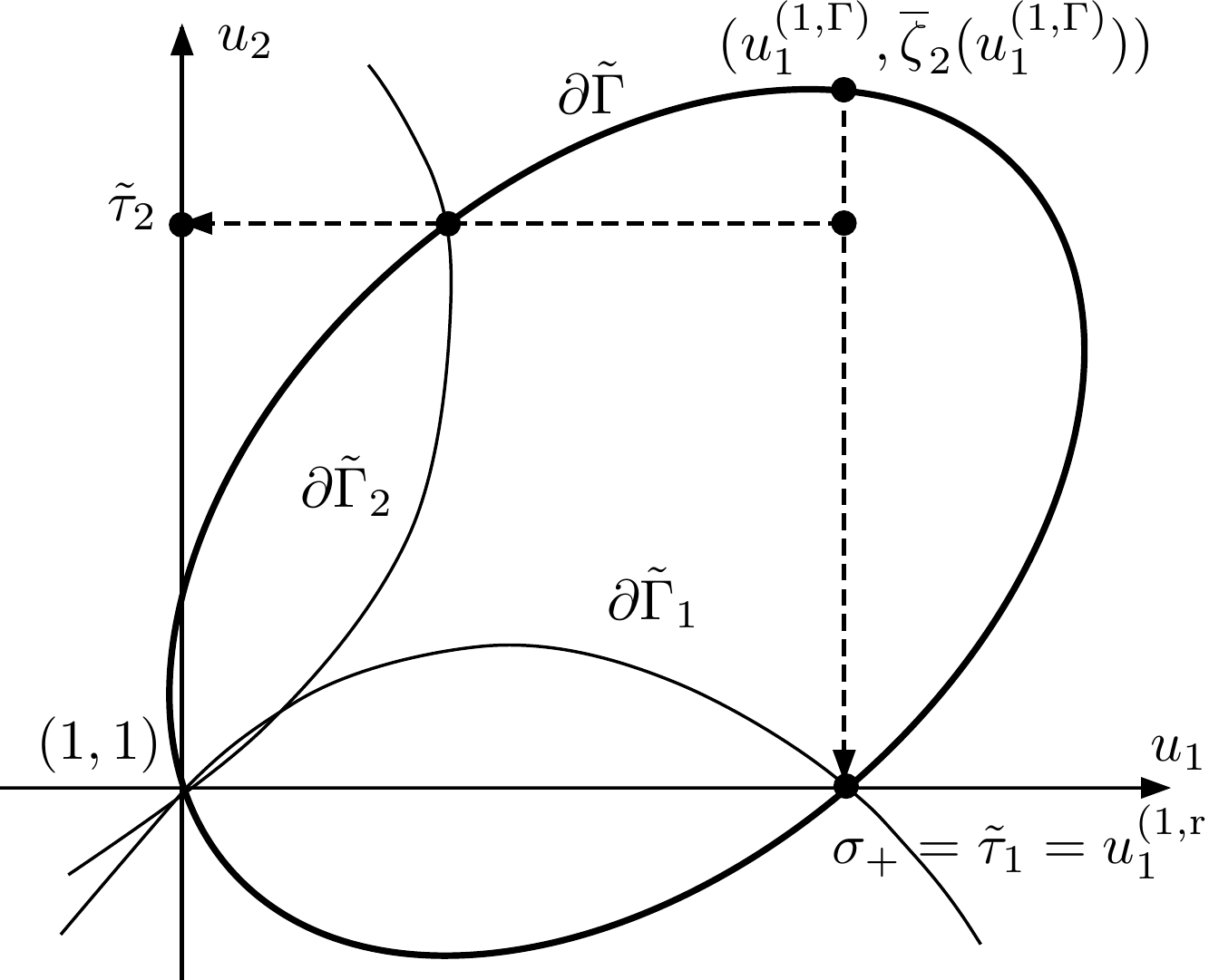} \hspace{2ex}
	\includegraphics[height=4cm]{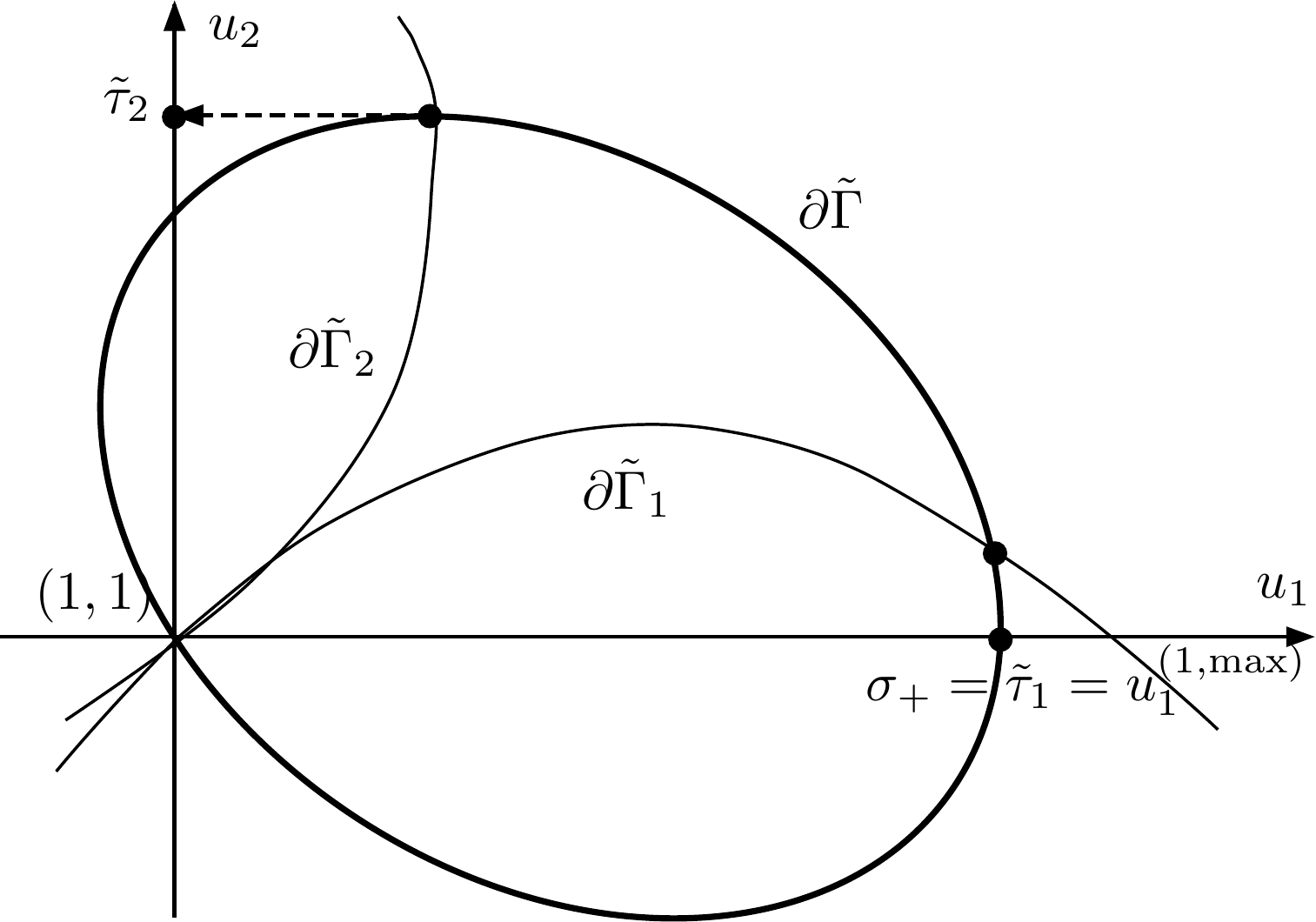} 
\end{figure}

We first consider the tail asymptotics for $\vc{c} = (1,0)$ under the non-arithmetic condition (v-a). From \eqn{psi 10}, the singularity of $\psi_{11}(z)$ on the circle $|z| = \rho_{(1,0)}$ occurs by either that of $\tilde{\varphi}_{1}(z)$ or the solution of the following equation:
\begin{eqnarray}
\label{eqn:gamma + 1}
  \tilde{\gamma}_{+}(z,1) = 1.
\end{eqnarray}
   Since this equation is quadratic and the domain $\tilde{\sr{D}}$ contains vectors $\vc{x} > \vc{1} \equiv (1,1)$, the equation \eqn{gamma + 1} has a unique real solution greater than 1. We denote it by $\sigma_{+}$. We then have the following asymptotics (see also \fig{marginal 1}). 
   
\begin{theorem} {\rm
\label{thr:exact marginal 10}
  Under the conditions (i)--(iv) and (v-a), let $h_{1}(n)$ be the exact asymptotic function given in Theorems \thrt{exact nu 10n-1}~and~\thrt{exact nu 10n-2}, then $P(L_{1} \ge n)$ has the following exact asymptotic $g_{1}(n)$ as $n \to \infty$.
\begin{mylist}{0}
\item [(a)] If $\ol{\zeta}_{2}(u^{(1,\Gamma)}_{1}) < 1$, then $g_{1}(n) = \sigma_{+}^{-n}$.
\item [(b)] If $\ol{\zeta}_{2}(u^{(1,\Gamma)}_{1}) > 1$ and $\ul{\zeta}_{2}(u^{(1,\Gamma)}_{1}) \ne 1$, then $g_{1}(n) = h_{1}(n)$.
\item [(c)] If $\ol{\zeta}_{2}(u^{(1,\Gamma)}_{1}) > 1 = \ul{\zeta}_{2}(u^{(1,\Gamma)}_{1})$, then $g_{1}(n) = \tilde{\tau}_{1}^{-n}$.
\item [(d)] If $\ol{\zeta}_{2}(u^{(1,\Gamma)}_{1}) = 1 = \ul{\zeta}_{2}(u^{(1,\Gamma)}_{1})$, then $g_{1}(n) = \tilde{\tau}_{1}^{-n}$.
\item [(e)] If $\ol{\zeta}_{2}(u^{(1,\Gamma)}_{1}) = 1 > \ul{\zeta}_{2}(u^{(1,\Gamma)}_{1})$, then $g_{1}(n) = n \tilde{\tau}_{1}^{-n}$.
\end{mylist}
}\end{theorem}
\begin{remark} {\rm
\label{rem:exact marginal 10}
  The corresponding but less complete results are obtained using moment generating functions in Corollary 4.3 of \cite{Miyazawa 2009}. 
}\end{remark}

\begin{figure}[h]
	\caption{Left: $\sigma_{\rm d} < \tilde{\tau}_{1}$, Right: $\sigma_{\rm d} > \tilde{\tau}_{1}$} 
 	\centering
	\includegraphics[height=4.5cm]{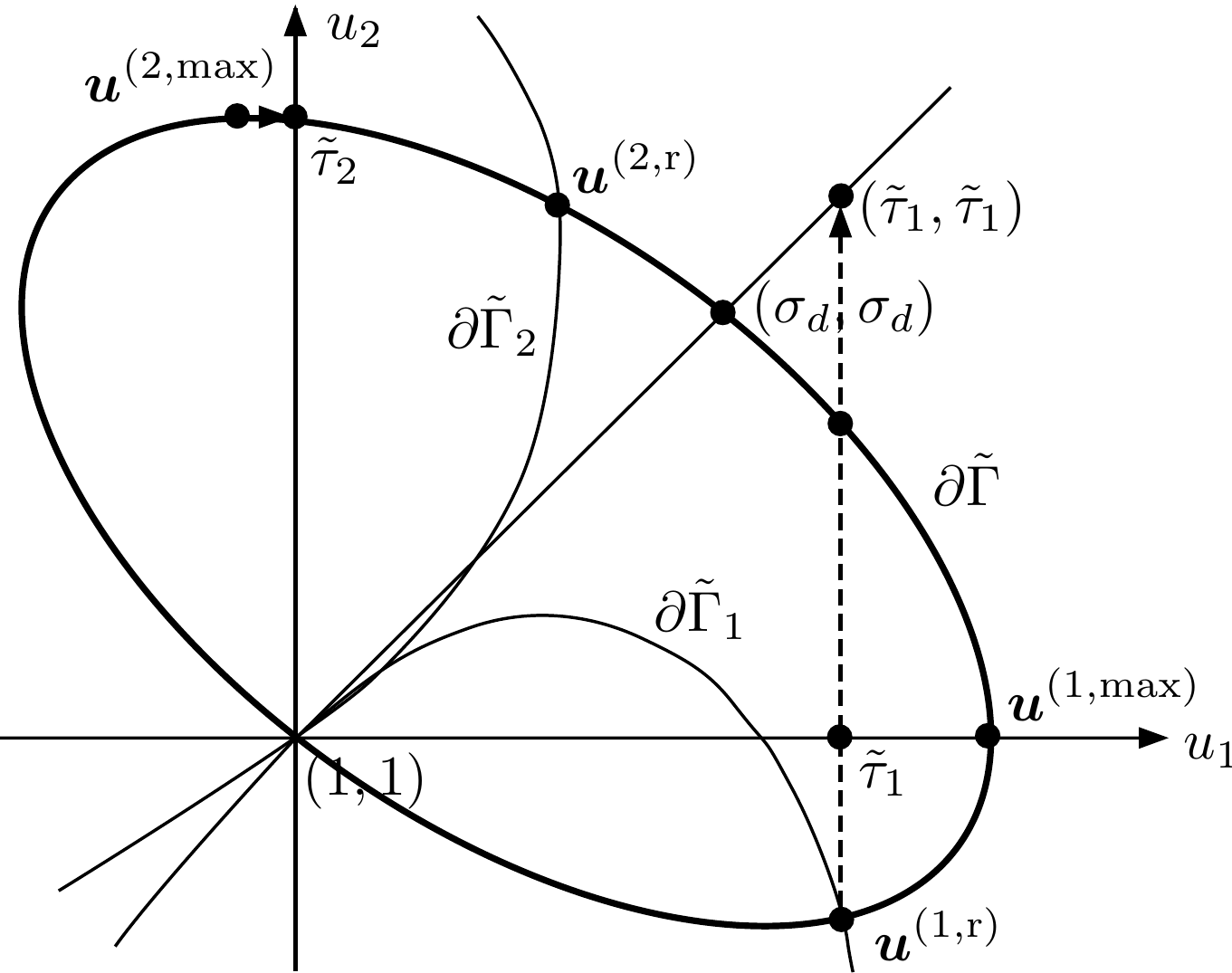} \hspace{2ex}
	\includegraphics[height=4.5cm]{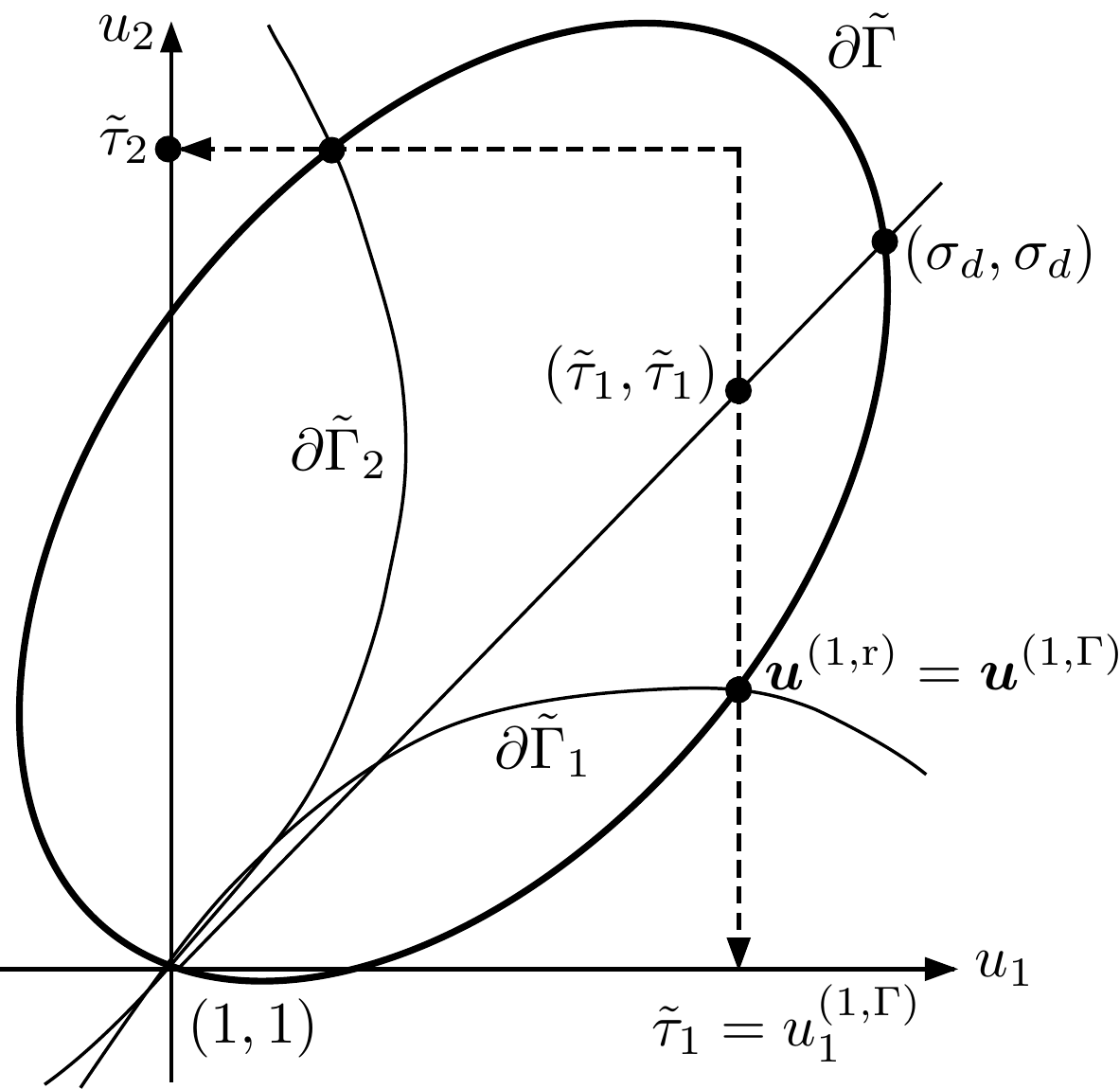} \vspace{2ex}\\
	\label{fig:diagonal 1}
	\caption{Left: $\sigma_{\rm d} = \tilde{\tau}_{1} = u^{(1,\max)}_{1} = \tilde{\tau}_{2}$, Right: $\sigma_{\rm d} = \tilde{\tau}_{1} = u^{(1,\max)}_{1} \ne \tilde{\tau}_{2}$} 
 	\centering
	\includegraphics[height=4.5cm]{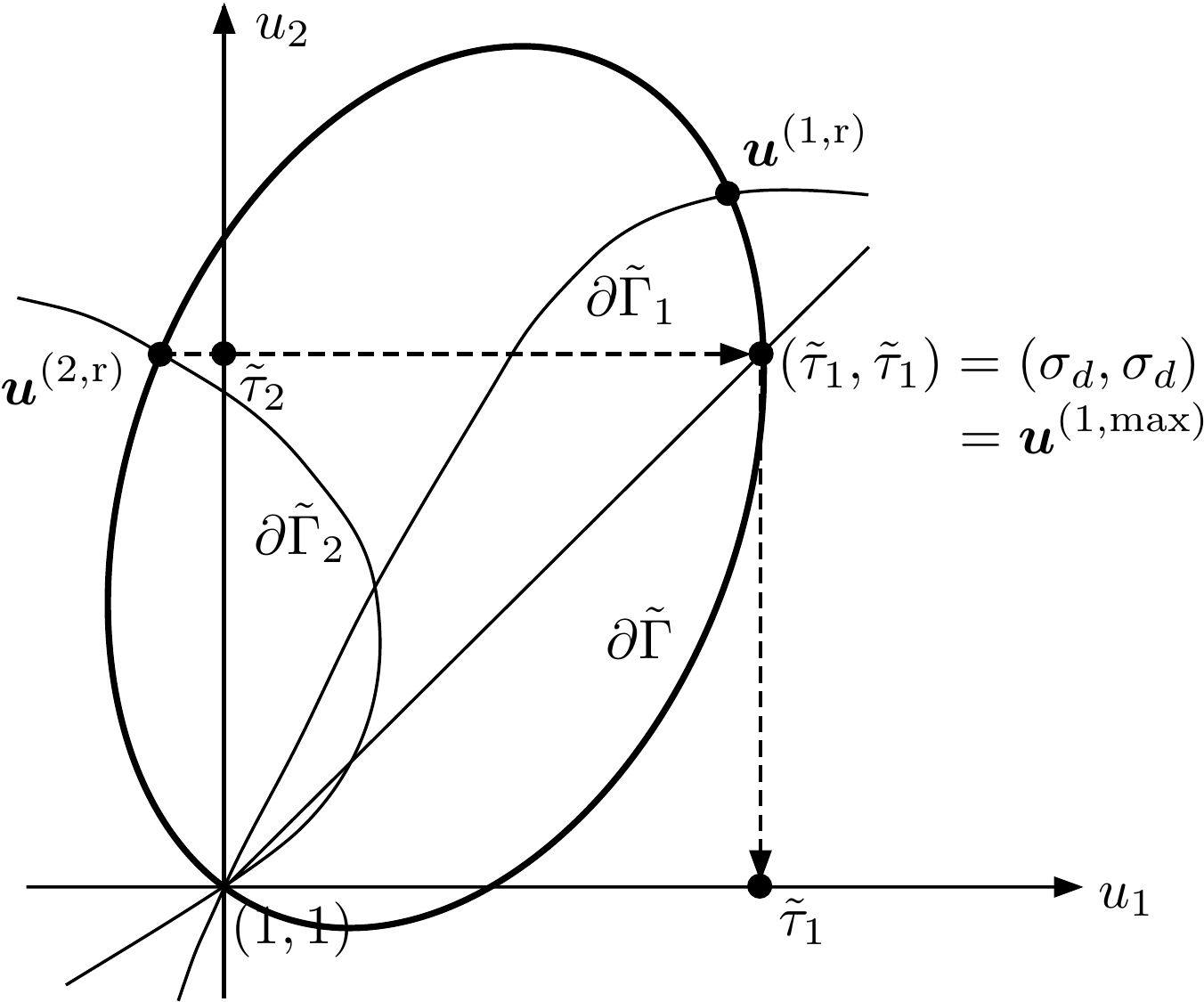} \hspace{2ex}
	\includegraphics[height=4.5cm]{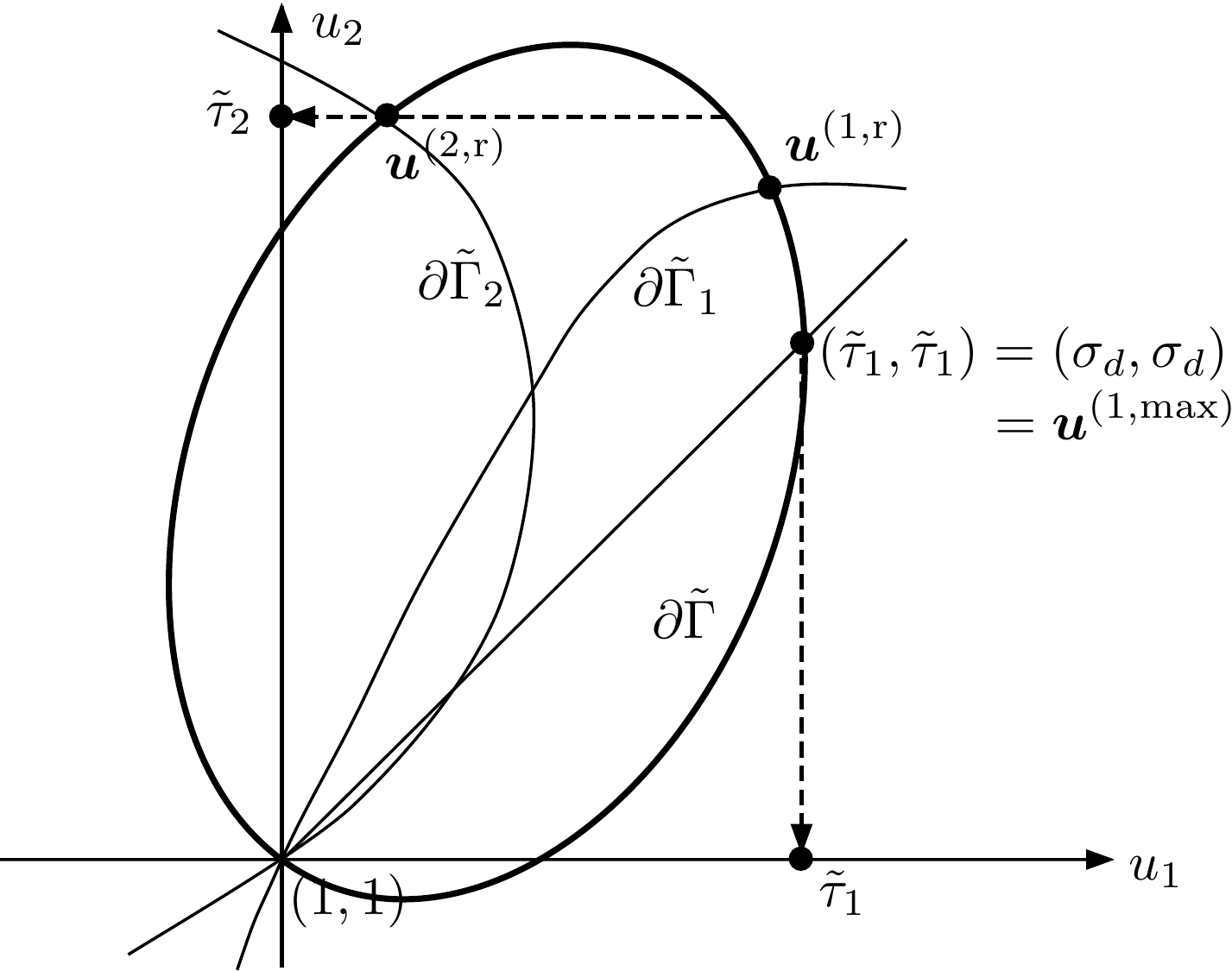} 
	\label{fig:diagonal 2}
\end{figure}

  Before proving this theorem, we present asymptotics for the marginal distribution in the diagonal direction. Let $\sigma_{\rm d} $ be the real solution of
\begin{eqnarray*}
  \tilde{\gamma}_{+}(u,u) = 1, \qquad u > 1,
\end{eqnarray*}
  which can be shown to be unique (see \fig{diagonal 1}).  Because of symmetry, we assume without loss of generality that $\tilde{\tau}_{1} \le \tilde{\tau}_{2}$. See \fig{marginal 2} for the location of this point.

\begin{theorem} {\rm
\label{thr:exact marginal 11}
  Under the conditions (i)--(iv), (v-a) and $\tilde{\tau}_{1} \le \tilde{\tau}_{2}$, let $h_{1}(n)$ be the exact asymptotic function given in Theorems \thrt{exact nu 10n-1}~and~\thrt{exact nu 10n-2}, then $P(L_{1} + L_{2} \ge n)$ has the following exact asymptotic $g_{+}(n)$ as $n \to \infty$.
\begin{mylist}{0}
\item [(a)] If $\sigma_{\rm d} < \tilde{\tau}_{1}$, then $g_{+}(n) = \sigma_{\rm d}^{-n}$.
\item [(b)] If $\sigma_{\rm d} > \tilde{\tau}_{1}$, then $g_{+}(n) = h_{1}(n)$.
\item [(c)] If $\sigma_{\rm d} = \tilde{\tau}_{1} \ne u^{(1,\max)}_{1}$, then $g_{+}(n) = n \sigma_{\rm d}^{-n}$.
\item [(d)] If $\sigma_{\rm d} = \tilde{\tau}_{1} = u^{(1,\max)}_{1} = \tilde{\tau}_{2}$, then $g_{+}(n) = n \sigma_{\rm d}^{-n}$
\item [(e)] If $\sigma_{\rm d} = \tilde{\tau}_{1} = u^{(1,\max)}_{1} \ne \tilde{\tau}_{2}$, then $g_{+}(n) = \sigma_{\rm d}^{-n}$.
\end{mylist}

}\end{theorem}

  In what follows, we prove \thr{exact marginal 10}. The proof of \thr{exact marginal 11} is similar, so we only shortly outline it.

\begin{proof*}{Proof of \thr{exact marginal 10}}
  Let
\begin{eqnarray*}
  \xi(z) = (\tilde{\gamma}_{2}(z,1) - \tilde{\gamma}_{+}(z,1)) \tilde{\varphi}_{2}(1) + (\tilde{\gamma}_{0}(z,1) - \tilde{\gamma}_{+}(z,1)) \nu(\vc{0}),
\end{eqnarray*}
  then \eqn{psi 10} can be written as
\begin{eqnarray}
\label{eqn:psi 10 2}
 \psi_{10}(z) = \left(1+\frac{\tilde{\gamma}_{1}(z,1) - 1} {1 - \tilde{\gamma}_{+}(z,1) }\right) \tilde{\varphi}_{1}(z) + \frac{\xi(z)} {1 - \tilde{\gamma}_{+}(z,1)} .
\end{eqnarray}
   Since $\tilde{\gamma}_{2}(u,1) > 1, \tilde{\gamma}_{0}(u,1) > 1$ for $u > 0$ and
\begin{eqnarray*}
\left. \frac {\partial} {\partial u} \tilde{\gamma}_{1}(u,1) \right|_{u=\sigma_{1}} < 0, \qquad \left. \frac {\partial} {\partial u} \tilde{\gamma}_{+}(u,1)\right|_{u=\sigma_{+}} > 0 \; \mbox{ if $\ul{\zeta}_{2}(\sigma_{+}) = 0$},
\end{eqnarray*}
  where $\sigma_{1}$ is a positive number satisfying that $\tilde{\gamma}_{1}(\sigma_{1},1)=1$, $\xi(\sigma_{+}) > 0$ and $\sigma_{+} = \sigma_{1}$ implies that the prefactor of $\tilde{\varphi}_{1}(z)$ is positive at $z = \sigma_{+}$ if $\ul{\zeta}_{2}(\sigma_{+}) = 0$. Having these observations in mind, we prove each cases.
  
\noindent (a) Assume that $\ol{\zeta}_{2}(u^{(1,\Gamma)}_{1}) < 1$. This occurs if and only if $\sigma_{+} = \rho_{10} < \tilde{\tau}_{1}$  (see the left picture of \fig{marginal 1}). In this case, $\psi_{10}(z)$ must be singular at $z = \sigma_{+}$ because it is one the boundary of the convergence domain $\tilde{\sr{D}}$. Hence, it has a simple pole at $z = \sigma_{+}$, and therefore we have the exact geometric asymptotic.

\noindent (b) Assume that $\ol{\zeta}_{2}(u^{(1,\Gamma)}_{1}) > 1$ and $\ul{\zeta}_{2}(u^{(1,\Gamma)}_{1}) \ne 1$. This case occurs if and only if $\sigma_{+} \ne \rho_{10} = \tilde{\tau}_{1}$ (see the right picture of \fig{marginal 1}). In this case, $\tilde{\gamma}_{1}(\tilde{\tau}_{1},1) \ne 1$, $\tilde{\gamma}_{+}(\tilde{\tau}_{1},1) \ne 1$ and $\tilde{\gamma}_{1}(\tilde{\tau}_{1},1) - 1$ has the same sign as $1 - \tilde{\gamma}_{+}(\tilde{\tau}_{1},1)$. Hence, the prefactor of $\tilde{\varphi}_{1}(z)$ is analytic at $z = \tilde{\tau}_{1}$, and the singularity of $\psi_{10}(z)$ is determined by $\tilde{\varphi}_{1}(z)$. Thus, we have the same asymptotics as in Theorems \thrt{exact nu 10n-1} and \thrt{exact nu 10n-2}. 

\noindent (c) Assume that $\ol{\zeta}_{2}(u^{(1,\Gamma)}_{1}) > 1$ and $\ul{\zeta}_{2}(u^{(1,\Gamma)}_{1}) = 1$ (see the left figure of \fig{marginal 2}). In this case, $\tilde{\gamma}_{+}(\tilde{\tau}_{1},1) = \tilde{\gamma}_{1}(\tilde{\tau}_{1},1) = 1$ and category \rmn{2} is impossible, and therefore, from \eqn{psi 10 2} and \thr{exact nu 10n-1}, we have the exact geometric asymptotic.

\noindent (d) Assume that $\ol{\zeta}_{2}(u^{(1,\Gamma)}_{1}) = 1 = \ul{\zeta}_{2}(u^{(1,\Gamma)}_{1})$ (see the right figure of \fig{marginal 2}). In this case, $\tilde{\tau}_{1} = \sigma_{+} = u^{(1,\max)}_{1}$, and therefore $\tilde{\gamma}_{+}(\tilde{\tau}_{1},1) = 1$. We need to consider two subcases, $u^{(1,\edge)}_{1} = u^{(1,\max)}_{1}$ and $u^{(1,\edge)}_{1} \ne u^{(1,\max)}_{1}$. If $u^{(1,\edge)}_{1} = u^{(1,\max)}_{1}$, then $\tilde{\gamma}_{1}(\tilde{\tau}_{1},1) = 1$ and $\tilde{\varphi}_{1}(z) \sim (\tilde{\tau}_{1} - z)^{-\frac 12}$ by \thr{exact nu 10n-1}. Thus, we have $\psi_{10}(z) \sim (\tilde{\tau}_{1} - z)^{-1}$ due to the second term of \eqn{psi 10 2}. Otherwise, if $u^{(1,\edge)}_{1} \ne u^{(1,\max)}_{1}$, then $\tilde{\gamma}_{1}(\tilde{\tau}_{1},1) \ne 1$ implies that the prefactor of $\tilde{\varphi}_{1}(z)$ in \eqn{psi 10} has a single pole at $z = \tilde{\tau}_{1}$ and that $\tilde{\varphi}_{1}(z) - \tilde{\varphi}_{1}(u^{(1,\max)}_{1}) \sim (\tilde{\tau}_{1} - z)^{\frac 12}$. Again from \eqn{psi 10 2}, we have $\psi_{10}(z) \sim (\tilde{\tau}_{1} - z)^{-1}$. Thus, we have the exact geometric asymptotic in both cases.

\noindent (e) Assume that $\ol{\zeta}_{2}(u^{(1,\Gamma)}_{1}) = 1 > \ul{\zeta}_{2}(u^{(1,\Gamma)}_{1})$. In this case, $\tilde{\tau}_{1} = \sigma_{+} = u^{(1,\edge)}_{1} < u^{(1,\max)}_{1}$ and we must have category \rmn{1} or \rmn{3}. Since $\tilde{\gamma}_{+}(\tilde{\tau}_{1},1) = 1$, $\tilde{\gamma}_{1}(\tilde{\tau}_{1},1) > 1$ and $\tilde{\varphi}_{1}(z)$ has a single pole at $z = \tilde{\tau}_{1}$, $\psi_{10}(z)$ in \eqn{psi 10} has a double pole at $z = \tilde{\tau}_{1}$. This yields the desired asymptotic.
\end{proof*}

\medskip

  The proof of \thr{exact marginal 11} is more or less similar to that of $P(L_{1} \ge n)$. From Figures \figt{diagonal 1} and \figt{diagonal 2}, we can see how the dominant singular point is located. Since its derivation is routine, we omit detailed proof.

\section{Exact tail asymptotics for the arithmetic case} 
\label{sect:exact asymptotic not v-a}

  Throughout this section, we assume that (v-a) does not hold. As in \sectn{arithmetic case}, we separately consider two cases: (B1) either (v-b) or $m^{(1)}_{2} = 0$ holds, and (B2) neither (v-b) nor $m^{(1)}_{2} = 0$ holds, according to \tab{period}. In some cases, we need: (C1) either (v-c) or $m^{(2)}_{1} = 0$ holds, and (C2) neither (v-c) nor $m^{(2)}_{1} = 0$ holds.

\subsection{The boundary probabilities for the arithmetic case with (B1)}
\label{sect:boundary not v-a (B1)}

In this case, we have the following asymptotics.

\begin{theorem} {\rm
\label{thr:exact nu 10A-1}
  Under the conditions (i)--(iv) and (B1), if (v-a) dose not hod, then for categories \rmn{1} and \rmn{3}, $\tilde{\tau}_{1} = u^{(1,\cp)}_{1}$, and $P(L_{1} = n, L_{2} = 0)$ has the following exact asymptotic $h_{2}(n)$. For some constant $b \in [-1, 1]$,
\begin{eqnarray}
\label{eqn:exact nu 10A-1}
  h_{2}(n) = \left\{\begin{array}{ll}
  \tilde{\tau}_{1}^{-n}, \quad & u^{(1,\cp)}_{1} \not= u^{(1,\max)}_{1},\\
  n^{- \frac 12} \tilde{\tau}_{1}^{-n}, \quad & u^{(1,\cp)}_{1} = u^{(1,\max)}_{1} = u^{(1,\edge)}_{1},\\
  n^{- \frac 32} \left( 1+ b (-1)^{n} \right) \tilde{\tau}_{1}^{-n}, \quad & u^{(1,\cp)}_{1} = u^{(1,\max)}_{1} \ne u^{(1,\edge)}_{1}.
  \end{array} \right.
\end{eqnarray}
  By symmetry, the corresponding results are also obtained for $P(L_{1} = 0, L_{2} = n)$ for categories \rmn{1} and \rmn{2}.
}\end{theorem}

\begin{theorem} {\rm
\label{thr:exact nu 10A-2}
  Under the conditions (i)--(iv) and (B1), if (v-a) dose not hod, then, for category \rmn{2}, $\tilde{\tau}_{1} = \ul{\zeta}_{2}(\tilde{\tau}_{2})$, $\tilde{\tau}_{2} = u^{(2,\edge)}_{2}$, and $P(L_{1} = n, L_{2} = 0)$ has the following exact asymptotic $h_{2}(n)$. For some constant $b \in [-1, 1]$,
\begin{eqnarray}
\label{eqn:exact nu 10A-2}
  h_{2}(n) = \left\{\begin{array}{ll}
  \tilde{\tau}_{1}^{-n}, \quad & \tilde{\tau}_{1} < u^{(1,\cp)}_{1} \mbox{ or }  \\
  & \tilde{\tau}_{1} = u^{(1,\cp)}_{1} = u^{(1,\max)}_{1} = u^{(1,\edge)}_{1}, \\
  n \tilde{\tau}_{1}^{-n}, \quad & \tilde{\tau}_{1} = u^{(1,\cp)}_{1} \not= u^{(1,\max)}_{1},\\
  n^{- \frac 12} \tilde{\tau}_{1}^{-n}, \; & \tilde{\tau}_{1} = u^{(1,\cp)}_{1} = u^{(1,\max)}_{1} \ne u^{(1,\edge)}_{1},\\ 
  & \mbox{and (C1) holds}.\\
  n^{- \frac 12} (1 + b (-1)^{n}) \tilde{\tau}_{1}^{-n}, \; & \tilde{\tau}_{1} = u^{(1,\cp)}_{1} = u^{(1,\max)}_{1} \ne u^{(1,\edge)}_{1},\\  
  & \mbox{and (C2) holds}.\\
\end{array} \right.
\end{eqnarray}
  By symmetry, the corresponding results are also obtained for $P(L_{1} = 0, L_{2} = n)$ for categories \rmn{3}.
}\end{theorem}

\begin{remark} {\rm
\label{rem:exact nu 10A 1}
  As we will see in the proofs of these theorems, the asymptotics can be refined for those with the same geometric decay term $\tilde{\tau}_{1}^{-n}$. There is no difficulty to find them, but they are cumbersome because we need further cases. Thus, we omit their details.
}\end{remark}

\begin{remark} {\rm
\label{rem:exact nu 10A 2}
  One may wonder whether $b = \pm 1$ can occur in Theorems \thrt{exact nu 10A-1} and \thrt{exact nu 10A-2}. If it is the case, then the tail asymptotics are purely periodic. Closely look at the coefficients of the asymptotic expansion of the terms in \eqn{stationary equation 4}, this unlikely occurs because $|\tilde{\varphi}_{2}(-\ul{\zeta}_{2}(\tilde{\tau}_{1}))| < \tilde{\varphi}_{2}(\ul{\zeta}_{2}(\tilde{\tau}_{1}))$. Thus, we conjecture that $|b| < 1$ is always the case.
}\end{remark}

By \tab{period}, $\tilde{\varphi}_{1}(z)$ may be singular at $z = - \tilde{\tau}_{1}$ on $|z| = \tilde{\tau}_{1}$. On the other hand, $\tilde{\varphi}_{1}(z)$ has the same singularity at $z = \tilde{\tau}_{1}$ as in the non-arithmetic case, so we can only focus on the singularity at $z = -\tilde{\tau}_{1}$. We note that $z=-u^{(1,\edge)}_{1}$ can not be the solution of \eqn{no solution 1} under the assumptions of Theorems \thrt{exact nu 10A-1} and \thrt{exact nu 10A-2}. Having these in mind, we give proofs.

\medskip

\begin{proof*}{Proof of \thr{exact nu 10A-1}}
We consider the singularity of $\tilde{\varphi}_{1}(z)$ at $z = - \tilde{\tau}_{1}$ by \eqn{stationary equation 3} using the arguments in Sections \sect{arithmetic case} and \sect{exact asymptotic v-a}. Note that $\tilde{\tau}_{1} = u^{(1,\Gamma)}_{1}$ because the category is either \rmn{1} or \rmn{3}. We need to consider the following three cases.
\begin{mylist}{2}
\item [(\sect{exact asymptotic not v-a}a)] $u^{(1,\cp)}_{1} \ne u^{(1,\max)}_{1}$: This case is equivalent to $u^{(1,\cp)}_{1} < u^{(1,\max)}_{1}$, and it follows from \eqn{stationary equation 3} that $\tilde{\varphi}_{1}(z)$ is analytic at $z = -u^{(1,\edge)}_{1}$. Hence, there is no singularity contribution by $z = -u^{(1,\edge)}_{1}$.

\item [(\sect{exact asymptotic not v-a}b)] $u^{(1,\cp)}_{1} = u^{(1,\max)}_{1}$, $u^{(1,\edge)}_{1} = u^{(1,\max)}_{1}$: In this case, as $z \to -u^{(1,\max)}_{1}$ in such a way that $z \in \tilde{\sr{G}}^{+}_{\delta}(-u^{(1,\max)}_{1})$ for some $\delta > 0$,
\begin{eqnarray*}
  \tilde{\varphi}_{2}(\ul{\zeta}_{2}(z)) - \tilde{\varphi}_{2}(\ul{\zeta}_{2}(-u^{(1,\max)}_{1})) \sim (-u^{(1,\max)}_{1} - z)^{\frac 12},
\end{eqnarray*}
  but $1 - \tilde{\gamma}_{1}(z,\ul{\zeta}_{2}(z))$ does not vanish at $z = -u^{(1,\max)}_{1}$, and therefore
\begin{eqnarray*}
  \tilde{\varphi}_{1}(z) - \tilde{\varphi}_{1}(-u^{(1,\max)}_{1}) \sim (-u^{(1,\max)}_{1} - z)^{\frac 12}.
\end{eqnarray*}
  This yields the asymptotic function $n^{-\frac 32} \tilde{\tau}_{1}^{-n}$, but this function is dominated by the slower asymptotic function $n^{-\frac 12} \tilde{\tau}_{1}^{-n}$ due to the singularity at $z = u^{(1,\max)}_{1}$. 
  
\item [(\sect{exact asymptotic not v-a}c)]$ u^{(1,\cp)}_{1} = u^{(1,\max)}_{1}$, $u^{(1,\edge)}_{1} \not= u^{(1,\max)}_{1}$: In this case, the solution of \eqn{no solution 1} has no essential role, so $\tilde{\varphi}_{1}(z)$ has the same analytic behavior at $z = - u^{(1,\max)}_{1}$ as at $z = u^{(1,\max)}_{1}$ in (\sect{exact asymptotic v-a}c) in the proof of \thr{exact nu 10n-1}.
\end{mylist}
Thus, combining with the asymptotics in \thr{exact nu 10n-1}, we complete the proof.
\end{proof*}

\medskip

\begin{proof*}{Proof of \thr{exact nu 10A-2}}
Because of category \rmn{2}, $\tau_{2} = \ul{\zeta}_{2}(\tilde{\tau}_{1})$, and therefore $\ul{\zeta}_{2}(-\tilde{\tau}_{1}) = - \ul{\zeta}_{2}(\tilde{\tau}_{1}) = - \tilde{\tau}_{2}$ by the assumption that (v-a) does not hold. We consider the singularity at $z = -\tilde{\tau}_{1}$ for the following cases with this in mind.

\begin{mylist}{2}

\item [(\sect{exact asymptotic not v-a}a')] $\tilde{\tau}_{1} < u^{(1,\cp)}_{1}$: This case is included in (\sect{Analytic function}c'-2). Hence, if (C1) holds, $\tilde{\varphi}_{2}(\ul{\zeta}_{2}(z))$ and therefore $\tilde{\varphi}_{1}(z)$ are analytic at $z = -u^{(1,\edge)}_{1}$. Otherwise, if (C2) holds, $\tilde{\varphi}_{2}(\ul{\zeta}_{2}(z))$ has a simple pole at $z = -u^{(1,\edge)}_{1}$. However, in \eqn{stationary equation 3}, $\tilde{\varphi}_{2}(\ul{\zeta}_{2}(z))$ has the prefactor, $\tilde{\gamma}_{2}(z, \ul{\zeta}_{2}(z)) - 1$, which vanishes at $z = -u^{(1,\edge)}_{1}$ because of (C2). Hence, the pole of $\tilde{\varphi}_{2}(\ul{\zeta}_{2}(z))$ is cancelled, and therefore $\tilde{\varphi}_{1}(z)$ is analytic at $z = -u^{(1,\edge)}_{1}$. Thus, either case has no contribution by $z = -u^{(1,\edge)}_{1}$. 

\item [(\sect{exact asymptotic not v-a}b')] $\tilde{\tau}_{1} = u^{(1,\cp)}_{1} \ne u^{(1,\max)}_{1}$: In this case, $\tilde{\tau}_{1} = u^{(1,\edge)}_{1}$. If (C2) holds, then $\tilde{\varphi}_{2}(z)$ has a simple pole at $z = - \tilde{\tau}_{2}$, and therefore as in (\sect{exact asymptotic v-a}b'-3-1), 
\begin{eqnarray*}
  \tilde{\varphi}_{2}(\ul{\zeta}_{2}(z)) \sim (-u^{(1,\max)}_{1} - z)^{-\frac 12},
\end{eqnarray*}
  but (\sect{exact asymptotic v-a}b'-3-2) is not the case, and therefore this yields the asymptotic function $n^{-\frac 12} \tilde{\tau}_{1}^{-n}$. However, this asymptotic term is again dominated by $\tilde{\tau}_{1}^{-n}$ due to the singularity at $z = u^{(1,\max)}_{1}$. On the other hand, if (C1) holds, then there is no singularity contribution by $z = - \tilde{\tau}_{1}$. Hence, we have the same asymptotics as in the corresponding case of \thr{exact nu 10n-2}.
  
\item [(\sect{exact asymptotic not v-a}c')]$\tilde{\tau}_{1} = u^{(1,\cp)}_{1} = u^{(1,\max)}_{1}$: This is the case of (\sect{Analytic function}c'-3). As we discussed there, if (C2) holds, $\tilde{\varphi}_{2}(\ul{\zeta}_{2}(z)) \sim (-u^{(1,\max)}_{1} - z)^{-\frac 12}$ around $z = - u^{(1,\max)}_{1}$. Because of (B1), there is no other singularity contribution in \eqn{stationary equation 3}, and therefore we also have $\tilde{\varphi}_{1}(z) \sim (-u^{(1,\max)}_{1} - z)^{-\frac 12}$ around $z = - u^{(1,\max)}_{1}$. This results the asymptotic $n^{-\frac 12} \tilde{\tau}_{1}^{-n}$. On the other hand, if (C1) holds, we similarly have $\tilde{\varphi}_{1}(z) - \tilde{\varphi}_{1}(-u^{(1,\max)}_{1}) \sim (-u^{(1,\max)}_{1} - z)^{\frac 12}$. This implies the asymptotic $n^{-\frac 32} \tilde{\tau}_{1}^{-n}$. To combine this with the corresponding asymptotics obtained in \thr{exact nu 10n-2}, we consider two subcases.
\begin{mylist}{0}
\item [(\sect{exact asymptotic not v-a}c'-1)] $u^{(1,\max)}_{1} = u^{(1,\edge)}_{1}$: In this case, the asymptotics caused by $z = \tilde{\tau}_{1}$ is $n \tilde{\tau}_{1}^{-n}$, and therefore the asymptotics due to $z = - u^{(1,\max)}_{1}$ is ignorable.
\item [(\sect{exact asymptotic not v-a}c'-2)] $u^{(1,\max)}_{1} \ne u^{(1,\edge)}_{1}$: In this case, the asymptotics caused by $z = \tilde{\tau}_{1}$ is $n^{-\frac 12} \tilde{\tau}_{1}^{-n}$. Hence, we have two different cases. If (C1) holds, the contribution by $z = -\tilde{\tau}_{1}$ is ignorable. Otherwise, if (C2) holds, then we have additional asymptotic term: $(-1)^{n} n^{-\frac 12} \tilde{\tau}_{1}^{-n}$.
\end{mylist}
\end{mylist}
Thus, the proof is completed.
\end{proof*}

\subsection{The boundary probabilities for the arithmetic case with (B2)}
\label{sect:boundary not v-a (B2)}

We next consider case (B2). As noted in \sectn{arithmetic case}, in this case, $\tilde{\varphi}_{1}(z)$ is singular at $z = \pm u^{(1,\edge)}_{1}$, and both singular points have essentially the same properties. Thus, we have the following theorems.

\begin{theorem} {\rm
\label{thr:exact nu 10B-1}
  Under the conditions (i)--(iv) and (B2), if (v-a) dose not hod, then for categories \rmn{1} and \rmn{3}, $\tilde{\tau}_{1} = u^{(1,\cp)}_{1}$, and $P(L_{1} = n, L_{2} = 0)$ has the following exact asymptotic $h_{3}(n)$. For some constants $b_{i} \in [-1, 1]$ for $i=1, 2, 3$,
\begin{eqnarray}
\label{eqn:exact nu 10B-1}
  h_{3}(n) = \left\{\begin{array}{ll}
  (1 + b_{1} (-1)^{n}) \tilde{\tau}_{1}^{-n}, \quad & u^{(1,\cp)}_{1} \not= u^{(1,\max)}_{1},\\
  n^{- \frac 12} (1 + b_{2} (-1)^{n}) \tilde{\tau}_{1}^{-n}, \quad & u^{(1,\cp)}_{1} = u^{(1,\max)}_{1} = u^{(1,\edge)}_{1},\\
  n^{- \frac 32} (1 + b_{3} (-1)^{n}) \tilde{\tau}_{1}^{-n}, \quad & u^{(1,\cp)}_{1} = u^{(1,\max)}_{1} \ne u^{(1,\edge)}_{1}. 
  \end{array} \right.
\end{eqnarray}
  By symmetry, the corresponding results are also obtained for $P(L_{1} = 0, L_{2} = n)$ for categories \rmn{1} and \rmn{2}.
}\end{theorem}

\begin{theorem} {\rm
\label{thr:exact nu 10B-2}
  Under the conditions (i)--(iv) and (B2), if (v-a) dose not hod, then, for category \rmn{2}, $\tilde{\tau}_{2} = u^{(2,\edge)}_{2}$, and $P(L_{1} = n, L_{2} = 0)$ has the following exact asymptotic $h_{3}(n)$. For some constants $b_{i} \in [-1, 1]$ for $i=1, 2, 3$,
\begin{eqnarray}
\label{eqn:exact nu 10B-2}
  h_{3}(n) = \left\{\begin{array}{ll}
  (1 + b_{1} (-1)^{n}) \tilde{\tau}_{1}^{-n}, \quad & \tilde{\tau}_{1} < u^{(1,\cp)}_{1} \mbox{ or }  \\
  & \tilde{\tau}_{1} = u^{(1,\cp)}_{1} = u^{(1,\max)}_{1} = u^{(1,\edge)}_{1}, \\
  n (1 + b_{2} (-1)^{n}) \tilde{\tau}_{1}^{-n}, \quad & \tilde{\tau}_{1} = u^{(1,\cp)}_{1} \not= u^{(1,\max)}_{1},\\
  n^{- \frac 12} (1 + b_{3} (-1)^{n}) \tilde{\tau}_{1}^{-n}, \; & \tilde{\tau}_{1} = u^{(1,\cp)}_{1} = u^{(1,\max)}_{1} \ne u^{(1,\edge)}_{1}.  \end{array} \right.
\end{eqnarray}
  By symmetry, the corresponding results are also obtained for $P(L_{1} = 0, L_{2} = n)$ for categories \rmn{3}.
}\end{theorem}

\subsection{The marginal distributions for the arithmetic case}
\label{eqn:marginal distributions not v-a}

Under the arithmetic condition that (v-a) does not hold, we consider the tail asymptotics of the marginal distributions. Basically, the results are the same as in Theorems \thrt{exact marginal 10} and \thr{exact marginal 11} in which Theorems \thrt{exact nu 10n-1}~and~\thrt{exact nu 10n-2} should be replaced by Theorems \thrt{exact nu 10A-1}~and~\thrt{exact nu 10A-2} for the case (B1) and Theorems \thrt{exact nu 10B-1}~and~\thrt{exact nu 10B-2} for the case (B2). Thus, we omit their details.

\section{Application to a network with simultaneous arrivals} 
\label{sect:Application to}
 In this section, we apply the asymptotic results to a queueing network with two nodes numbered as $1$ and $2$. Assume that customers simultaneously arrive at both nodes from the outside subject to the Poisson process with rate $\lambda$. For $i=1,2$, service times at node $i$ are independent and identically distributed with the exponential distribution with mean $\mu_{i}^{-1}$. Customers who have finished their services at node $1$ go to node $2$ with probability $p$. Similarly, departing customers from queue $2$ go to queue $1$ with probability $q$. These routing is independent of everything else. Customers what are not routed to the other queue leave the network. We refer to this queueing model as a two node Jackson network with simultaneously arrival.
 
  Obviously, this network is stable, that is, has the stationary distribution, if and only if
\begin{eqnarray}
\label{eqn:stable}
\frac{\lambda(1 + q)}{1-pq} < \mu_{1}, \qquad \frac{\lambda(1 + p)}{1-pq} < \mu_{2}.
\end{eqnarray}
  This fact also can be checked by the stability condition (iv).

  We are interested in how the tail asymptotics of the stationary distribution of this network are changed. If $p=q=0$, this model is studied in \cite{Flatto-Hahn 1984,Flatto-McKean 1977}. As we will see below, this model can be described by a double QBD process, and therefore we know solutions to the tail asymptotics problem. However, this does not mean that the solutions are analytically tractable. Thus, we will consider what kind of difficulty arises in applications of our tail asymptotic results.
  
Let $L_{i}(t)$ be the number of customers at node $i$ at time $t$. It is easy to see that $\{(L_{1}(t), L_{2}(t)); t \in \mathbb{R}_{+}\}$ is a continuous time Markov chain. Because the transition rates of this Markov chain are uniformly bounded, we can construct a discrete time Markov chain given by uniformization, which has the same stationary distribution. We denote this discrete time Markov chain by $\{\vc{L}_{n} = (L_{1\ell},L_{2\ell});\ell \in \mathbb{Z}_{+}\}$, where it is assumed without loss of generality that
\begin{eqnarray*}
  \lambda + \mu_{1} + \mu_{2} = 1.
\end{eqnarray*}
  Obviously, $\{\vc{L}_{n};\ell \in \mathbb{Z}_{+}\}$ is a double QBD process. We denote a random vector subject to the stationary distribution of this process by $\vc{L} \equiv (L_{1}, L_{2})$ as we did in \sectn{Double QBD}.

For applying our asymptotic results, we first compute generating functions. For $\vc{u} = (u_{1},u_{2}) \in \mathbb{R}^{2}$,
\begin{eqnarray}
\label{eqn:gf X+}
 &&\hspace{-6ex}\tilde{\gamma}_{+}(\vc{u}) = \lambda u_1 u_2  + \mu _1 pu_1^{-1} u_2  + \mu _2 qu_1 u_2^{ - 1} + \mu_1 (1 - p)u_1^{-1}  + \mu_2 (1 - q)u_2^{-1},  \\ 
\label{eqn:gf X1} 
&&\hspace{-6ex}\tilde{\gamma}_1 (\vc{u}) = \lambda u_1 u_2  + \mu_1 pu_1^{-1} u_2  + \mu_1 (1 - p)u_1^{-1}  + \mu_2,  \\ 
\label{eqn:gf X2}
 &&\hspace{-6ex}\tilde{\gamma}_2 (\vc{u}) = \lambda u_1 u_2  + \mu _2 qu_1 u_2^{-1} + \mu _2 (1 - q)u_2^{-1}  + \mu_1. \qquad
\end{eqnarray}

 We next find the extreme point $\vc{u}^{(1,\edge)} = (u^{(1,\edge)}_{1}, u^{(1,\edge)}_{2})$. This is obtained as the solution of the equations:
\begin{eqnarray*}
  \tilde{\gamma}_{+}(\vc{u})=\tilde{\gamma}_{1}(\vc{u})=1
\end{eqnarray*}
  Applying \eqn{gf X+} and \eqn{gf X1} to the first equation, we have
\begin{eqnarray}
\label{eqn:gamma+ and gamma1}
u_{2} = u_{1}q + (1 -q).
\end{eqnarray}
Substituting \eqn{gamma+ and gamma1} into $\tilde{\gamma}_{1}(\vc{u})=1$, we have
\begin{eqnarray*}
 \lambda u_1^{2} (u_1 q + 1 - q) + \mu_{1} p ( u_1q + 1 - q ) + \mu_{1} ( 1 - p )  + \mu _2 u_{1} = u_{1}.
\end{eqnarray*}
  Assume that $q > 0$. Then, $u_{1}$ has the following solutions.
\begin{eqnarray*}
u_1  = 1, \frac{ - \lambda \pm \sqrt{\lambda^2  + 4\lambda q\mu_1 ( 1 - pq )} }{2\lambda q}.
\end{eqnarray*}
We are only interested in the solution $u_{1} > 1$, which must be $u^{(1,\edge)}_{1}$, that is,
\begin{eqnarray}
 u^{(1,\edge)}_{1} = \frac{ - \lambda + \sqrt{\lambda^2  + 4\lambda q\mu_1 ( 1 - pq )} }{2\lambda q}.
\end{eqnarray}

We next consider the maximal point $\vc{u}^{(1,\max)}$ of $\tilde{\gamma}(\vc{u}) = 1$. This can be obtained to solve the equations:
\begin{eqnarray*}
  \tilde{\gamma}_{+}(\vc{u}) = 1, \quad \frac {d u_{1}} {d u_{2}} = 0.
\end{eqnarray*}
  These equations are equivalent to
\begin{eqnarray}
\label{eqn:simultaneous 1}
&& \lambda u_1  + \mu _1 pu_1^{ - 1}  - \mu _2 qu_1 u_2^{-2}  - \mu_2 ( 1 - q)u_2^{-2}  = 0,\\
\label{eqn:simultaneous 2}
&& \lambda u_1 u_2  + \mu _1 pu_1^{-1} u_2  + \mu _2 qu_1 u_2^{ - 1} + \mu_1 (1 - p)u_1^{-1}  + \mu_2 (1 - q)u_2^{-1}=1.
\end{eqnarray}
  Theoretically we know that these equations have two solutions such that $\vc{u} > \vc{0}$, which must be $\vc{u}^{(1,\min)}$ and $\vc{u}^{(1,\max)}$. We can numerically obtain them, but their analytic expressions are not easy to get. Furthermore, even if they are obtained, they would be analytically intractable.
  
\begin{figure}[h]
	\caption{Effect of the arrival rate: $\lambda$ is changed from $1$ to $1.2$ and $1.5$ (thicker curves) while $\mu_{1} = 5$, $\mu_{2} = 4$, $p=0.25$, $q=0.4$ are unchanged} 
 	\centering
	\includegraphics[height=6cm]{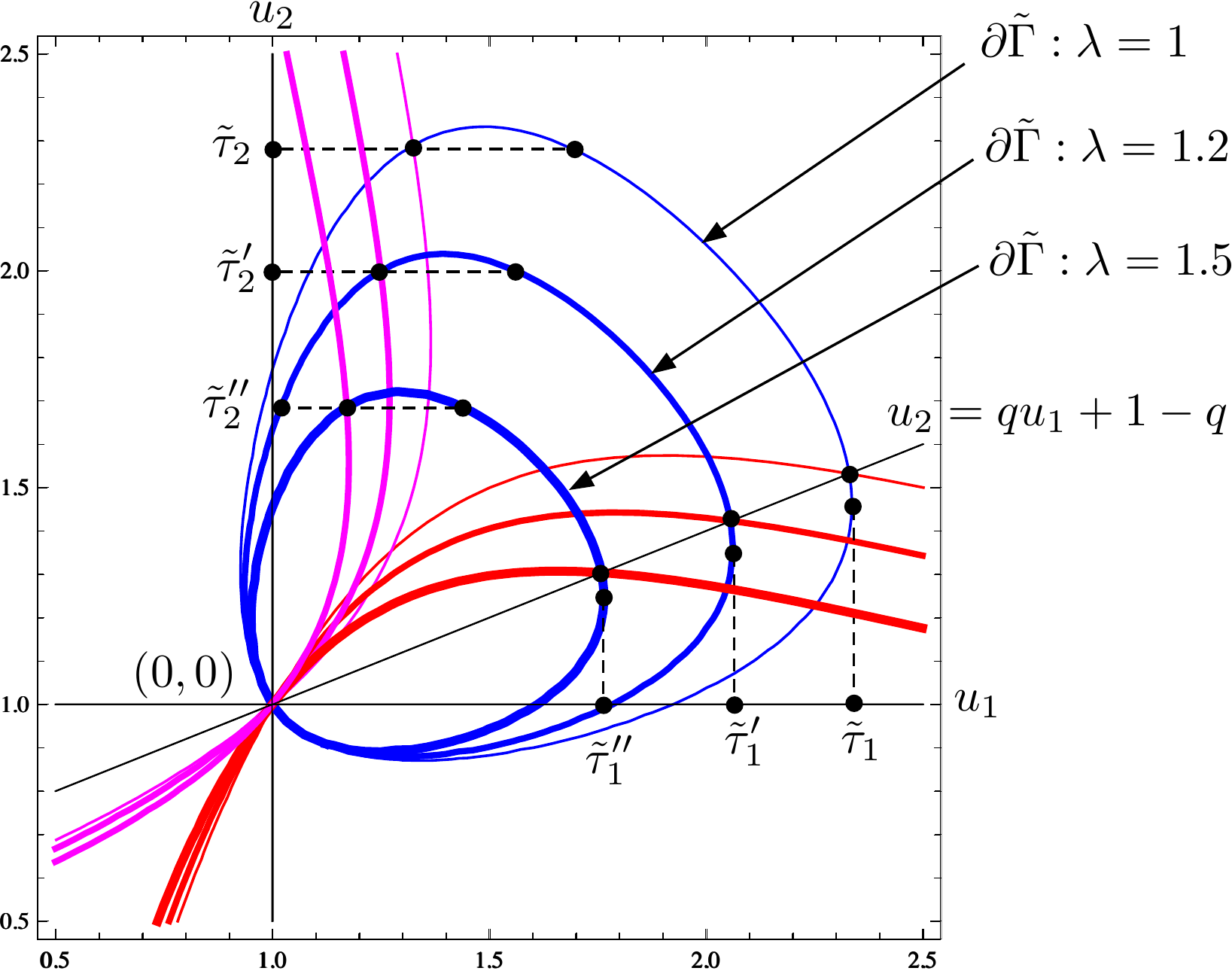} 
	\label{fig:S-arrival 1}
\end{figure}
  
  To circumvent this difficulty, we propose to draw figures. Nowadays we have excellent software such as Mathematica to draw two dimensional figures. Then, we can manipulate figures, and may find how modeling parameters change the tail asymptotics. This is essentially the same as numerical computations. However, figures are more informative to see how changes occur (see, e.g., \fig{S-arrival 1}).
  
  We finally consider a simpler case to find analytically tractable results. Assume that $q=0$ but $p > 0$. $q=0$ implies that
\begin{eqnarray*}
  \vc{u}^{(1,\edge)} = (\rho_{1}^{-1}, 1),
\end{eqnarray*}
  where $\rho_{1} = \frac {\lambda_{1}} {\mu_{1}}$. Obviously, $\rho_{1}$ must the decay rate of $P(L_{1} \ge n)$. This can be also verified by \thr{exact marginal 10}. However, it may not be the decay rate of $P(L_{1} \ge n, L_{2} = 0)$. In fact, we can derive on the curve $\tilde{\gamma}_{+}(\vc{u}) = 1$,
\begin{eqnarray*}
  \left. \frac {d u_{1}} {d u_{2}} \right|_{\vc{u} = \vc{u}^{(1,\edge)}} = \frac {\mu_{2} - (\mu_{1} + \lambda p)} {\lambda (1 - \rho_{1}^{-1})}.
\end{eqnarray*}
  Hence, $u^{(1,\cp)}_{1} = u^{(1,\edge)}_{1}$ if and only if
\begin{eqnarray}
\label{eqn:positive 2}
  \mu_{2} \ge \mu_{1} + \lambda p_{1}.
\end{eqnarray}
  Thus, if \eqn{positive 2} holds, then $P(L_{1} \ge n, L_{2} = 0)$ has the exactly geometric asymptotic. Otherwise, we have, by \thr{exact nu 10n-1},
\begin{eqnarray}
\label{eqn:transient 3}
  \lim_{n \to \infty} n^{-\frac 32} (u^{(1,\max)}_{1})^{-n} P(L_{1} \ge n, L_{2} = 0) = b.
\end{eqnarray}
  We can see that $\rho_{1}^{-1} < u^{(1,\max)}_{1}$, but $u^{(1,\max)}_{1}$ is only numerically obtained by solving \eqn{simultaneous 1} and \eqn{simultaneous 2}.

\section{Concluding remarks} 
\label{sect:Concluding remarks}

We derived the exact asymptotics for the stationary distribution applying the analytic function method based on the convergence domain. We here discuss which problems can be studied by this method and what are needed for further developing it.

\begin{mylist}{3}
\item [({\it Technical issue})] In the analytic function method, a key ingredient is that the function $\ul{\zeta}_{2}(z)$ is analytic and suitably bounded for an appropriate region as we have shown in \lem{extension 1}. For this, we use the fact that $\ul{\zeta}_{2}(z)$ is the solution of a quadratic equation, which is equivalent for the random walk to be skip free in the interior of the quadrant. The quadratic equation (or polynomial equation in general) is also a key for the alternative approach based on analytic extension on Riemann surface. If the random walk is not skip free, it would be harder to get a right analytic function. However, the non skip free case is also interesting. Thus, it is challenging to overcome this difficulty. One here may need a completely different approach.

\smallskip 

\item [({\it Probabilistic interpretation})] We have employed the purely analytic method, and gave no stochastic interpretations except a few although the asymptotic results are stochastic. However, probabilistic interpretations may be helpful. For example, one may wonder what are probabilistic meaning of the function $\ul{\zeta}_{2}$ and the equation \eqn{stationary equation 3}. We do believe something should be here. If they are well answered, then we may better explain \lem{extension 1}, and may resolve the technical issues discussed above.

\smallskip

\item [({\it Modeling extensions})] We think the present approach is applicable for a higher dimensional model as well as a generalized reflecting random walk proposed in \cite{Miyazawa 2011} as long as the skip free assumption is satisfied. One may also consider to relax the irreducibility condition on the random walk in the interior of the quadrant. However, this is essentially equivalent to reducing the dimension, so there should be no difficulty to consider it. Another extension is to modulate the double QBD or multidimensional reflecting random walk in general by a background Markov chain. The tail asymptotic problem becomes harder, but there should be a way to use the present analytic function approach at least for the two dimensional case with finitely many background states. Related discussions can be found in \cite{Miyazawa 2011}.

\smallskip 

\item [({\it Applicability})] As we have seen in \sectn{Application to}, analytic results on the tail asymptotics may not be easy to apply for each specific application because they are not analytically tractable. To fill this gap between theory and application, we have proposed to use geometric interpretations instead of analytic formulas. However, this is currently more or less similar to have numerical tables. We here should make clear what we want to do using the tail asymptotics. Once a problem is set up, we may consider to solve it using geometric interpretations. Probably, there would be a systematic way for this not depending on a specific problem. This is also challenging.
\end{mylist}

\begin{acknowledgement}
We are grateful to Mark S. Squillante for his encouragement to complete this work. We are also thankful to anonymous three referees. This research was supported in part by Japan Society for the Promotion of Science under grant No.\ 21510165. 
\end{acknowledgement}

\appendix

\section*{Appendix}
\addcontentsline{toc}{section}{Appendix}

\subsection*{A. Proof of \lem{D root}}

  Note that $u^{2} D_{2}(u)$ is polynomial with order 2 at least and order 4 at most. For $k=1, 3$, let $c_{k}$ be the coefficients of $u^{k}$ in the polynomial $u^{2} D_{2}(u)$. Then, 
\begin{eqnarray*}
 && c_{1} = - 2(1-p_{00})p_{(-1)0} - 4(p_{(-1)(-1)} p_{01} + p_{(-1)1} p_{0(-1)}) \le 0,\\
 && c_{3} = - 2(1-p_{00})p_{10} - 4(p_{1(-1)} p_{01} + p_{11} p_{0(-1)}) \le 0.
\end{eqnarray*}
  Hence, if both of $u^{(1,\max)}_{1}$ and $-u^{(1,\max)}_{1}$ are the solutions of $u^{2} D_{2}(u) = 0$, then
\begin{eqnarray*}
  2(c_{1} u^{(1,\max)}_{1} + c_{3} (u^{(1,\max)}_{1})^{3}) = (u^{(1,\max)}_{1})^{2} (D_{2}(u^{(1,\max)}_{1}) - D_{2}(-u^{(1,\max)}_{1})) = 0.
\end{eqnarray*}
  Since $u^{(1,\max)}_{1} > 0$, this holds true if and only if $c_{1} = c_{3} = 0$, which is equivalent to that $p_{01}= p_{0(-1)} = p_{(-1)0} = p_{10} = 0$ because $p_{00} = 1$ is impossible. Hence, $u^{2} D_{2}(u) = 0$ has the two solutions $u^{(1,\max)}_{1}$ and $-u^{(1,\max)}_{1}$ if and only if (v-a) does not hold. In this case, we have $c_{1} = c_{3} = 0$, which implies that $u^{2} D_{2}(u)$ is an even function. Since $u^{2} D_{2}(u) = 0$ has only real solutions including $u^{(1,\max)}_{1}$ by Lemmas \lemt{D null points} and \lemt{real root}, we complete the proof. \pend

\subsection*{B. Proof of \lem{arithmetic condition 1}}

By \eqn{arithmetic condition 1}, we have
\begin{eqnarray*}
  \sum_{i \in \{-1,0,1\}} \sum_{j \in \{-1,0,1\}} p_{ij} x^{i} y^{j} = 1, \qquad \sum_{i \in \{-1,0,1\}} \sum_{j \in \{-1,0,1\}} (-1)^{i+j} p_{ij} x^{i} y^{j} = 1.
\end{eqnarray*}
  Subtracting both sides of these equations, we have
\begin{eqnarray*}
  p_{10} x + p_{01} y + p_{0(-1)} y^{-1} + p_{(-1)0} x^{-1} = 0.
\end{eqnarray*}
  Since $x,y$ are positive, this equation holds true if and only if
\begin{eqnarray*}
  p_{10} = p_{01} = p_{0(-1)} = p_{(-1)0} = 0.
\end{eqnarray*}
  This is the condition that (v-a) does not hold. \pend

\end{document}